%
%
%
%
%
\input amstex
\UseAMSsymbols
\loadbold
\magnification=\magstep1   
\vcorrection{-1cm}
\input pictex
\input xy \xyoption{all}
\input graphicx
\nopagenumbers 
\parindent=8mm
\font\Gross=cmr10 scaled\magstep2   
\font\gross=cmr10 scaled\magstep1   
\font\sc=cmcsc10 
\font\srm=cmr7
\font\rmk=cmr8  \font\itk=cmti8 \font\bfk=cmbx8   \font\ttk=cmtt8
\NoBlackBoxes
\mathsurround=1pt

 \def \Ker{\operatorname{Ker}}
 
 \def \Im{\operatorname{Im}}

 \def \rad{\operatorname{rad}}
 
 \def \soc{\operatorname{soc}} 
 \def \sub{\subseteq}
 
 \def \Hom{\operatorname{Hom}}

 \def \can{\operatorname{can}}

 \def \ind{\operatorname{ind}}

 \def \mod{\operatorname{mod}}
 
 \def \supp{\operatorname{supp}}

 \def \lto#1{\;\mathop{\longrightarrow}\limits^{#1}\;} 
 \def \sto#1{\;\mathop{\to}\limits^{#1}\;}
 \def \lamod{\Lambda\text-\mod}

 \def\arr#1#2{\arrow <2mm> [0.25,0.75] from #1 to #2}
 \def\sq{\plot 0 0  1 0  1 1  0 1  0 0 /}
 \def\smallsq#1{\plot 0 0  0.#1 0  0.#1 0.#1  0 0.#1  0 0 /}
 \def\t#1{\;\text{#1}\;}
 
 \def\T#1{\qquad\text{#1}\qquad}
 \def\E#1{\hrule\medskip\noindent#1\medskip\hrule}
 \def\qed{\phantom{m.} $\!\!\!\!\!\!\!\!$\nolinebreak\hfill\checkmark}
 \headline{\ifnum\pageno=1\hfill %
    \else\ifodd\pageno \Rechts \else \Links \fi  \fi}
    \def\Links{\rmk \the\pageno\hfil Schmidmeier\hfil}
    \def\Rechts{\rmk 
        \hfil Bounded Submodules of Modules\hfil\the\pageno}
{\noindent\rmk [s-rel1e-arxiv, August 13, 2004]}

        \vglue3truecm
\centerline{\Gross Bounded Submodules of Modules}
        \bigskip
\centerline{\gross Markus Schmidmeier}
        \medskip
\centerline{\itk Mathematical Sciences, Florida Atlantic University,
        Boca Raton, Florida 33431-0991}
        \smallskip
\centerline{{\itk E-mail address:\/} {\ttk markus\@math.fau.edu}}
        \bigskip
\E{{\baselineskip=9pt {\bfk Abstract}
\newline\smallskip
{\srm Let $\ssize m$, $\ssize n$ be positive integers such that $\ssize m\leq n$.  We consider
all pairs $\ssize (B,A)$ where 
$\ssize B$ is a finite dimensional $\ssize T^n$-bounded $\ssize k[T]$-module and 
$\ssize A$ is a submodule of $\ssize B$ which is $\ssize T^m$-bounded.
They form the objects of the submodule category $\ssize\Cal S_m(k[T]/T^n)$
which is a Krull-Schmidt category with Auslander-Reiten sequences.  
The case $\ssize m=n$ deals with submodules of $\ssize k[T]/T^n$-modules and has been
studied well.
In this manuscript we determine the representation
type of the categories $\ssize \Cal S_m(k[T]/T^n)$ also for the 
cases where $\ssize m<n$: 
It turns out that there are only finitely many indecomposables 
in $\ssize \Cal S_m(k[T]/T^n)$ if either $\ssize m<3$, $\ssize n<6$, or $\ssize (m,n)=(3,6)$;
the category is tame if $\ssize (m,n)$ is one of the pairs $\ssize (3,7)$, $\ssize (4,6)$, $\ssize (5,6)$,
or $\ssize (6,6)$; otherwise, $\ssize \Cal S_m(k[T]/T^n)$ has wild representation type.
Moreover, in each of the finite or tame cases we describe the 
indecomposables and picture the Auslander-Reiten quiver.
}
}}

\bigskip
For a ring $\Lambda$, the {\it submodule category\/} 
$\Cal S(\Lambda)$ has as objects all pairs $(A\sub B)$
where $B$ is a finite length $\Lambda$-module and $A$ a submodule
of $B$.  The morphisms from $(A\sub B)$ to $(A'\sub B')$
are given by the $\Lambda$-linear maps $f\:B\to B'$ which map 
$A$ into $A'$. 
For $m$ a natural number, the 
category $\Cal S_m(\Lambda)$ of {\it bounded\/} 
submodules is the full subcategory of $\Cal S(\Lambda)$
of all pairs $(A\sub B)$ for which $\rad^m\!\!A=0$ holds.

\smallskip Categories of type $\Cal S_m(\Lambda)$ are exact Krull-Schmidt
categories, so every object has a decomposition as a direct sum
of indecomposables; this decomposition is unique up to isomorphism and
reordering.  Aim is the the classification
of the indecomposable objects in the case where $\Lambda=k[T]/T^n$.
The following theorem describes the complexity of this classification problem.
It complements known results for the situation where $m=n$.

\smallskip\noindent {\bf Theorem 1.} {\it
Let $m\leq n$ be positive integers and $k$ be a field.

\smallskip
\item{\rm 1.} The category $\Cal S_m(k[T]/T^n)$ is representation finite
if $n< 6$, or if $m<3$, or if $(m,n)=(3,6)$.

\smallskip
\item{\rm 2.} The category $\Cal S_m(k[T]/T^n)$ has tame infinite 
representation type in case $(m,n)$ is one of the
pairs $(3,7)$, $(4,6)$, $(5,6)$, or $(6,6)$. 

\smallskip
\item{\rm 3.} 
Otherwise the category $\Cal S_m(k[T]/T^n)$ is of wild representation type.

}

\medskip
First we relate our results to Birkhoff's problem on the classification
of subgroup embeddings and list further related results. 
We observe that the objects in $\Cal S_m(k[T]/T^n)$ are just 
the possible configurations 
of a vectorspace together with a subspace that 
is invariant under the action of a nilpotent endomorphism. 
Next we present our methods: Homological algebra,
in particular Auslander-Reiten theory, and covering theory. Finally,
we investigate the submodule 
category $\Cal S_m(k[T]/T^n)$ in each of the following
cases: (1)~$m=1$; (2)~$m=2$; (3)~$m=n-1$; (4)~$m=3$, $n=5$; 
(5)~$m=3$, $n=6$, which is the largest among the representation finite cases;
(6)~$m=3$, $n=7$, where we comment on Birkhoff's family of indecomposables; 
(7)~$m=3$, $n=8$ with a discussion of the wild cases; 
and (8)~$m=4$, $n=6$ with a discussion of $\Cal S_4(k[T]/T^6)$ as a
subcategory of $\Cal S(k[T]/T^6)$. 

\bigskip\centerline{\sc A Comment on Birkhoff's Problem}

\medskip
There is a corresponding problem where $\Lambda=\Bbb Z/p^n$, so the
objects in the category $\Cal S(\Bbb Z/p^n)$ are the pairs $(A\sub B)$
where $B$ is a finite abelian $p^n$-bounded group and $A$ a subgroup
of $B$.  These subgroup embeddings have attracted a lot of attention since 
Birkhoff's work in 1934.  He observes in [B, Corollary~15.1] that the 
number of isomorphism types of embeddings $(A\sub B)$ 
where $B=\Bbb Z/p^6\oplus \Bbb Z/p^4\oplus \Bbb Z/p^2$ tends to infinity
with $p$.
In fact, if $B$ is generated by three elements $x$, $y$, $z$ of 
order $p^6$, $p^4$, $p^2$, and if $A_\lambda$ is the subgroup
of $B$ generated by $u=p^2x+py+z$ and $v=p^2y+\lambda pz$ for some
parameter $\lambda=0,\ldots,p-1$, as pictured below, then the pairs
$(A_\lambda\sub B)_\lambda$ form a family of 
indecomposable and pairwise nonisomorphic
objects in $\Cal S_4(\Bbb Z/p^6)$.  
This is the first such family in the sense that $n=6$ is minimal.
However, $m=4$ is not minimal; in Corollary~15 we show that there is a family
$(C_\lambda\sub D)$ indexed by $\lambda\in k$
of indecomposable and pairwise nonisomorphic objects in the category
in $\Cal S_3(k[T]/T^7)$, and hence a corresponding family
in $\Cal S_3(\Bbb Z/p^7)$, given by the following diagram.
$$
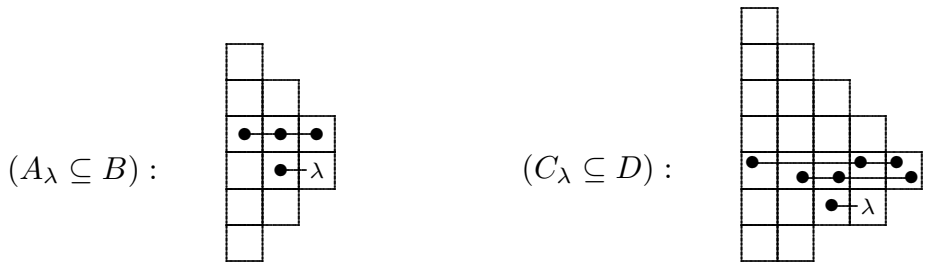
\vbox{\beginpicture 
\setcoordinatesystem units <.4cm,.4cm>
\put{$(A_\lambda\sub B):$} at -4 2
\multiput{\sq} at 0 5  0 4  0 3  0 2  0 1  0 0  1 4  1 3  1 2  1 1  2 3  2 2 /
\put{$\bullet$} at 0.5 3
\put{$\bullet$} at 1.5 3 
\put{$\bullet$} at 2.5 3 
\put{$\bullet$} at 1.5 2 
\put{$\ssize \lambda$} at 2.5 2 
\plot 0.5 3  2.5 3 /
\plot 1.5 2  2.2 2 /
\endpicture}
\qquad\qquad\qquad
\vbox{\beginpicture 
\setcoordinatesystem units <.4cm,.4cm>
\put{$(C_\lambda\sub D):$} at -4 2
\multiput{\sq} at 0 6  0 5  0 4  0 3  0 2  0 1  0 0  1 5  1 4  1 3  1 2  1 1  1 0  
                  2 4  2 3  2 2  2 1  3 3  3 2  3 1  4 2 /
\put{$\bullet$} at 0.3 2.2
\put{$\bullet$} at 3.3 2.2
\put{$\bullet$} at 4.3 2.2
\put{$\bullet$} at 1.7 1.8
\put{$\bullet$} at 2.7 1.8 
\put{$\bullet$} at 4.7 1.8 
\put{$\bullet$} at 2.5 1 
\put{$\ssize \lambda$} at 3.5 1 
\plot 0.3 2.2  4.3 2.2 /
\plot 1.7 1.8  4.7 1.8 /
\plot 2.5 1  3.2 1 /
\endpicture}
$$
Thus, $D$ has five generators $x_1,\ldots,x_5$ of order
$p^7$, $p^6$, $p^4$, $p^3$, and $p$, and the subgroup $C_\lambda$ is
generated by $u_1=p^4x_1+px_4+x_5$, $u_2=p^3x_2+p^2x_3+x_5$, and
$u_3=p^3x_3+p^2\lambda x_4$.

\smallskip
Let us mention some further results.
For the representation finite cases, the classification of the indecomposables
has been completed in [4] in case $m=n\leq 5$ and in [9] in case $(m,n)$ is
one of the pairs  $(3,5)$, $(3,6)$, or $(2,n)$ for some $n$. 
A detailed analysis of the category $\Cal S(k[T]/T^6)$ is given in [6].
The categories $\Cal S_4(\Bbb Z/p^8)$ and $\Cal S_5(\Bbb Z/p^{10})$ are
presented as examples of infinite and of wild type in [1, Examples 8.2.5
and 8.2.6].

\bigskip\centerline{\sc An Equivalent Problem: Invariant Subspaces of a Vector
Space}

\medskip We consider all triples $(V; U, f)$ where
$V$ is a finite dimensional vector space 
over a field $k$, 
$f:V\to V$ a $k$-linear map
and $U$ a subspace of $V$ which is invariant under the action of $f$. 
The morphisms from $(V; U, f)$ to $(V'; U', f')$ are given by
$k$-linear maps $h:V\to V'$ which preserve the subspace, 
that is,  $h(U)\sub U'$ holds,
and which commute with the action of $f$:  $hf'=fh$.
$$
\hbox{\beginpicture
\setcoordinatesystem units <0.7cm,0.5cm>
\put{} at 0 -1
\put{} at 0 3
\put{$\ssize V$} at 0 0
\put{$\ssize U$} at 0 2
\circulararc 160 degrees from 0 3.2 center at 0 2.7
\circulararc -150 degrees from 0 3.2 center at 0 2.7
\circulararc 150 degrees from 0 -1.2 center at 0 -.7
\circulararc -160 degrees from 0 -1.2 center at 0 -.7
\arr{0 1.6}{0 0.4}
\put{$\ssize f|_U$} at -1.1 2.6
\put{$\ssize f$} at -.7 -.6
\put{$\ssize \text{incl}$} at -0.7 1
\arr{0.3 2.45} {0.2 2.3}  
\arr{0.3 -.45} {0.2 -.3}
\put{$\ssize {}\;V'$} at 3 0
\put{$\ssize U'$} at 3 2
\circulararc 160 degrees from 3 3.2 center at 3 2.7
\circulararc -150 degrees from 3 3.2 center at 3 2.7
\circulararc 150 degrees from 3 -1.2 center at 3 -.7
\circulararc -160 degrees from 3 -1.2 center at 3 -.7
\arr{3 1.6}{3 0.4}
\put{$\ssize f'|_{U'}$} at 4.2 2.6
\put{$\ssize f'$} at 3.9 -.6
\put{$\ssize \text{incl}$} at 3.8 1 
\arr{3.3 2.45} {3.2 2.3}  
\arr{3.3 -.45} {3.2 -.3}
\arr{1 0} {2 0}
\arr{1 2} {2 2}
\put{$\ssize h$} at 1.5 -0.4
\put{$\ssize h|_{U,U'}$} at 1.5 2.5
\endpicture}
$$
Assuming that the map $f$ acts nilpotently, then the triple
$(V; U, f)$ has two isomorphism invariants, namely the nilpotency indices $m$
and $n$ of the action of $f$ on $U$ and on $V$, respectively. 
Denote by $\Cal S(m',n')$ the category of all
triples for which the nilpotency indices satisfy $m\leq m'$ and $n\leq n'$.

\medskip\noindent{\bf Lemma 2.} {\it The categories $\Cal S(m,n)$ 
and $\Cal S_m(k[T]/T^n)$ are equivalent. }

\smallskip\noindent{\it Proof:} The action of the variable $T$ in
$k[T]/T^n$ is given by $f$. \qed

        \bigskip\centerline
{\sc Some Homological Properties}
        \medskip

Let $\Lambda$ be an artin algebra, or more general, a locally bounded 
associative $k$-algebra over a local commutative artinian ring $k$.  
Then $\Lambda$ may not have a unit element
but it is assumed that there is a complete set $\{e_i:i\in I\}$ 
of pairwise orthogonal primitive idempotents such that each of the 
indecomposable projective modules $e_i\Lambda$ and $\Lambda e_i$ 
has finite length (see for example [3]).  

\smallskip
Many properties of the category $\Cal S_m(\Lambda)$ of bounded
submodules are obtained from a correponding category $\Cal H_m(\Lambda)$
of morphisms between bounded modules, similar to the interplay
between $\Cal S(\Lambda)$ and $\Cal H(\Lambda)$ studied in [7].
The objects in $\Cal H_m(\Lambda)$ are the $\Lambda$-linear
homomorphisms $(f\: A\to B)$ where $B$ is a finite length
$\Lambda$-module and
$A$ a finite length $\Lambda/\rad^m\!\!\Lambda$-module; 
homomorphisms in this category 
of homomorphisms are given by the commutative diagrams.  
Thus $\Cal H_m(\Lambda)$
is equivalent to the category of modules over the triangular matrix ring
$T=\big({\Lambda/\rad^m\!\!\Lambda \atop 0}\;
        {\Lambda/\rad^m\!\!\Lambda \atop \Lambda}\big)$.  
With the exact structure inherited from $T$-mod, 
the category $\Cal S_m(\Lambda)$ of $\rad^m\!\!\Lambda$-bounded 
submodules becomes a 
{\it full exact subcategory of a module category;} it is 
{\it closed under subobjects and extensions.}
As all objects in $\Cal H_m(\Lambda)$ have finite length, 
it follows that with
$\Cal H_m(\Lambda)$ also $\Cal S_m(\Lambda)$ is a 
{\it Krull-Schmidt category.}

\smallskip
{\it The category $\Cal S_m(\Lambda)$ has
Auslander-Reiten sequences:}  
Let $(f\:A\to B)$ be an object in $\Cal H_m(\Lambda)$, and let
$e'\:\Ker f\to I$ be an injective envelope and $e\:A\to I$ be an
extension of $e'$ to $A$.  Then the two canonical maps
$$(A\sto fB) \;\longrightarrow \;(\Im f\sub B) \T{and}
(A\lto{\big({f \atop e}\big)}B\oplus I)\;\longrightarrow (A\sto fB)$$
modulo the kernel of $f$ and modulo $I$, respectively, 
are a minimal left approximation and a minimal right approximation
of $(f\:A\to B)$ in $\Cal S_m(\Lambda)$.  
The existence of Auslander-Reiten sequences follows from [2, Theorem 2.4].

\smallskip
We describe the indecomposable 
projective and the indecomposable injective
objects in $\Cal S_m(\Lambda)$.

\smallskip\noindent
{\bf Proposition 3.} {\it 
Let $P$ be an indecomposable  projective $\Lambda$-module and $I$
an indecomposable injective $\Lambda$-module.}

\smallskip
\item{1.} {\it The module $(0\sub P)$ is a projective
object in $\Cal S_m(\Lambda)$ with sink map the inclusion map 
$(0\sub\rad P)\to (0\sub P)$.}

\item{2.} {\it  The module $(P/\rad^mP\sub P/\rad^mP)$ is a
projective object in $\Cal S_m(\Lambda)$; the inclusion
map $(\rad P/\rad^mP\sub P/\rad^mP)\to(P/\rad^mP\sub P/\rad^mP)$
is a sink map.}  

\item{3.} {\it The module $(0\sub I)$ is 
a (relatively) injective object in $\Cal S_m(\Lambda)$
with source map the inclusion map
$(0\sub I) \to (\soc I\sub I)$.}

\item{4.} {\it The module $(\soc^mI\sub I)$ is 
injective in $\Cal S_m(\Lambda)$ and a source map is given by the 
canonical map $(\soc^mI\sub I)\to(\soc^mI/\soc I\sub I/\soc I)$.}

\item{5.} {\it Each projective and each injective 
indecomposable 
representation in $\Cal S_m(\Lambda)$ has this form.}

\smallskip\noindent{\it Proof:\/} 
The modules $(0\sub P)$ and $(P/\rad^mP\sub P/\rad^mP)$
are indecomposable projective modules in the
category $\Cal H_m(\Lambda)$, and each indecomposable projective module
in this category has this form.  Also the sink maps are in fact sink
maps in the category $\Cal H_m(\Lambda)$. 

The module $(0\sub I)$ is (relatively) injective in the category
$\Cal S(\Lambda)$, and hence in $\Cal S_m(\Lambda)$, and the source
map is as specified.  

The module $(\soc^mI\sub I)$ is injective in $\Cal H_m(\Lambda)$
and hence in $\Cal S_m(\Lambda)$.  The source map in $\Cal S_m(\Lambda)$
is obtained by taking the composition of the source map 
$(\soc^mI\sub I)\to (\soc^mI\to I/\soc I)$ in $\Cal H_m(\Lambda)$
with the left approximation
$(\soc^mI\to I/\soc I)\to (\soc^mI/\soc I\sub I/\soc I)$ 
in $\Cal S_m(\Lambda)$. 

It remains to verify the last assertion for injective modules.
Let $(A\sub B)$ in $\Cal S_m(\Lambda)$ 
be an indecomposable relatively injective module, and let
$u\:A\to I$ and $w\:B/A\to J$ 
be injective envelopes in the category $\lamod$. 
Since $\rad^mA=0$, the map $u\:A\to I$ restricts to a map $u'\:A\to \soc^mI$.
Also, since $I$ is an injective $\Lambda$-module, $u$ extends to a map
$u''\:B\to I$.  This morphism maps $A$ into $\soc^mI$ and hence gives
rise to a cokernel map $u'''\:B/A\to I/\soc^mI$.
Then $v=\big({u''\atop \can\circ w}\big)\:B\to I\oplus J$
makes the upper part of the following diagram commutative.
$$\CD 0 @>>> A @>>> B @>\can>> B/A @>>> 0 \cr
        @. @Vu'VV @VVvV  @VV{\big({u'''\atop w}\big)}V \cr
      0 @>>> \soc^mI @>{\big({\iota\atop 0}\big)}>> I\oplus J 
         @>{\big({\can\atop 0}{0\atop1}\big)}>> I/\soc^mI \oplus J @>>> 0 \cr
        @. @VVV @VVV @VVV @. \cr
        @. \soc^mI/A @>f>> {(I\oplus J)/\Im v} @>>> X @>>> 0 \endCD $$
The row at the bottom is the cokernel sequence. 
Since $\big({u'''\atop w}\big)$ is a monomorphism, 
so is $f$, by the snake lemma.  Hence the
first two columns of the diagram define a short exact sequence 
in the category $\Cal S_m(\Lambda)$:
$$0\longrightarrow (A\sub B) \longrightarrow (\soc^mI\sub I\oplus J)
\longrightarrow (\soc^mI/A\sub (I\oplus J)/\Im v)\longrightarrow 0$$
By assumption, the object $(A\sub B)$ is relatively injective, and hence
this sequence is split exact. Thus, $(A\sub B)$ is isomorphic to one 
of the indecomposable summands of the middle term, by the Krull-Remak-Schmidt
theorem. This is to say that
$(A\sub B)$ is isomorphic to one of the indecomposable injective
modules in our list. \qed

\smallskip\noindent{\it Example 1.} In case $\Lambda=k[T]/T^n$,
there are the following indecomposable projective and injective modules.
The object $Y=(0\sub\Lambda)$ is projective and relatively injective,
$P=(\Lambda/\rad^m\!\!\Lambda\sub\Lambda/\rad^m\!\!\Lambda)$ is projective, and 
$I=(\rad^{n-m}\Lambda\sub\Lambda)$ is injective in $\Cal S_m(\Lambda)$. 
As usual in this manuscript, given an object $(A\sub B)$ in 
$\Cal S_m(k[T]/T^n)$, we picture each generator of 
the $\Lambda$-module $B$ by a column of boxes
and the image of each generator of $A$ in $B$ by a dot or a dotted line:
$$Y= \;\left.
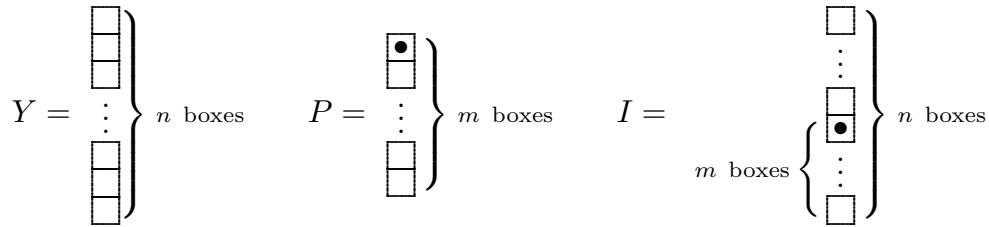
\beginpicture 
        \setcoordinatesystem units <0.3cm,0.3cm>
        \multiput{\sq} at 0 3.75  0 2.75  0 1.75  0 -1.25  0 -2.25  0 -3.25 /
        \put{\vdots} at .5 .5 \endpicture\;
     \right\}{\ssize\text{$n$ boxes}} \qquad
  P= \;\left.\beginpicture 
        \setcoordinatesystem units <0.3cm,0.3cm>
        \multiput{\sq} at 0 2.75  0 1.75  0 -1.25  0 -2.25 /
        \put{$\bullet$} at .5 2.75 
        \put{\vdots} at .5 .5 \endpicture\;
     \right\}{\ssize\text{$m$ boxes}} \qquad
  I= \;\left.\beginpicture 
        \setcoordinatesystem units <0.3cm,0.3cm>
        \multiput{\sq} at 0 3.75  0 .75  0 -.25  0 -3.25 /
        \put{$\bullet$} at .5 -.25 
        \multiput{\vdots} at .5 2.5  .5 -1.5 /  
        \put{$\Bigg\{$} at -.7 -1.75 
        \put{$\ssize\text{$m$ boxes}$} at -3.2 -1.75 \endpicture\;
     \right\}{\ssize\text{$n$ boxes}} $$
Clearly, the modules $P$ and $I$ are isomorphic in case $m=n$.

\smallskip
An indecomposable object is called {\it stable} 
if there is no projective or injective module in its 
orbit under the Auslander-Reiten translation. 
The following result provides an overview for our findings 
in this manuscript.

\medskip\noindent{\bf Theorem 4.} {\it The table below specifies
the stable part of the Auslander-Reiten 
quiver for each of the categories
$\Cal S_m(k[T]/T^n)$ which have finite or tame representation type.

\medskip
\centerline{
\vbox{\offinterlineskip
\halign{\strut#&#\hfil&\vrule#& \;\hfil#\hfil&\;\hfil#\hfil&\;\hfil#\hfil&
          \;\hfil#\hfil&\;\hfil#\hfil&\;\hfil#\hfil&\;\hfil#\hfil\cr
& $m=6$ &height14pt& {} & {} & {} & {} 
        & {\beginpicture\setcoordinatesystem units <0.5cm,0.5cm>\put{$\ssize E_8$} at 0 0
                 \put{} at -0.5 0 \circulararc 300 degrees from 0.4 0 center at -.05 0.05
                  \endpicture} 
        & & \cr
& $\quad5$ &height14pt& {} & {} & {} & ${\ssize \Bbb Z\Bbb E_8/\tau^6}$ 
        & {\beginpicture\setcoordinatesystem units <0.5cm,0.5cm>
                \put{$\ssize E_8$} at 0 0
                \put{} at -0.5 0 \circulararc 
                        300 degrees from 0.4 0 center at -.05 0.05
            \endpicture} 
        & & \cr
& $\quad4$ &height14pt& {} & {} & ${\ssize \Bbb Z\Bbb E_6/\tau^3\sigma}$ 
        & ${\ssize \Bbb Z\Bbb E_7/\tau^6}$ 
        & ${\ssize \Bbb Z\tilde{\Bbb E}_8\vee \Cal T_{5,3,2}}$ & & \cr
& $\quad3$ &height14pt& {} & ${\ssize \Bbb Z\Bbb D_4/\tau^2\rho}$ 
        & ${\ssize \Bbb Z\Bbb D_4/\tau^3\sigma}$
        & ${\ssize \Bbb Z\Bbb D_6/\tau^5}$ 
        & ${\ssize \Bbb Z\Bbb E_8/\tau^{10}}$ 
        & ${\ssize \Bbb Z\tilde{\Bbb E}_8\vee \Cal T_{5,3,2}}$ & \cr
& $\quad2$ &height14pt& ${\ssize \Bbb Z\Bbb A_2/\tau^{\frac 32}\sigma}$ 
        & ${\ssize \Bbb Z\Bbb A_1/\tau^2}$
        & ${\ssize \Bbb Z\Bbb A_4/\tau^{\frac 52}\sigma}$ 
        & ${\ssize \Bbb Z\Bbb A_3/\tau^3\sigma}$ 
        & ${\ssize \Bbb Z\Bbb A_6/\tau^{\frac 72}\sigma}$ 
        & ${\ssize \Bbb Z\Bbb A_5/\tau^4\sigma}$ 
        & ${\ssize \Bbb Z\Bbb A_8/\tau^{\frac 92}\sigma}$ \cr
& $\quad1$ &height14pt& ${\ssize \emptyset}$ & ${\ssize \emptyset}$ & ${\ssize \emptyset}$ 
        & ${\ssize \emptyset}$ 
        & ${\ssize \emptyset}$ & ${\ssize \emptyset}$ & ${\ssize \emptyset}$ \cr
\noalign{\hrule}
& \omit &height14pt& $n=2$ & $3$ & $4$ & $5$ & $6$ & $7$ & $8$ \cr
}}}

\medskip\noindent
For $D$ a Dynkin or a Euclidean diagram, $\tau$ denotes the translation
$\tau:\Bbb Z D\to \Bbb ZD, (\ell, d)\mapsto (\ell-1,d)$;
the map $\sigma$ is an automorphism of order two on $\Bbb Z\Bbb D$ 
which is induced by an automorphism of order two on the underlying 
Dynkin diagram $D=\Bbb A_n$, $\Bbb D_4$, or $\Bbb E_6$; similarly,
$\rho$ is the map on $\Bbb Z\Bbb D_4$ induced from the rotation on $\Bbb D_4$.
By $\Cal T_{a,b,c}$ we denote a family of tubular components of the 
Auslander-Reiten quiver; the family is indexed by $\Bbb P_1(k)$, 
the set of simple $k[T]$-modules, up to isomorphism; all tubes have circumference
one, with the exception of three tubes of circumference $a$, $b$, and $c$.
The symbol \beginpicture\setcoordinatesystem units <0.5cm,0.5cm>
                \put{$\ssize E_8$} at 0 .5
                \put{} at -0.5 0.5 \circulararc 
                        300 degrees from 0.4 0.5 center at -.05 0.55
            \endpicture\,
indicates a tubular category of type $\Bbb E_8$.}

\bigskip\centerline{\sc Galois Coverings}

\medskip
We use covering theory to obtain information 
about the indecomposable objects and about the shape of the
Auslander-Reiten quiver, in particular 
for the representation finite and tame categories of type 
$\Cal S_m(k[T]/T^n)$.

\smallskip
The category $\Cal S_m(k[T]/T^n)$ is just the category of representations
of the quiver
$$\hbox{\beginpicture
\setcoordinatesystem units <0.5cm,0.5cm>
\put{} at -4 -1
\put{} at 1 3.2
\put{$\circ$} at 0 0
\put{$\circ$} at 0 2
\circulararc 160 degrees from 0 3.2 center at 0 2.7
\circulararc -150 degrees from 0 3.2 center at 0 2.7
\circulararc 150 degrees from 0 -1.2 center at 0 -.7
\circulararc -160 degrees from 0 -1.2 center at 0 -.7
\arr{0 1.6}{0 0.4}
\put{$\ssize \alpha$} at -1 2.6
\put{$\ssize \beta$} at -1 -.6
\put{$\ssize \iota$} at 0.4 1
\arr{0.3 2.31} {0.2 2.2}  
\arr{0.3 -.31} {0.2 -.2}
\put{$Q\:$} at -4 1
\endpicture}
$$
which satisfy the relations $\alpha\iota=\iota\beta$, $\alpha^m=0$,
$\beta^n=0$, and the extra condition that the map $\iota$ is a 
monomorphism.
The Galois covering of $Q$ is the quiver $\widetilde Q$.
$$
\hbox{\beginpicture
\setcoordinatesystem units <0.5cm,0.5cm>
\put{} at 0 0
\put{} at 5 2
\put{$\circ$} at 0 0
\put{$\circ$} at 2 0
\put{$\circ$} at 4 0
\put{$\circ$} at 0 2
\put{$\circ$} at 2 2
\put{$\circ$} at 4 2
\arr{-0.4 0}{-1.6 0}
\arr{-0.4 2}{-1.6 2}
\arr{1.6 0}{0.4 0}
\arr{1.6 2}{0.4 2}
\arr{3.6 0}{2.4 0}
\arr{3.6 2}{2.4 2}
\arr{5.6 0}{4.4 0}
\arr{5.6 2}{4.4 2}

\arr{0 1.6}{0 0.4}
\arr{2 1.6}{2 0.4}
\arr{4 1.6}{4 0.4}
\put{$\cdots$} at -3 0
\put{$\cdots$} at -3 2

\put{$\cdots$} at 7 0
\put{$\cdots$} at 7 2
\put{$\ssize \alpha_0$} at -1 2.4 
\put{$\ssize \alpha_1$} at 1 2.4
\put{$\ssize \alpha_2$} at 3 2.4
\put{$\ssize \alpha_3$} at 5 2.4
\put{$\ssize \beta_0$} at -1 -.4 
\put{$\ssize \beta_1$} at 1 -.4
\put{$\ssize \beta_2$} at 3 -.4
\put{$\ssize \beta_3$} at 5 -.4
\put{$\ssize \iota_0$} at 0.5 1 
\put{$\ssize \iota_1$} at 2.5 1
\put{$\ssize \iota_2$} at 4.5 1
\put{$\ssize 0'$} at  0 2.5
\put{$\ssize 1'$} at  2 2.5
\put{$\ssize 2'$} at  4 2.5

\put{$\ssize 0$} at 0 -.5
\put{$\ssize 1$} at 2 -.5
\put{$\ssize 2$} at 4 -.5
\put{$\widetilde Q\:$} at -6 1
\endpicture}$$
Let $\Bbb A_\infty^\infty$ be the double infinite linear quiver, it
has vertex set $\Bbb Z$ and for each $i\in\Bbb Z$, there is an arrow
$\alpha_i\:(i-1)\gets i$. By $(\alpha^n)$ we denote the ideal generated 
by all compositions of $n$ arrows.  The associative $k$-algebra 
$k\Bbb A_\infty^\infty/\alpha^n$ is locally bounded, more precisely it 
is serial in the sense that each 
indecomposable module has a unique composition series,
which has length at most $n$.
Then the category $\Cal S_m(k\Bbb A_\infty^\infty/\alpha^n)$
of $\alpha^m$-bounded submodules of $k\Bbb A_\infty^\infty/\alpha^n$-modules
consists of those representations of $\widetilde Q$ which satisfy
the commutativity relations $\alpha_i\iota_{i-1}=\iota_i\beta_i$,
the nilpotency relations $\alpha^m=0$ and $\beta^n=0$ for all
compositions of $m$ and $n$ arrows $\alpha_i$ and $\beta_j$, respectively,
and the extra condition that all the maps $\iota_i$ are monomorphisms.

\smallskip
There are two functors defined on
$\Cal S_m(k\Bbb A_\infty^\infty/\alpha^n)$.  
The {\it shift\/}
maps an object $M=\Big((M'_i)_i\lto{(\iota_i)_i}(M_i)_i\Big)$ to 
$M[1]=\Big((M'_{i-1})_i\lto{(\iota_{i-1})_i}(M_{i-1})_i\Big)$ 
and is a selfequivalence
on $\Cal S_m(k\Bbb A_\infty^\infty/\alpha^n)$. The 
{\it covering functor\/} 
$F\:\Cal S_m(k\Bbb A_\infty^\infty/\alpha^n)\to \Cal S_m(k[T]/T^n)$
assigns to $M$
the module $F_M=\Big(\bigoplus_iM'_i\sto\iota\bigoplus_iM_i\Big)$ where $\iota$
is the diagonal map given by the $\iota_i$, and the action of $T$
on $\bigoplus_iM'_i$ and on $\bigoplus_iM_i$ is given by the maps for 
the arrows $\alpha_i$ and $\beta_i$, respectively. 
Obviously, $F$ is invariant under the shift:  $F_{M[1]}\cong F_M$.

\smallskip\noindent{\it Example 2.}
We picture the projective and the injective indecomposable representations
in the category $\Cal S_m(k\Bbb A_\infty^\infty/\alpha^n)$.
It follows from Proposition~3 that there are 
three orbits under the shift of projective or injective indecomposables:
The modules $P[i]$, pictured below by their dimension type,
are indecomposable projective and not relatively injective, the $I[i]$
are injective, and the $Y[i]$ projective and relatively injective.
The index $i$ is chosen such that the leftmost 
(and in the case of $P[i]$ the rightmost) nonzero entry in the bottom row 
is at position $i$ in $\widetilde Q$. 
$$P[i]=\;{\ssize{\cdots\atop\cdots}{0\atop0}{0\atop0}\!
        \underbrace{\ssize{1\atop1}{\cdots\atop\cdots}{1\atop1}}_{%
                \text{\!\!\!$m$ ones\!\!\!}}\!
        {0\atop0}{0\atop0}{\cdots\atop\cdots} } \quad
I[i]=\;{\ssize {\cdots\atop\cdots}{0\atop0}{0\atop0}\!
        \underbrace{\ssize \overbrace{\ssize{1\atop1}{\cdots\atop\cdots}{1\atop1}}^{%
                                \text{\!\!\!$m$ ones\!\!\!}}
                \!{0\atop1}{\cdots\atop\cdots}{0\atop1}}_{%
                        \text{$n$ ones}}\!
                        {0\atop0}{0\atop0}{\cdots\atop\cdots} }\quad
Y[i]=\;{\ssize{\cdots\atop\cdots}{0\atop0}{0\atop0}\!
        \underbrace{\ssize{0\atop1}{0\atop1}{\cdots\atop\cdots}{0\atop1}{0\atop1}}_{%
                \text{$n$ ones}}\!{0\atop0}{0\atop0}{\cdots\atop\cdots}} $$
Observe that the covering functor preserves projective and injective 
objects:  With the notation from Example~1 we have
$F_{P[i]}=P$, $F_{I[i]}=I$ and $F_{Y[i]}=Y$.

\smallskip\noindent{\it Definition.}
The {\it support} of 
$M\in\Cal S_m(k\Bbb A_\infty^\infty/\alpha^n)$ 
is the full subcategory $\supp M$
of $\widetilde Q$ given by all points $x$ with $M_x\neq0$. We say that 
$\Cal S_m(k\Bbb A_\infty^\infty/\alpha^n)$ is
{\it locally support finite\/} if for all points $x$
the union is a finite set.
$$\bigcup\Big\{\supp M:M\in \ind \Cal S_m(k\Bbb A_\infty^\infty/\alpha^n)
\t{such that}M_x\neq 0\Big\}$$

We recall the following result from [6, \S 2]:

\medskip\noindent{\bf Theorem 5.} {\it
\smallskip\item{\rm 1.} 
The covering functor $F$ induces a one-to-one 
correspondence between the orbits
under the shift of the isoclasses of indecomposable objects in
$\Cal S_m(k\Bbb A_\infty^\infty/\alpha^n)$ and the isoclasses of 
indecomposable objects in $\Cal S_m(k[T]/T^n)$.
\smallskip\item{\rm 2.} The covering functor preserves Auslander-Reiten 
sequences and induces a one-to-one map from the factor of the 
Auslander-Reiten quiver for $\Cal S_m(k\Bbb A_\infty^\infty/\alpha^n)$
modulo the shift into the Auslander-Reiten quiver for
$\Cal S_m(k[T]/T^n)$. 
\smallskip\item{\rm 3.} If the category 
        $\Cal S_m(k\Bbb A_\infty^\infty/\alpha^n)$  is locally support finite
        then the one-to-one maps in $1.$ and $2.$ are bijective. 
\smallskip\item{\rm 4.} For $M,N\in\Cal S_m(k\Bbb A_\infty^\infty/\alpha^n)$
        the following formula holds.
        $$\qquad\Hom_{\Cal S_m(k[T]/T^n)}(FM,FN)\;=\;\bigoplus_{i\in\Bbb Z}
        \Hom_{\Cal S_m(k\Bbb A_\infty^\infty/\alpha^n)}(M[i],N)\qquad\text{\qed} $$

}

\smallskip\noindent{\it Notation.}
Given $m<n$, let $\Cal S_m(n)$ be an abbreviation for 
$\Cal S_m(k\Bbb A_\infty^\infty/\alpha^n)$. In order to compute the 
indecomposables in $\Cal S_m(n)$ and to determine the Auslander-Reiten
quiver for $\Cal S_m(n)$, we will usually consider some full subquiver
$Q^+$ of $\widetilde Q$, for example 
$$
\hbox{\beginpicture
\setcoordinatesystem units <0.5cm,0.5cm>
\put{} at 0 0
\put{} at 5 2
\put{$Q^+\:$} at -2 1
\put{$\bullet$} at 0 0
\put{$\bullet$} at 0 2
\put{$\bullet$} at 2 0
\put{$\bullet$} at 4 0
\put{$\bullet$} at 2 2
\put{$\bullet$} at 4 2

\arr{1.6 0}{0.4 0}

\arr{1.6 2}{0.4 2}
\arr{3.6 0}{2.4 0}
\arr{5.6 0}{4.4 0}
\arr{2 1.6}{2 0.4}
\arr{0 1.6}{0 0.4}
\arr{3.6 2} {2.4 2}
\arr{5.6 2} {4.4 2}
\arr{4 1.6} {4 0.4}

\setdots<2pt>
\plot 0.5 .5  1.5 1.5 /
\plot 2.5 .5  3.5 1.5 /
\plot 4.5 .5  5.5 1.5 /

\multiput{$\cdots$} at 7 0  7 2  /

\put{$\ssize 0'$} at  0 2.5
\put{$\ssize 1'$} at  2 2.5
\put{$\ssize 2'$} at  4 2.5

\put{$\ssize 0$} at 0 -.5
\put{$\ssize 1$} at 2 -.5
\put{$\ssize 2$} at 4 -.5
\endpicture}
$$
and form the factor algebra $B^+$ of the path algebra $kQ^+$
modulo the ideal given by the commutativity relations 
$\alpha_i\iota_{i-1}=\iota_i\beta_i$ and the nilpotency relations
$\alpha^n=0$ and $\beta ^m=0$. 
By $S_i, P_i, I_i$ and $S'_i, P'_i, I'_i$ we denote the simple, projective,
and injective indecomposable $B^+$-modules
corresponding to the points $i$ and $i'$,
respectively.  
Depending on the index $i$, the objects in $k\widetilde Q$-mod given by
$P_i'$ and $P[i]$ may or may not coincide, and similarly for 
$P_i$ and $Y[i]$ and for $I_i$ and $I[i]$. 
In order to detect which $B^+$-modules correspond
to objects in $\Cal S_m(n)$, the following easy result will be useful.

\medskip\noindent{\bf Lemma 6.} {\it 
The following assertions are equivalent
for $M\in\mod B^+$.}
\smallskip\item{1.} $M\in \Cal S_m(n)$
\item{2.} {\it All maps $M(\iota_i)$ are monomorphisms.}
\item{3.} {\it $\Hom_{B^+}(S_i',M)=0$ for all $i\in \Bbb Z$.}\qed

\medskip
We can now consider the various choices for $m$ and $n$. 

\bigskip\centerline{\sc The Case $m=1$.}

\medskip
For each of the categories
$\Cal S_1(n)$ and $\Cal S_1(k[T]/T^n)$ we determine the indecomposable 
representations and obtain the Auslander-Reiten quiver.

\smallskip
Let $M$ be an indecomposable object in $\Cal S_1(n)$.  We may assume 
that the vector spaces $M_i$ are all zero for $i<0$ and that $M_0\neq 0$; 
otherwise
replace $M$ by a translate under the shift. Thus, $M$ is a representation
of the following subquiver of $Q$ such that $M_0\neq 0$.
$$
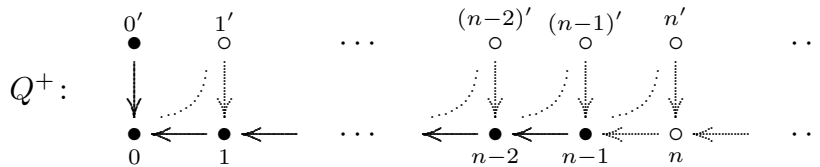
\hbox{\beginpicture
\setcoordinatesystem units <0.5cm,0.5cm>
\put{$Q^+\:$} at 0 1
\put{} at 12 2
\put{$\bullet$} at 2 0
\put{$\bullet$} at 4 0

\put{$\bullet$} at 10 0
\put{$\bullet$} at 12 0

\put{$\bullet$} at 2 2
\multiput{$\circ$} at 4 2  10 2  12 2  14 2  14 0 /
\arr{3.6 0}{2.4 0}
\arr{5.6 0}{4.4 0}
\arr{9.6 0}{8.4 0}
\arr{11.6 0}{10.4 0}
\arr{2 1.6}{2 0.4}

\multiput{$\cdots$} at 7 0  7 2  17 0  17 2 /

\put{$\ssize 0'$} at  2 2.5
\put{$\ssize 1'$} at  4 2.5
\put{$\ssize\strut (n-2)'$} at  10 2.6
\put{$\ssize\strut (n-1)'$} at  12 2.5
\put{$\ssize\strut n'$} at  14 2.6
\put{$\ssize 0$} at 2 -.5
\put{$\ssize 1$} at 4 -.5
\put{$\ssize n-2$} at 10 -.5
\put{$\ssize n-1$} at 12 -.5
\put{$\ssize n$} at 14 -.5
\setdots<1pt>
\arr{4 1.6}{4 0.4}
\arr{10 1.6}{10 0.4}
\arr{12 1.6}{12 0.4}
\arr{14 1.6}{14 0.4}
\arr{13.6 0}{12.4 0}
\arr{15.6 0}{14.4 0}
\setdots<2pt>
\setquadratic
\plot 3.6 1.4  3.3 .7  2.6 .4 /
\plot 9.6 1.4  9.3 .7  8.6 .4 /
\plot 11.6 1.4  11.3 .7  10.6 .4 /
\plot 13.6 1.4  13.3 .7  12.6 .4 /
\endpicture}
$$
In  $\Cal S_1(n)$ the arrows in the upper row in $\widetilde Q$
represent the zero map, so the commutativity relations degenerate
to zero relations as indicated.
Moreover, any composition of $n$ horizontal arrows
is zero.  Let $B^+$ be the factor algebra of $kQ^+$ modulo these relations. 
We are going to show that $M$ has support in the solid part of 
$Q^+$, which forms a diagram of Dynkin type $\Bbb A_n$.  
Consider the left hand part of the Auslander-Reiten quiver
for $B^+$.  As usual in this manuscript, we display dimension vectors in such
a way that the two leftmost entries represent the positions $0$ and $0'$.

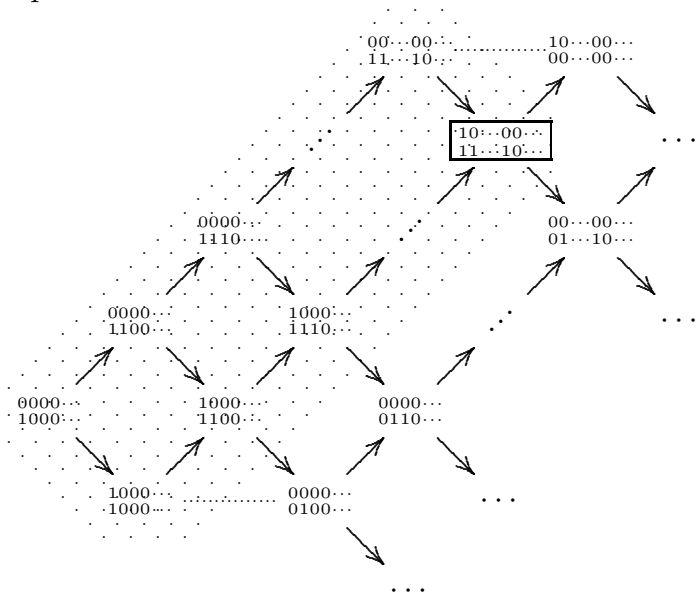
\def\udots{\beginpicture\setcoordinatesystem units <1mm,1mm>
        \multiput{.} at -1 -1  0 0  1 1 / \endpicture}
$$\beginpicture
\setcoordinatesystem units <1cm,1cm>
\put{$\sssize {0000\cdots \atop 1000\cdots}$} at 0 2
\put{$\sssize {1000\cdots \atop 1000\cdots}$} at 1 1
\put{$\sssize {0000\cdots \atop 1100\cdots}$} at 1 3
\put{$\sssize {1000\cdots \atop 1100\cdots}$} at 2 2
\put{$\sssize {0000\cdots \atop 1110\cdots}$} at 2 4
\put{$\sssize {0000\cdots \atop 0100\cdots}$} at 3 1
\put{$\sssize {1000\cdots \atop 1110\cdots}$} at 3 3
\put{$\sssize {0000\cdots \atop 0110\cdots}$} at 4 2
\multiput{$\udots$} at 3 5  4 4  5 3 /
\multiput{$\cdots$} at 4 0  5 1  7 3  7 5 /
\put{$\sssize {00\cdots00\cdots \atop 11\cdots10\cdots}$} at 4 6
\put{\frame{$\ssize\strut{10\cdots00\cdots \atop 11\cdots10\cdots}$}} at 5 5
\put{$\sssize {00\cdots00\cdots \atop 01\cdots10\cdots}$} at 6 4
\put{$\sssize {10\cdots00\cdots \atop 00\cdots00\cdots}$} at 6 6
\arr{0.3 1.7} {0.7 1.3}
\arr{1.3 2.7} {1.7 2.3}
\arr{2.3 3.7} {2.7 3.3}
\arr{4.3 5.7} {4.7 5.3}
\arr{2.3 1.7} {2.7 1.3}
\arr{3.3 2.7} {3.7 2.3}
\arr{5.3 4.7} {5.7 4.3}
\arr{3.3 0.7} {3.7 0.3}
\arr{4.3 1.7} {4.7 1.3}
\arr{6.3 3.7} {6.7 3.3}
\arr{6.3 5.7} {6.7 5.3}
\arr{0.3 2.3} {0.7 2.7}
\arr{1.3 3.3} {1.7 3.7}
\arr{2.3 4.3} {2.7 4.7}
\arr{3.3 5.3} {3.7 5.7}

\arr{1.3 1.3} {1.7 1.7}
\arr{2.3 2.3} {2.7 2.7}
\arr{3.3 3.3} {3.7 3.7}
\arr{4.3 4.3} {4.7 4.7}
\arr{5.3 5.3} {5.7 5.7}
\arr{3.3 1.3} {3.7 1.7}
\arr{4.3 2.3} {4.7 2.7}
\arr{5.3 3.3} {5.7 3.7}

\arr{6.3 4.3} {6.7 4.7}
\setdots<2pt>
\plot 1.5 1  2.5 1 /
\plot 4.5 6  5.5 6 / 
\setshadegrid span <1.5mm>
\vshade -.5   1.5 2.5  <,z,,> 
         .5    .5 3.5  <z,z,,>
        1.5    .5 4.5  <z,z,,>
        3.5   2.5 6.5  <z,z,,>
        4.5   3.5 6.5  <z,z,,>
        5.6   4.6 5.4  /
\endpicture$$
The module in the box is $I[0]$.  
Since $M_0\neq0$, there is a nonzero morphism $M\to I[0]$, and hence
$M$ is one of the $2n$ indecomposable and pairwise nonisomorphic 
modules in the picture which have a path 
to $I[0]$.  Note that the module $S_0'$ occurs as 
an irreducible successor of $I[0]$, so no module of type $S_i'$ has a 
nonzero map to $M$ and hence 
each of the $2n$ modules is in $\Cal S_1(n)$, by Lemma~6. 
The $2n$ predecessors of $I[0]$ form the fundamental domain $\Cal D$ for
the shift in $\Cal S_1(n)$, as indicated by the shaded region:
Each indecomposable object $M$ in $\Cal S_1(n)$ occurs
in exactly one of the sets $\Cal D[i]$.  Here the index is given by
$i=\min\{j\:M_j\neq 0\}$.  

\smallskip
In order to determine the Auslander-Reiten structure for $\Cal S_1(n)$,
we compute the source map for each object $M$ in $\Cal D$.
We have already seen that the source map for $I[0]=(\soc I_0\sub I_0)$ 
is the map $(\soc I_0\sub I_0)\to (0\sub I_0/\soc I_0)$.  
We claim that otherwise, the source map $M\to N$ in $B^+$-mod
is also the source map in the category $\Cal S_1(n)$:  
Clearly, $N$ is an object in $\Cal S_1(n)$, and the factorization property
for a source map holds:
Take an indecomposable 
module $T$ which occurs in $\Cal D[i]$, say, and a nonisomorphism $t\:M\to T$. 
If $i\geq 0$
then $\Hom_{\Cal S}(M,T)=\Hom_{B^+}(M,T)$ and $t$ factors through $N$;
if $i<0$ then there is no 
nonzero map in $\Hom_{\Cal S}(M,T)$ since 
$\Hom_{B^+}(M[-i],T[-i])=0$.  Thus, in any case, $t$ factors through $N$.

\smallskip
The modules on the upper diagonal in $\Cal D$ coincide
with those in the diagonal underneath $\Cal D$, up to the shift and with
the exception of the last one, $Y[0]$.  However, the module $\rad Y[1]$ 
does occur in the diagonal underneath $\Cal D$ and 
we can restore the missing irreducible morphism as
$\rad Y[1]\to Y[1]$, obtaining
the upper diagonal in the domain $\Cal D[1]$.  
Attaching all the fundamental domains $\Cal D[i]$ to each other we obtain
the Auslander-Reiten quiver for $\Cal S_1(n)$; it
is repetitive, up to the shift, as expected.   
For example, if $n=3$ we obtain:
$$\quad
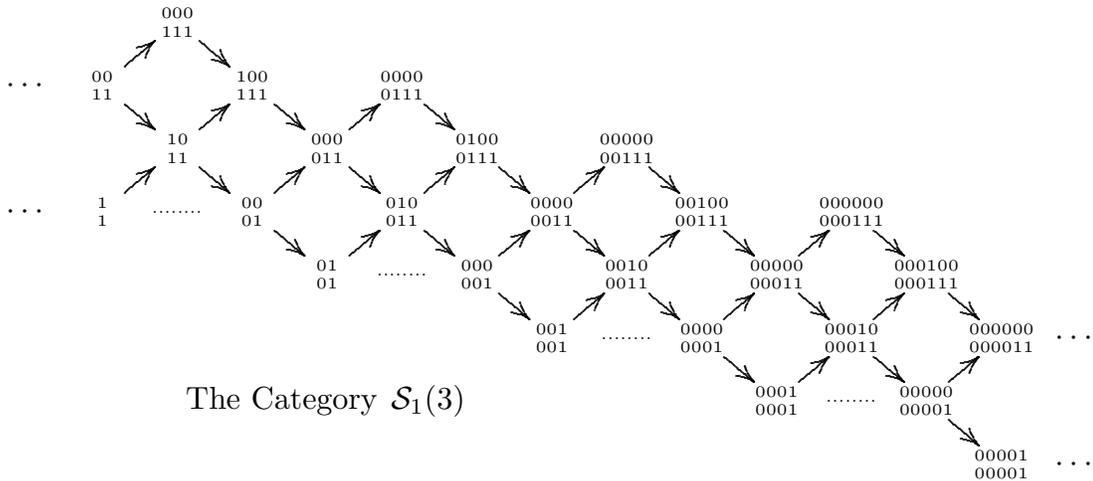
\hbox{\beginpicture
\setcoordinatesystem units <.83cm,0.7cm>
\put{} at 0 0
\put{} at 3 3
\put{The Category $\Cal S_1(3)$} at 3 0
\put{$\ssize {00\atop 11}$}   at 0 5
\put{$\ssize {1\atop 1}$}   at 0 3
\put{$\ssize {10\atop 11}$}   at 1 4
\put{$\ssize {00\atop 01}$} at 2 3
\put{$\ssize {000\atop 111}$} at 1 6
\put{$\ssize {100\atop 111}$} at 2 5
\put{$\ssize {01 \atop 01}$}  at 3 2
\put{$\ssize {000\atop 011}$}    at 3 4
\put{$\ssize {010\atop 011}$}   at 4 3
\put{$\ssize {000\atop 001}$} at 5 2
\put{$\ssize {0000\atop 0111}$} at 4 5
\put{$\ssize {0100\atop 0111}$} at 5 4
\put{$\ssize {001 \atop 001}$}  at 6 1
\put{$\ssize {0000\atop 0011}$}    at 6 3
\put{$\ssize {0010\atop 0011}$}   at 7 2
\put{$\ssize {0000\atop 0001}$} at 8 1
\put{$\ssize {00000\atop 00111}$} at 7 4
\put{$\ssize {00100\atop 00111}$} at 8 3
\put{$\ssize {0001 \atop 0001}$}  at 9 0
\put{$\ssize {00000\atop 00011}$}    at 9 2
\put{$\ssize {00010\atop 00011}$}   at 10 1
\put{$\ssize {00000\atop 00001}$} at 11 0
\put{$\ssize {000000\atop 000111}$} at 10 3
\put{$\ssize {000100\atop 000111}$} at 11 2
\put{$\ssize {00001 \atop 00001}$}  at 12 -1
\put{$\ssize {000000\atop 000011}$}    at 12 1

\arr{0.3 4.7} {0.7 4.3}
\arr{0.3 3.3} {0.7 3.7} 
\arr{1.3 3.7} {1.7 3.3} 
\arr{2.3 3.3} {2.7 3.7} 
\arr{2.3 2.7} {2.7 2.3} 

\arr{0.3 5.3} {0.7  5.7}
\arr{1.3 4.3}  {1.7 4.7}
\arr{1.3 5.7} {1.7  5.3}
\arr{2.3 4.7}  {2.7 4.3}

\arr{3.3 3.7} {3.7 3.3}
\arr{3.3 2.3} {3.7 2.7} 
\arr{4.3 2.7} {4.7 2.3} 
\arr{5.3 2.3} {5.7 2.7} 
\arr{5.3 1.7} {5.7 1.3} 

\arr{3.3 4.3} {3.7  4.7}
\arr{4.3 3.3}  {4.7 3.7}
\arr{4.3 4.7} {4.7  4.3}
\arr{5.3 3.7}  {5.7 3.3}

\arr{6.3 2.7} {6.7 2.3}
\arr{6.3 1.3} {6.7 1.7} 
\arr{7.3 1.7} {7.7 1.3} 
\arr{8.3 1.3} {8.7 1.7} 
\arr{8.3  .7} {8.7  .3} 

\arr{6.3 3.3} {6.7  3.7}
\arr{7.3 2.3}  {7.7 2.7}
\arr{7.3 3.7} {7.7  3.3}
\arr{8.3 2.7}  {8.7 2.3}

\arr{9.3 1.7} {9.7 1.3}
\arr{9.3 0.3} {9.7 .7} 
\arr{10.3 0.7} {10.7 .3} 
\arr{11.3 0.3} {11.7 .7} 
\arr{11.3 -.3} {11.7 -.7} 

\arr{9.3  2.3} {9.7  2.7}
\arr{10.3 1.3}  {10.7 1.7}
\arr{10.3 2.7} {10.7  2.3}
\arr{11.3 1.7}  {11.7 1.3}

\setdots<2pt>
\plot 0.7 3  1.3 3 /
\plot 3.7 2  4.3 2 /
\plot 6.7 1  7.3 1 /
\plot 9.7 0  10.3 0 /

\multiput{$\cdots$} at -1 3  -1 5  13 1  13 -1 /
\endpicture} 
$$
We will usually picture only one copy of the fundamental domain:
$$
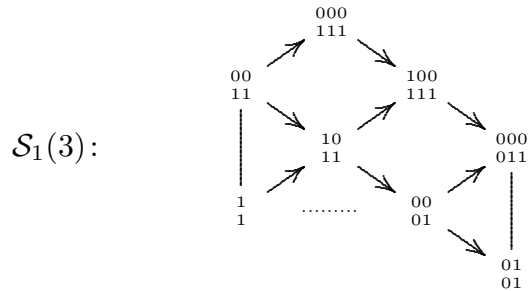
\hbox{\beginpicture
\setcoordinatesystem units <1cm,0.7cm>
\put{$\Cal S_1(3)\:$} at -2 2
\put{} at 3 3
\put{$\ssize {00\atop 11}$}   at 0 3
\put{$\ssize {1\atop 1}$}   at 0 1
\put{$\ssize {10\atop 11}$}   at 1 2
\put{$\ssize {00\atop 01}$} at 2 1
\put{$\ssize {01 \atop 01}$}  at 3 0
\put{$\ssize {000\atop 011}$}    at 3 2

\put{$\ssize {000\atop 111}$} at 1 4
\put{$\ssize {100\atop 111}$} at 2 3

\arr{0.3 2.7} {0.7 2.3}
\arr{0.3 1.3} {0.7 1.7} 
\arr{1.3 1.7} {1.7 1.3} 
\arr{2.3 1.3} {2.7 1.7} 
\arr{2.3 0.7} {2.7 0.3} 

\arr{0.3 3.3} {0.7 3.7}
\arr{1.3 2.3}  {1.7 2.7}
\arr{1.3 3.7} {1.7 3.3}
\arr{2.3 2.7}  {2.7 2.3}

\setdots<2pt>
\plot 0.7 1  1.3 1 /

\setsolid
\plot 0 2.6  0 1.4 /
\plot 3 1.6  3 0.4 /
\endpicture} $$

For the corresponding categories of type $\Cal S_1(k[T]/T^n)$ we use
covering theory to obtain:

\medskip\noindent{\bf Proposition 7.} {\it
Let $m=1$.  The category $\Cal S_m(k[T]/T^n)$ 
has $2n$ indecomposables; they are of the form
$(0\sub k[T]/T^i)$ or $(T^{i-1}k[T]/T^i\sub k[T]/T^i)$ for 
some $1\leq i\leq n$.  The Auslander-Reiten quiver
for  $\Cal S_m(k[T]/T^n)$ has the shape of a bounded tube with one coray
and two rays, or conversely.
The Auslander-Reiten translation shifts objects along a helix; the two
orbits are finite and nonperiodic of length 1 and $2n-1$. 
The examples below are obtained from the Auslander-Reiten quivers
for $\Cal S_1(n)$ by identifying the objects modulo the shift.} \qed
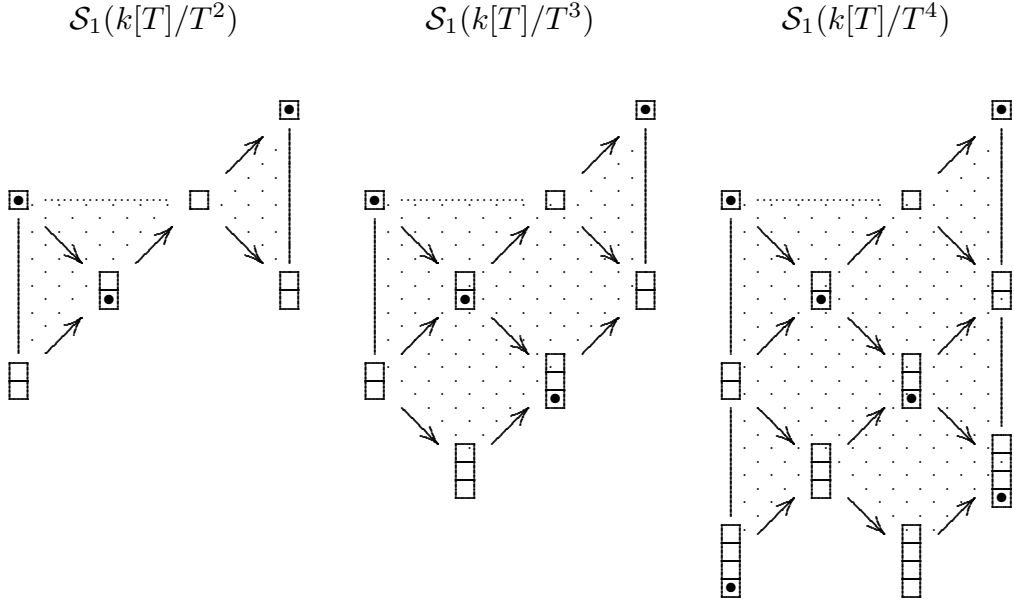
$$
\hbox{\beginpicture
\setcoordinatesystem units <1cm,1cm>
\put{} at 0 0
\put{} at 3 4
\put{$\Cal S_1(k[T]/T^2)$} at 1.5 4
\multiput{$\smallsq2$} at -.1 -.1  -.1 .1 /
\put{$\smallsq2$} at -.1 2
\put{$\ssize \bullet$} at 0 2
\multiput{$\smallsq2$} at .9 .9  .9 1.1 /
\put{$\ssize \bullet$} at 1 .9
\put{$\smallsq2$} at 1.9 2
\put{$\smallsq2$} at 2.9 3 
\put{$\ssize\bullet$}  at 3 3
\multiput{$\smallsq2$} at 2.9 .9 2.9 1.1 /

\arr{0.3 0.3} {0.7 0.7}
\arr{0.3 1.7} {0.7 1.3} 
\arr{1.3 1.3} {1.7 1.7} 
\arr{2.3 1.7} {2.7 1.3} 
\arr{2.3 2.3} {2.7 2.7} 

\setdots<2pt>
\plot 0.3 2 1.7 2 /

\setdashes <2cm>
\plot 0 0.3  0 1.8 /
\plot 3 1.3  3 2.8 /

\setshadegrid span <1.5mm>
\vshade   0   0 2  <,z,,> 
          2   2 2  <z,z,,>
          3   1 3 /
\endpicture} 
\qquad
\hbox{\beginpicture
\setcoordinatesystem units <1cm,1cm>
\put{} at 0 -1
\put{} at 3 4
\put{$\Cal S_1(k[T]/T^3)$} at 1.5 4
\multiput{$\smallsq2$} at -.1 -.1  -.1 .1 /
\put{$\smallsq2$} at -.1 2
\put{$\ssize\bullet$}  at 0 2
\multiput{$\smallsq2$} at .9 .9  .9 1.1 /
\put{$\ssize\bullet$}   at 1 .9
\put{$\smallsq2$} at 1.9 2
\put{$\smallsq2$} at 2.9 3
\put{$\ssize\bullet$}  at 3 3
\multiput{$\smallsq2$}  at 2.9 .9  2.9 1.1 /
\multiput{$\smallsq2$}  at 0.9  -1.2  .9  -1  .9 -.8 /
\multiput{$\smallsq2$} at  1.9 -.2  1.9 0  1.9 .2 /
\put{$\ssize\bullet$} at 2 -.2

\arr{0.3 0.3} {0.7 0.7}
\arr{0.3 1.7} {0.7 1.3} 
\arr{1.3 1.3} {1.7 1.7} 
\arr{2.3 1.7} {2.7 1.3} 
\arr{2.3 2.3} {2.7 2.7} 

\arr{0.3 -0.3} {0.7 -0.7}
\arr{1.3 0.7}  {1.7 0.3}
\arr{1.3 -0.7} {1.7 -0.3}
\arr{2.3 0.3}  {2.7 0.7}

\setdots<2pt>
\plot 0.3 2  1.7 2 /

\setdashes <2cm>
\plot 0 0.3  0 1.8 /
\plot 3 1.3  3 2.8 /
\setshadegrid span <1.5mm>
\vshade   0   0 2  <,z,,> 
          1   -1 2 <z,z,,>
          2   0 2  <z,z,,>
          3   1 3 /
\endpicture} 
\qquad
\hbox{\beginpicture
\setcoordinatesystem units <1cm,1cm>
\put{} at 0 -2
\put{} at 3 4
\put{$\Cal S_1(k[T]/T^4)$} at 1.5 4
\multiput{$\smallsq2$} at -.1 -.1  -.1 .1 /
\put{$\smallsq2$} at -.1 2
\put{$\ssize\bullet$}  at 0 2
\multiput{$\smallsq2$} at .9 .9  .9 1.1 /
\put{$\ssize\bullet$}   at 1 .9
\put{$\smallsq2$} at 1.9 2
\put{$\smallsq2$} at 2.9 3
\put{$\ssize\bullet$}  at 3 3
\multiput{$\smallsq2$}  at 2.9 .9  2.9 1.1 /
\multiput{$\smallsq2$}  at 0.9  -1.2  .9  -1  .9 -.8 /
\multiput{$\smallsq2$} at  1.9 -.2  1.9 0  1.9 .2 /
\put{$\ssize\bullet$} at 2 -.2
\multiput{$\smallsq2$} at -.1 -2.3  -.1 -2.1  -.1 -1.9  -.1 -1.7 /
\put{$\ssize\bullet$} at 0 -2.3
\multiput{$\smallsq2$} at 1.9 -2.3  1.9 -2.1  1.9 -1.9  1.9 -1.7 /
\multiput{$\smallsq2$} at 2.9 -1.3  2.9 -1.1  2.9 -.9  2.9 -.7 /
\put{$\ssize\bullet$} at 3 -1.3

\arr{0.3 0.3} {0.7 0.7}
\arr{0.3 1.7} {0.7 1.3} 
\arr{1.3 1.3} {1.7 1.7} 
\arr{2.3 1.7} {2.7 1.3} 
\arr{2.3 2.3} {2.7 2.7} 

\arr{0.3 -0.3} {0.7 -0.7}
\arr{1.3 0.7}  {1.7 0.3}
\arr{1.3 -0.7} {1.7 -0.3}
\arr{2.3 0.3}  {2.7 0.7}

\arr{0.3 -1.7} {0.7 -1.3}
\arr{1.3 -1.3} {1.7 -1.7}
\arr{2.3 -1.7} {2.7 -1.3}
\arr{2.3 -0.3} {2.7 -0.7}

\setdots<2pt>
\plot 0.3 2  1.7 2 /

\setdashes <2cm>
\plot 0 0.3  0 1.8 /
\plot 0 -1.5 0 -0.3 /
\plot 3 1.3  3 2.8 /
\plot 3 -0.5 3 0.7 /
\setshadegrid span <1.5mm>
\vshade   0   -2 2  <,z,,> 
          1   -1 2 <z,z,,>
          2   -2 2  <z,z,,>
          3   -1 3 /
\endpicture} 
$$

\bigskip
\centerline{\sc The Case $m=2$.}

\medskip We determine the indecomposables in the 
categories $\Cal S_2(n)$ and  $\Cal S_2(k[T]/T^n)$.  
There is a subtle difference in the orbit structure of the 
Auslander-Reiten quivers, depending on whether $n$ is even or odd.

\smallskip
Consider the following quiver
$$
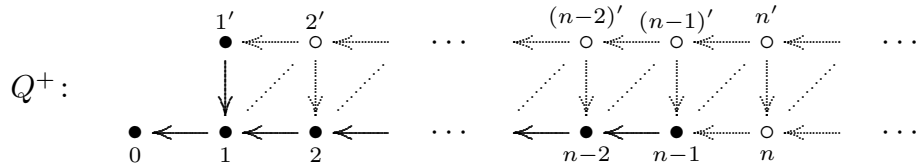
\hbox{\beginpicture
\setcoordinatesystem units <0.5cm,0.5cm>
\put{} at 0 0
\put{} at 5 2
\put{$Q^+\:$} at -2 1
\put{$\bullet$} at 0 0
\put{$\bullet$} at 2 0
\put{$\bullet$} at 4 0

\put{$\bullet$} at 10 0
\put{$\bullet$} at 12 0
\put{$\circ$} at 14 0

\put{$\bullet$} at 2 2
\put{$\circ$} at 4 2
\put{$\circ$} at 10 2
\put{$\circ$} at 12 2
\put{$\circ$} at 14 2

\arr{1.6 0}{0.4 0}
\arr{3.6 0}{2.4 0}
\arr{5.6 0}{4.4 0}

\arr{9.6 0}{8.4 0}
\arr{11.6 0}{10.4 0}

\arr{2 1.6}{2 0.4}

\setdots<1pt>
\arr{13.6 0} {12.4 0}
\arr{15.6 0} {14.4 0}
\arr{3.6 2} {2.4 2}
\arr{5.6 2} {4.4 2}
\arr{9.6 2} {8.4 2}
\arr{11.6 2}{10.4 2}
\arr{13.6 2}{12.4 2}
\arr{15.6 2}{14.4 2}
\arr{4 1.6} {4 0.4}
\arr{10 1.6} {10 0.4}
\arr{12 1.6} {12 0.4}
\arr{14 1.6} {14 0.4}

\setdots<2pt>
\plot 2.5 .5  3.5 1.5 /
\plot 4.5 .5  5.5 1.5 /
\plot 8.5 .5  9.5 1.5 /
\plot 10.5 .5  11.5 1.5 /
\plot 12.5 .5  13.5 1.5 /
\plot 14.5 .5  15.5 1.5 /

\multiput{$\cdots$} at 7 0  7 2  17 0  17 2 /

\put{$\ssize 1'$} at  2 2.5
\put{$\ssize 2'$} at  4 2.5

\put{$\ssize 0$} at 0 -.5
\put{$\ssize 1$} at 2 -.5
\put{$\ssize 2$} at 4 -.5

\put{$\ssize n-2$} at 10 -.5
\put{$\ssize n-1$} at 12 -.5
\put{$\ssize n$} at  14 -.5

\put{$\ssize \strut(n-2)'$} at 10 2.6
\put{$\ssize \strut(n-1)'$} at 12 2.5
\put{$\ssize \strut n'$} at  14 2.6

%
\endpicture}
$$
and the corresponding algebra $B^+$ given by 
the usual commutativity and nilpotence relations. 
We will see that with one exception, 
the indecomposables in $\Cal S_2(n)$, up to the shift,
correspond to the representations of $Q^+$ which have support 
in the Dynkin diagram $\Bbb D_{n+1}$ indicated by the solid dots and lines.  
The correspondence is given by mapping the $B^+$-module $M$
to the representation $M^*$ of $\widetilde Q$, as follows.
$$M\:\quad \xymatrix@1@=5mm{ & \ssize M_1' \ar[d]^{j_1} 
                & \ssize M_2' \ar[l]_{a_2} \ar[d]^{j_2} & \cdots \ar[l] \\
        \ssize M_0 & \ssize M_1 \ar[l]^{b_1} 
                & \ssize M_2 \ar[l]^{b_2} & \cdots \ar[l]}
\qquad M^*\:\quad \xymatrix@1@=5mm{ \ssize M_0\ar[d]^1 
                & \ssize M_1' \ar[d]^{j_1} \ar[l]_{j_1b_1} 
                & \ssize M_2' \ar[l]_{a_2} \ar[d]^{j_2} & \cdots \ar[l] \\
        \ssize M_0 & \ssize M_1 \ar[l]^{b_1} 
                & \ssize M_2 \ar[l]^{b_2} & \cdots \ar[l]}$$
The following Lemma is easy to verify.

\medskip\noindent{\bf Lemma 8.} {\it The following statements
are equivalent for a module $M\in \Cal S_2(n)$ which
is such that $M_i=0$ for $i<0$.}
\smallskip\item{1.}{\it The map $j_0$ is an isomorphism.}
\item{2.}{\it $M\cong N^*$ for some $B^+$-module $N$.}
\item{3.}{\it There is no nonzero map $M\to Y[0]$.}\qed

\smallskip
The left hand part of the Auslander-Reiten quiver for $B^+$
has the following form; an object $M$ is represented by the
dimension vector of the corresponding object $M^*$ in $\Cal S_2(n)$.
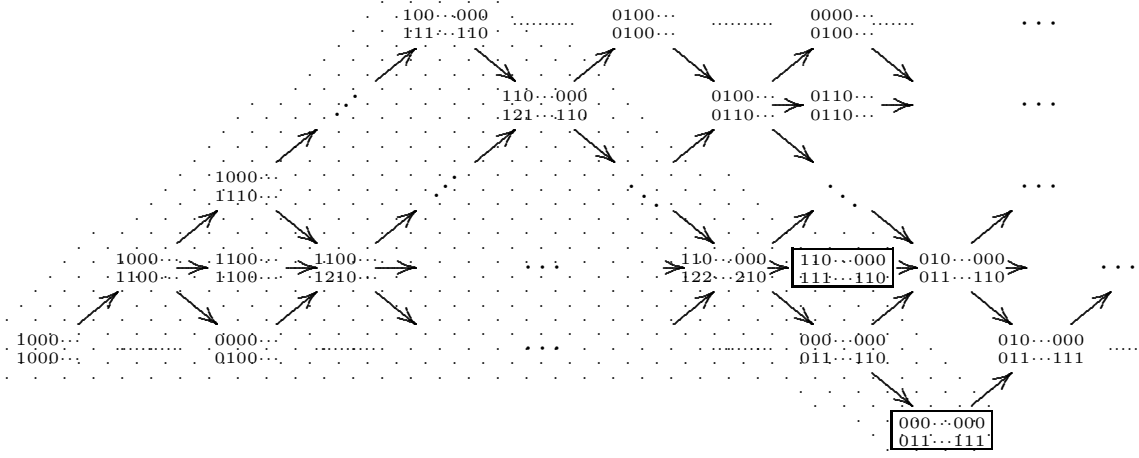
$$
\hbox{\beginpicture
\setcoordinatesystem units <1.1cm,0.9cm>
\put{} at 1 0
\put{$\ssize {1000\cdots\atop1000\cdots}$}   at 1 0
\put{$\ssize {1000\cdots\atop1100\cdots}$}   at 2 1
\put{$\ssize {1100\cdots\atop1100\cdots}$}    at 3 1
\put{$\ssize {0000\cdots\atop0100\cdots}$} at 3 0
\put{$\ssize {1000\cdots\atop1110\cdots}$}   at 3 2
\put{$\ssize {1100\cdots\atop1210\cdots}$} at 4 1

\put{$\ssize {100\cdots000\atop111\cdots110}$}  at 5 4
\put{$\ssize {110\cdots000\atop121\cdots110}$}   at 6 3
\put{$\ssize {0100\cdots\atop0100\cdots}$} at 7 4
\put{$\ssize {0100\cdots\atop0110\cdots}$} at 8 3
\put{$\ssize {0110\cdots\atop0110\cdots}$} at 9 3
\put{$\ssize {0000\cdots\atop0100\cdots}$}    at 9 4
\put{$\ssize {110\cdots000\atop122\cdots210}$} at 7.8 1
\put{$\ssize {000\cdots000\atop011\cdots110}$} at 9 0
\put{\frame{$\ssize\strut{110\cdots000\atop111\cdots110}$}}    at 9 1
\put{\frame{$\ssize\strut{000\cdots000\atop011\cdots111}$}}    at 10 -1
\put{$\ssize {010\cdots000\atop011\cdots110}$}    at 10.2 1
\put{$\ssize {010\cdots000\atop011\cdots111}$}    at 11 0
\arr{9.3 -0.3}{9.7 -0.7}
\arr{10.3 -0.7}{10.7 -0.3}
\arr{1.3 0.3} {1.7 0.7} 
\arr{2.3 0.7} {2.7 0.3} 
\arr{3.3 0.3} {3.7 0.7} 
\arr{4.3 0.7} {4.7 0.3} 
\arr{7.3 0.3} {7.7 0.7} 
\arr{8.3 0.7} {8.7 0.3} 
\arr{9.3 0.3} {9.7 0.7} 
\arr{10.3 0.7} {10.7 0.3}
\arr{11.3 0.3} {11.7 0.7}
\arr{2.3 1.3} {2.7 1.7} 
\arr{3.3 1.7} {3.7 1.3} 
\arr{4.3 1.3} {4.7 1.7} 

\arr{7.3 1.7} {7.7 1.3}
\arr{8.3 1.3} {8.7 1.7} 
\arr{9.3 1.7} {9.7 1.3}
\arr{10.3 1.3} {10.7 1.7}
\arr{2.3 1}{2.6 1}
\arr{3.4 1}{3.7 1}
\arr{4.3 1}{4.7 1}
\arr{7.2 1}{7.4 1}
\arr{8.25 1}{8.45 1}
\arr{9.55 1}{9.75 1}
\arr{10.65 1}{10.85 1}
\arr{3.3 2.3} {3.7 2.7} 

\arr{5.3 2.3} {5.7 2.7}
\arr{6.3 2.7} {6.7 2.3} 
\arr{7.3 2.3} {7.7 2.7} 
\arr{8.3 2.7} {8.7 2.3} 
\arr{8.3 3}{8.6 3}
\arr{9.4 3}{9.7 3}
\arr{4.3 3.3} {4.7 3.7} 
\arr{5.3 3.7} {5.7 3.3} 
\arr{6.3 3.3} {6.7 3.7} 
\arr{7.3 3.7} {7.7 3.3}
\arr{8.3 3.3} {8.7 3.7} 
\arr{9.3 3.7} {9.7 3.3}
\setdots<2pt>
\plot 1.7 0  2.3 0 /
\plot 3.7 0  4.3 0 /
\plot 7.7 0  8.3 0 /
\plot 11.7 0 12 0 /
\plot 5.7 4  6.3 4 /
\plot 7.7 4  8.3 4 /
\plot 9.7 4  9.3 4 /
\multiput{$\udots$} at 4 3  5 2 /
\multiput{$\ddots$} at 7 2  9 2 /
\multiput{$\cdots$} at 11.8 1  11 2  11 3  11 4  6 0  6 1 /
\setshadegrid span <1.6mm>
\vshade   .5  -.5 .5  <,z,,> 
          4.5 -.5 4.5 <z,z,,>
          5.5 -.5 4.5 <z,z,,>
          8.5 -.5 1.5 <z,z,,>
          9.5 -1.5 1.5 <z,z,,>
          9.6 -1.5 0.4 <z,z,,>
         10.5 -1.5 -.5 /
\endpicture} 
$$
The boxed modules are $I[0]$ and $Y[1]$ and the shaded area
will be the fundamental domain $\Cal D$. Note that in case
that $n$ is odd, the positions of the predecessor of $Y[1]$ and
of $I[0]$ have to be exchanged.  This is due to the operation
of Gabriel's Frobenius permutation: If $n+1$ is even, then the
projective modules $P_i$ corresponding to the point $i$ in the Dynkin
diagram $\Bbb D_{n+1}$ occur in the same $\tau$-orbits as the injective module 
$I_i$ for this point; if $n+1$ is odd, the modules $P_1'$ and $I_0=I[0]$
are in the same orbit, and so are $P_0$ and $I_1'$ (which is at the
position of the successor of $Y[1]$ in the Auslander-Reiten 
quiver of the Dynkin diagram).

\smallskip
Next, we determine the indecomposable representations of
$\Cal S_m(n)$, up to the shift.
Let $M$ be an indecomposable representation such that 
$M_i=0=M_{i'}$ for all $i<0$, and such that $M_0\neq 0$.
If $j_0$ is an isomorphism, then $M$ is a $B^+$-module
and has a non-zero map to the injective 
representation $I[0]$. 
If $j_0$ is not an isomorphism, then $M[1]$ is a 
$B^+$-module,
and the factor $N:=(M/\Im(j_i)_i)[1]$ has a non-zero map
to $Y[1]$. Thus, the set $\Cal D$ of predecessors of either $I[0]$
or $Y[1]$ as indicated by the shaded region in the above diagram,
contains a  full set of representatives for the 
indecomposable representations of $\Cal S_m(n)$, up to the shift.
Conversely, $\Cal D$ consists only of objects in $\Cal S_2(n)$ since
no module of dimension type $\ssize{\cdots010\cdots\atop\cdots000\cdots}$
has a nonzero map into one of the indecomposables in $\Cal D$.  
In conclusion, $\Cal D$ consists of those indecomposables $M$ in $\Cal S_2(n)$
which satisfy the following three conditions. (1) $M_i=0$ for $i<0$,
(2) $M(\iota_0)$ is an isomorphism, (3) if $M_0= 0$ then 
        $M(\iota_1)$ is not an isomorphism.  It follows that 
$\Cal D$ is a fundamental domain for the shift and we read off that
$\Cal D$ contains $1+2+\cdots+(n-1)+n\;+(n-2)+2=
n(n+1)/2+n=\frac n2(n+3)$ indecomposable objects.

\smallskip
Note that the modules in the antidiagonal on the right hand side of $\Cal D$
are exactly the shifted copies of the modules on the diagonal on the left
hand side within  $\Cal D$.  Thus, in order to 
obtain the Auslander-Reiten quiver for $\Cal S_2(n)$,
one verifies as in the case where $m=1$ that each of the source maps
of a module in $\Cal D$ in the category $B^+$-mod
is also a source map in the category $\Cal S_2(n)$.

\smallskip The Auslander-Reiten quiver for $\Cal S_2(n)$ is obtained by
``attaching'' the shifted copies $\Cal D[i]$ of the fundamental domain
to each other, as above in the case $m=1$. In the two examples below,
the dashed line separates the modules in $\Cal D$ from those in 
$\Cal D[1]$. 
$$
%
%
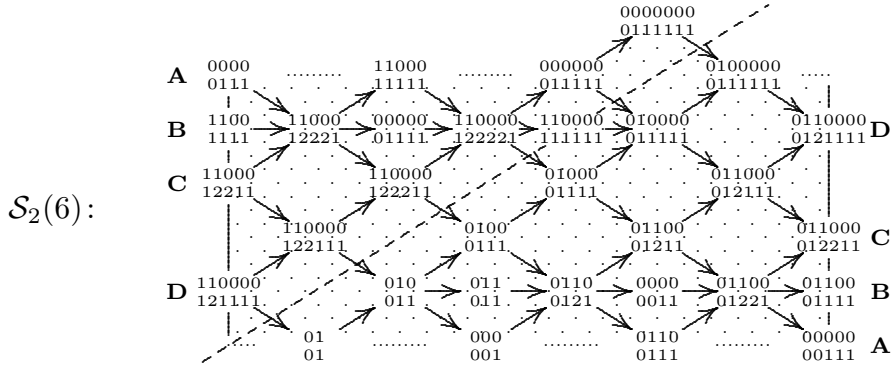
\hbox{\beginpicture
\setcoordinatesystem units <0.95cm,0.6cm>
\put{} at 1 0
\put{} at 8 6
\put{$\Cal S_2(6)\:$} at -1 2.5
\put{$\ssize {01 \atop 01    }$}   at 2 0
\put{$\ssize {000\atop 001}$}      at 4 0
\put{$\ssize {0110\atop 0111  }$}   at 6 0
\put{$\ssize {00000\atop 00111}$} at 8 0

\put{$\ssize {110000 \atop 121111}$}   at 1 1
\put{$\ssize {010\atop 011}$}       at 3 1
\put{$\ssize {011\atop 011}$}   at 4 1
\put{$\ssize {0110\atop 0121}$}       at 5 1
\put{$\ssize {0000 \atop 0011}$}   at 6 1
\put{$\ssize {01100\atop 01221}$}       at 7 1
\put{$\ssize {01100\atop 01111}$} at 8 1

\put{$\ssize {110000\atop 122111}$} at 2 2
\put{$\ssize {0100\atop 0111}$} at 4 2
\put{$\ssize {01100\atop 01211}$} at 6 2
\put{$\ssize {011000\atop 012211}$} at 8 2

\put{$\ssize {11000\atop 12211}$} at 1 3
\put{$\ssize {110000\atop 122211}$} at 3 3
\put{$\ssize {01000\atop 01111}$} at 5 3
\put{$\ssize {011000\atop 012111}$} at 7 3

\put{$\ssize {1100\atop 1111}$} at 1  4
\put{$\ssize {11000\atop 12221}$} at 2 4
\put{$\ssize {00000\atop 01111}$} at 3 4
\put{$\ssize {110000\atop 122221}$} at 4 4
\put{$\ssize {110000\atop 111111}$} at 5 4
\put{$\ssize {010000\atop 011111}$} at 6 4
\put{$\ssize {0110000 \atop 0121111}$}   at 8 4

\put{$\ssize {0000\atop 0111}$} at 1 5
\put{$\ssize {11000 \atop 11111}$} at 3 5
\put{$\ssize {000000 \atop 011111}$} at 5 5
\put{$\ssize {0100000 \atop 0111111}$} at 7 5

\put{$\ssize {0000000 \atop 0111111}$} at 6 6
\arr{1.3 0.7} {1.7 0.3} 
\arr{2.3 0.3} {2.7 0.7} 
\arr{3.3 0.7} {3.7 0.3} 
\arr{4.3 0.3} {4.7 0.7} 
\arr{5.3 0.7} {5.7 0.3} 
\arr{6.3 0.3} {6.7 0.7} 
\arr{7.3 0.7} {7.7 0.3} 
\arr{3.3 1} {3.7 1}
\arr{4.3 1} {4.7 1}
\arr{5.3 1} {5.7 1}
\arr{6.3 1} {6.7 1}
\arr{7.3 1} {7.7 1}
\arr{1.3 1.3} {1.7 1.7} 
\arr{2.3 1.7} {2.7 1.3} 
\arr{3.3 1.3} {3.7 1.7} 
\arr{4.3 1.7} {4.7 1.3} 
\arr{5.3 1.3} {5.7 1.7} 
\arr{6.3 1.7} {6.7 1.3} 
\arr{7.3 1.3} {7.7 1.7} 
\arr{1.3 2.7} {1.7 2.3} 
\arr{2.3 2.3} {2.7 2.7} 
\arr{3.3 2.7} {3.7 2.3} 
\arr{4.3 2.3} {4.7 2.7} 
\arr{5.3 2.7} {5.7 2.3} 
\arr{6.3 2.3} {6.7 2.7} 
\arr{7.3 2.7} {7.7 2.3} 
\arr{1.3 3.3} {1.7 3.7} 
\arr{2.3 3.7} {2.7 3.3} 
\arr{3.3 3.3} {3.7 3.7} 
\arr{4.3 3.7} {4.7 3.3} 
\arr{5.3 3.3} {5.7 3.7} 
\arr{6.3 3.7} {6.7 3.3} 
\arr{7.3 3.3} {7.7 3.7} 
\arr{1.3 4} {1.7 4}
\arr{2.3 4} {2.7 4}
\arr{3.3 4} {3.6 4}
\arr{4.4 4} {4.6 4}
\arr{5.4 4} {5.7 4}
\arr{1.3 4.7} {1.7 4.3} 
\arr{2.3 4.3} {2.7 4.7} 
\arr{3.3 4.7} {3.7 4.3} 
\arr{4.3 4.3} {4.7 4.7} 
\arr{5.3 4.7} {5.7 4.3} 
\arr{6.3 4.3} {6.7 4.7} 
\arr{7.3 4.7} {7.7 4.3} 
\arr{5.3 5.3} {5.7 5.7} 
\arr{6.3 5.7} {6.7 5.3} 
\setdots<2pt>
\plot   1 0  1.3 0 /
\plot 2.7 0  3.3 0 /
\plot 4.7 0  5.3 0 /
\plot 6.7 0  7.3 0 /
\plot 1.7 5  2.3 5 /
\plot 3.7 5  4.3 5 /
\plot 7.7 5  8   5 /
\setsolid
\plot 8 0.4  8 0.6 /
\plot 8 1.4  8 1.6 /
\plot 8 2.4  8 3.6 /
\plot 8 4.4  8 4.8 /
\plot 1 0.2  1 0.6 /
\plot 1 1.4  1 2.6 /
\plot 1 3.4  1 3.6 /
\plot 1 4.4  1 4.6 /
\multiput{$\ssize\text{\bf A}$} at .4  5  8.6 0 /  
\multiput{$\ssize\text{\bf B}$} at .4  4  8.6 1 /
\multiput{$\ssize\text{\bf C}$} at .4  3  8.6 2 /
\multiput{$\ssize\text{\bf D}$} at .4  1  8.6 4 /
\setdashes <1mm>
\plot .7 -.3  4.7 3.7 /
\plot 5.3 4.3  7.3 6.3 /
\setshadegrid span <1.5mm>
\vshade   1  0 5  <,z,,> 
          5  0 5  <z,z,,>
          6  0 6  <z,z,,>
          7  0 5  <z,,,>
          8  0 5 /
\endpicture} 
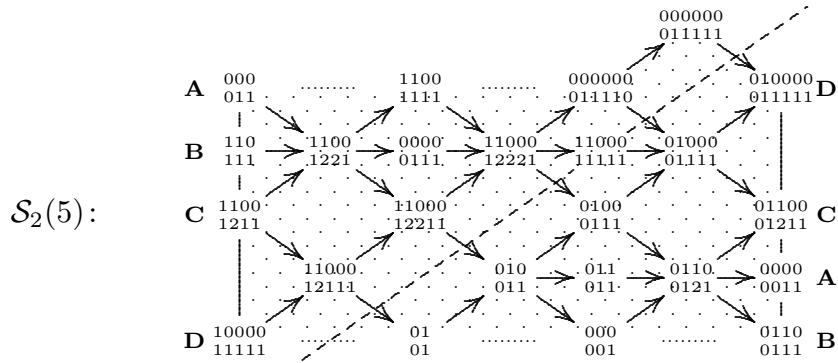
$$

$$
%
%
\hbox{\beginpicture
\setcoordinatesystem units <1cm,0.7cm>
\put{$\Cal S_2(5)\:$} at -2 2
\put{$\ssize {10000 \atop 11111}$} at 0 0
\put{$\ssize {01    \atop 01   }$} at 2 0
\put{$\ssize {000   \atop 001  }$} at 4 0
\put{$\ssize {0110  \atop 0111 }$} at 6 0
\put{$\ssize {11000 \atop 12111}$} at 1 1
\put{$\ssize {010   \atop 011  }$} at 3 1
\put{$\ssize {011   \atop 011  }$} at 4 1
\put{$\ssize {0110  \atop 0121 }$} at 5 1
\put{$\ssize {0000  \atop 0011 }$} at 6 1
\put{$\ssize {1100  \atop 1211 }$} at 0 2
\put{$\ssize {11000 \atop 12211}$} at 2 2
\put{$\ssize {0100  \atop 0111 }$} at 4 2
\put{$\ssize {01100 \atop 01211}$} at 6 2
\put{$\ssize {110   \atop 111  }$} at 0 3
\put{$\ssize {1100  \atop 1221 }$} at 1 3
\put{$\ssize {0000  \atop 0111 }$} at 2 3
\put{$\ssize {11000 \atop 12221}$} at 3 3
\put{$\ssize {11000 \atop 11111}$} at 4 3
\put{$\ssize {01000 \atop 01111}$} at 5 3
\put{$\ssize {000   \atop 011  }$} at 0 4
\put{$\ssize {1100  \atop 1111 }$} at 2 4
\put{$\ssize {000000\atop011110}$} at 4 4
\put{$\ssize {010000\atop011111}$} at 6 4
\put{$\ssize {000000\atop011111}$} at 5 5
\arr{0.3 0.3} {0.7 0.7} 
\arr{1.3 0.7} {1.7 0.3}
\arr{2.3 0.3} {2.7 0.7}
\arr{3.3 0.7} {3.7 0.3}
\arr{4.3 0.3} {4.7 0.7}
\arr{5.3 0.7} {5.7 0.3}

\arr{3.3 1  } {3.7 1  }
\arr{4.3 1  } {4.7 1  }
\arr{5.3 1  } {5.7 1  }

\arr{0.3 1.7} {0.7 1.3}
\arr{1.3 1.3} {1.7 1.7}
\arr{2.3 1.7} {2.7 1.3}
\arr{3.3 1.3} {3.7 1.7}
\arr{4.3 1.7} {4.7 1.3}
\arr{5.3 1.3} {5.7 1.7}

\arr{0.3 2.3} {0.7 2.7}
\arr{1.3 2.7} {1.7 2.3}
\arr{2.3 2.3} {2.7 2.7}
\arr{3.3 2.7} {3.7 2.3}
\arr{4.3 2.3} {4.7 2.7}
\arr{5.3 2.7} {5.7 2.3}

\arr{0.3 3  } {0.7 3  }
\arr{1.3 3  } {1.7 3  }
\arr{2.3 3  } {2.7 3  }
\arr{3.3 3  } {3.7 3  }
\arr{4.3 3  } {4.7 3  }

\arr{0.3 3.7} {0.7 3.3}
\arr{1.3 3.3} {1.7 3.7}
\arr{2.3 3.7} {2.7 3.3}
\arr{3.3 3.3} {3.7 3.7}
\arr{4.3 3.7} {4.7 3.3}
\arr{5.3 3.3} {5.7 3.7}

\arr{4.3 4.3} {4.7 4.7}
\arr{5.3 4.7} {5.7 4.3}
\setdots<2pt>
\plot 0.7 0  1.3 0 /
\plot 2.7 0  3.3 0 /
\plot 4.7 0  5.3 0 /
\plot 0.7 4  1.3 4 /
\plot 2.7 4  3.3 4 /
\setsolid
\plot 0 0.4  0 1.6 /
\plot 0 2.4  0 2.6 /
\plot 0 3.4  0 3.6 /
\plot 6 0.4  6 0.6 /
\plot 6 1.4  6 1.6 /
\plot 6 2.4  6 3.6 /
\multiput{$\ssize\text{\bf A}$} at -.5 4  6.5 1 /  
\multiput{$\ssize\text{\bf B}$} at -.5 3  6.5 0 /
\multiput{$\ssize\text{\bf C}$} at -.5 2  6.5 2 /
\multiput{$\ssize\text{\bf D}$} at -.5 0  6.5 4 /
\setdashes <1mm>
\plot .7 -.3  3.7 2.7 /
\plot 4.3 3.3  6.3 5.3 /
\setshadegrid span <1.5mm>
\vshade   0  0 4  <,z,,> 
          4  0 4  <z,z,,>
          5  0 5  <z,,,>
          6  0 4 /
\endpicture} 
$$

In order to obtain the Auslander-Reiten quiver for $\Cal S_2(k[T]/T^n)$,
identify the modules on the left edge with their shifted copies on the right.
The type of this identification depends on whether $n$ is even or odd,
and is indicated by the letters A--D. We conclude:

\medskip\noindent{\bf Proposition 9.} {\it
For $m=2$ and $n>2$, the category $\Cal S_m(k[T]/T^n)$ 
has $\frac n2(n+3)$ indecomposable representations, 
namely the following:
$$
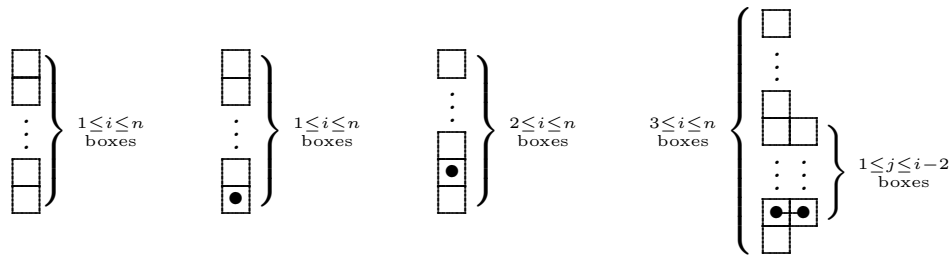
\left.
\hbox{\beginpicture 
        \setcoordinatesystem units <0.3cm,0.3cm>
        \multiput{\sq} at 0 2.75  0 1.75  0 -1.25  0 -2.25 /
        \put{\vdots} at .5 .5 \endpicture}\;
     \right\}{\ssize{1\leq i\leq n \atop\text{boxes}}} \qquad
\left.
\hbox{\beginpicture 
        \setcoordinatesystem units <0.3cm,0.3cm>
        \multiput{\sq} at 0 2.75  0 1.75  0 -1.25  0 -2.25  /
        \put{$\bullet$} at .5 -2.25
        \put{\vdots} at .5 .5 \endpicture}\;
     \right\}{\ssize{1\leq i\leq n\atop \text{boxes}}} \qquad
\left.
\hbox{\beginpicture
        \setcoordinatesystem units <0.3cm,0.3cm>
        \multiput{\sq} at 0 2.75   0 -.25  0 -1.25  0 -2.25 /
        \put{\vdots} at .5 1.5 
        \put{$\bullet$} at .5 -1.25 \endpicture}\;
     \right\}{\ssize{2\leq i\leq n\atop\text{boxes}}} \qquad
{\ssize{3\leq i\leq n\atop\text{boxes}}}\left\{\;
\hbox{\beginpicture
        \setcoordinatesystem units <0.3cm,0.3cm>
        \multiput{\sq} at 0 4.25  1 .25  0 1.25  1 -2.75  0 -3.75
                0 -2.75  0 .25 /
        \put{$\Bigg\}$} at 2.7 -1.25
        \multiput{$\bullet$} at .5 -2.75  1.5 -2.75 /
        \plot .5 -2.75  1.5 -2.75 /
        \multiput{\vdots} at .5 -1  1.5 -1  .5 3 / 
        \put{$\ssize{1\leq j\leq i-2\atop\text{boxes}}$} at 5.2 -1.25 
        \endpicture}\;\right.
     $$
In case $n$ is even, there are stable modules on the boundary of the Auslander-Reiten quiver.
The stable part of this quiver has type $\Bbb Z\Bbb A_n/\tau^{(n+1)/2}\sigma$,
where $\sigma$ is the quiver automorphism of the Dynkin diagram $\Bbb A_n$ of
order two, and hence for each stable representation $M$, 
the formula $M\cong\tau^{n+1}(M)$ holds.
There are two orbits attached to the stable part, one contains
$n-1$ modules between $P$ and $I$; the other one consists only of the 
projective injective representation $Y$.

\smallskip\noindent
In case $n$ is odd, there are no stable modules on the boundary 
of the Auslander-Reiten quiver.
The stable part has type $\Bbb Z \Bbb A_{n-2}/\tau^{(n+1)/2}\sigma$;
each module $M$ on the central axis satisfies $M\cong\tau^{\frac{n+1}2}(M)$,
each stable module  satisfies $M\cong\tau^{n+1}(M)$.
There is a non-stable $\tau$-orbit of length $2n$ 
containing the projective representation $P$
and the injective representation $I$; 
attached to this orbit is the projective injective module $Y$.} \qed

\smallskip
Here are two examples for illustration.

$$
%
%
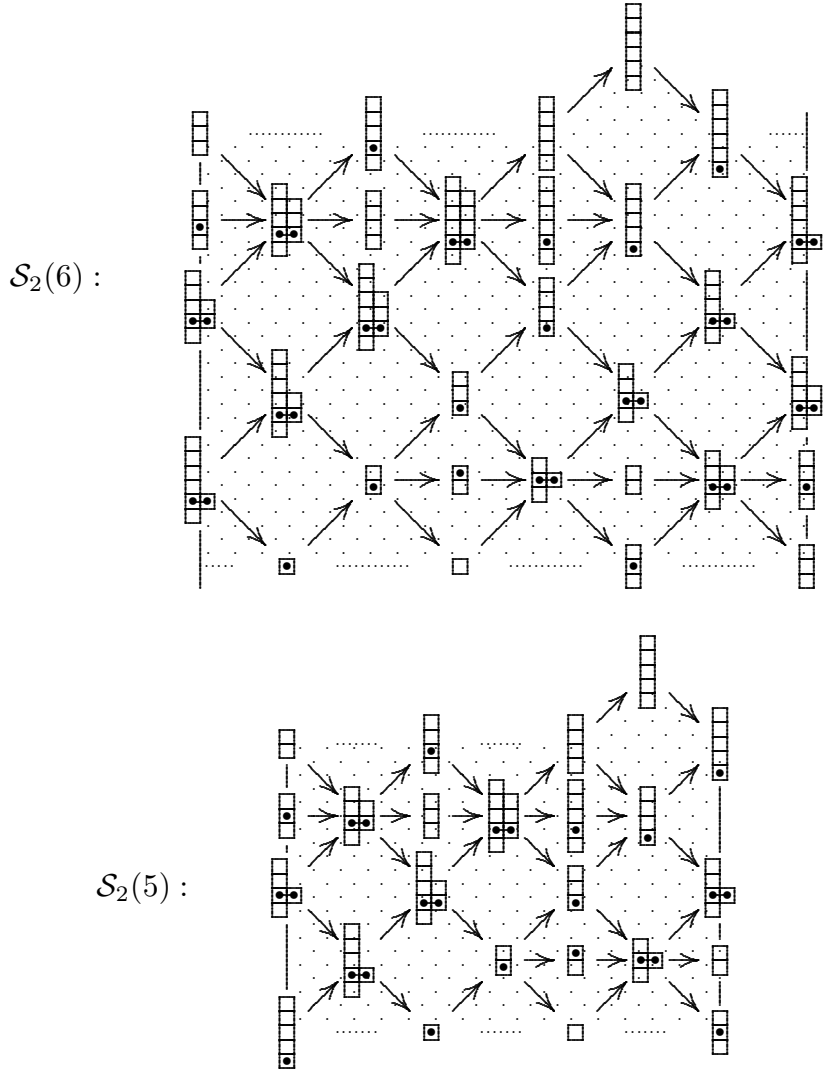
\hbox{\beginpicture
\setcoordinatesystem units <.8cm,.8cm>
\put{} at 1.2 0
\put{} at 9.6 7.2
\put{$\Cal S_2(6):$} at -.8 4

\put{$\smallsq2$} at 2.3 0
\put{$\sssize \bullet$} at 2.4 0
\put{$\smallsq2$} at 4.7 0
\multiput{$\smallsq2$} at 7.1 -.2  7.1 0 7.1 .2 /
\put{$\sssize \bullet$} at 7.2 0
\multiput{$\smallsq2$} at  9.5 -.2  9.5 0  9.5 .2 /

\multiput{$\smallsq2$} at 1 .7  1 .9  1 1.1  1 1.3  1 1.5  1 1.7  1.2 .9 /
\multiput{$\sssize \bullet$} at 1.1 .9  1.3 .9 /
\plot 1.1 .9  1.3 .9 /
\multiput{$\smallsq2$} at 3.5 1.1  3.5 1.3 /
\put{$\sssize \bullet$} at 3.6 1.1
\multiput{$\smallsq2$} at 4.7 1.1  4.7 1.3 /
\put{$\sssize \bullet$} at 4.8 1.3
\multiput{$\smallsq2$} at 5.8 1  5.8 1.2  5.8 1.4  6 1.2 /
\multiput{$\sssize \bullet$} at 5.9 1.2  6.1 1.2 /
\plot 5.9 1.2  6.1 1.2 /
\multiput{$\smallsq2$} at 7.1 1.1  7.1 1.3 /
\multiput{$\smallsq2$} at 8.2 .9  8.2 1.1  8.2 1.3  8.2 1.5  8.4 1.1  8.4 1.3 /
\multiput{$\sssize \bullet$} at 8.3 1.1  8.5 1.1 /
\plot 8.3 1.1  8.5 1.1 /
\multiput{$\smallsq2$} at 9.5 .9  9.5 1.1  9.5 1.3  9.5 1.5 /
\put{$\sssize \bullet$} at  9.6 1.1 

\multiput{$\smallsq2$} at 2.2 1.9  2.2 2.1  2.2 2.3  2.2 2.5  2.2 2.7  2.2 2.9  2.4 2.1  2.4 2.3 /
\multiput{$\sssize \bullet$} at 2.3 2.1  2.5 2.1 /
\plot 2.3 2.1  2.5 2.1 /
\multiput{$\smallsq2$} at 4.7 2.2  4.7 2.4  4.7 2.6 /
\put{$\sssize \bullet$} at 4.8 2.2
\multiput{$\smallsq2$} at 7 2.1  7 2.3  7 2.5  7 2.7  7.2 2.3 /
\multiput{$\sssize \bullet$} at  7.1 2.3  7.3 2.3 /
\plot 7.1 2.3  7.3 2.3 /
\multiput{$\smallsq2$} at 9.4 2  9.4 2.2  9.4 2.4  9.4 2.6  9.4 2.8  9.6 2.2  9.6 2.4 /
\multiput{$\sssize \bullet$} at 9.5 2.2  9.7 2.2 /
\plot 9.5 2.2  9.7 2.2 /

\multiput{$\smallsq2$} at 1 3.2  1 3.4  1 3.6  1 3.8  1 4  1.2 3.4  1.2 3.6 /
\multiput{$\sssize \bullet$} at 1.1 3.4  1.3 3.4 /
\plot 1.1 3.4  1.3 3.4 /
\multiput{$\smallsq2$} at 3.4 3.1  3.4 3.3  3.4 3.5  3.4 3.7  3.4 3.9  3.4 4.1  3.6 3.3  3.6 3.5  3.6 3.7 /
\multiput{$\sssize \bullet$} at 3.5 3.3  3.7 3.3 /
\plot 3.5 3.3  3.7 3.3 /
\multiput{$\smallsq2$} at 5.9 3.3  5.9 3.5  5.9 3.7  5.9 3.9 /
\put{$\sssize \bullet$} at 6 3.3
\multiput{$\smallsq2$} at 8.2 3.2  8.2 3.4  8.2 3.6  8.2 3.8  8.2 4  8.4 3.4 /
\multiput{$\sssize \bullet$} at 8.3 3.4  8.5 3.4 /
\plot 8.3 3.4  8.5 3.4 /

\multiput{$\smallsq2$} at 1.1 4.5  1.1 4.7  1.1 4.9  1.1 5.1 /
\put{$\sssize \bullet$} at 1.2  4.7
\multiput{$\smallsq2$} at 2.2 4.4  2.2 4.6  2.2 4.8  2.2 5.0  2.2 5.2  2.4 4.6  2.4 4.8  2.4 5 /
\multiput{$\sssize \bullet$} at 2.3 4.6  2.5 4.6 /
\plot 2.3 4.6  2.5 4.6 /
\multiput{$\smallsq2$} at 3.5 4.5  3.5 4.7  3.5 4.9  3.5 5.1 /
\multiput{$\smallsq2$} at 4.6 4.3  4.6 4.5  4.6 4.7  4.6 4.9  4.6 5.1  4.6 5.3  4.8 4.5  4.8 4.7  4.8 4.9  4.8 5.1 /
\multiput{$\sssize \bullet$} at 4.7 4.5  4.9 4.5 /
\plot 4.7 4.5  4.9 4.5 /
\multiput{$\smallsq2$} at 5.9 4.3  5.9 4.5  5.9 4.7  5.9 4.9  5.9 5.1  5.9 5.3 /
\put{$\sssize \bullet$} at 6 4.5
\multiput{$\smallsq2$} at 7.1 4.4  7.1 4.6  7.1 4.8  7.1 5  7.1 5.2 /
\put{$\sssize \bullet$} at 7.2 4.4
\multiput{$\smallsq2$} at 9.4 4.3  9.4 4.5  9.4 4.7  9.4 4.9  9.4 5.1  9.4 5.3  9.6 4.5 /
\multiput{$\sssize \bullet$} at 9.5 4.5  9.7 4.5 /
\plot 9.5 4.5  9.7 4.5 /

\multiput{$\smallsq2$} at 1.1  5.8  1.1 6  1.1 6.2 /
\multiput{$\smallsq2$} at 3.5 5.6  3.5 5.8  3.5 6  3.5 6.2  3.5 6.4 /
\put{$\sssize \bullet$} at  3.6 5.8 
\multiput{$\smallsq2$} at 5.9 5.6  5.9 5.8  5.9 6  5.9 6.2  5.9 6.4 /
\multiput{$\smallsq2$} at 8.3 5.5  8.3 5.7  8.3 5.9  8.3 6.1  8.3 6.3  8.3 6.5 /
\put{$\sssize \bullet$} at 8.4 5.5

\multiput{$\smallsq2$} at 7.1 6.7  7.1 6.9  7.1 7.1  7.1 7.3  7.1 7.5  7.1 7.7 /
\arr{1.5 0.9} {2.1 0.3} 
\arr{2.7 0.3} {3.3 0.9} 
\arr{3.9 0.9} {4.5 0.3} 
\arr{5.1 0.3} {5.7 0.9} 
\arr{6.3 0.9} {6.9 0.3} 
\arr{7.5 0.3} {8.1 0.9} 
\arr{8.7 0.9} {9.3 0.3} 
\arr{3.9 1.2} {4.5 1.2}
\arr{5.1 1.2} {5.7 1.2}
\arr{6.3 1.2} {6.9 1.2}
\arr{7.5 1.2} {8.1 1.2}
\arr{8.7 1.2} {9.3 1.2}
\arr{1.5 1.5} {2.1 2.1} 
\arr{2.7 2.1} {3.3 1.5} 
\arr{3.9 1.5} {4.5 2.1} 
\arr{5.1 2.1} {5.7 1.5} 
\arr{6.3 1.5} {6.9 2.1} 
\arr{7.5 2.1} {8.1 1.5} 
\arr{8.7 1.5} {9.3 2.1} 
\arr{1.5 3.3} {2.1 2.7} 
\arr{2.7 2.7} {3.3 3.3} 
\arr{3.9 3.3} {4.5 2.7} 
\arr{5.1 2.7} {5.7 3.3} 
\arr{6.3 3.3} {6.9 2.7} 
\arr{7.5 2.7} {8.1 3.3} 
\arr{8.7 3.3} {9.3 2.7} 
\arr{1.5 3.9} {2.1 4.5} 
\arr{2.7 4.5} {3.3 3.9} 
\arr{3.9 3.9} {4.5 4.5} 
\arr{5.1 4.5} {5.7 3.9} 
\arr{6.3 3.9} {6.9 4.5} 
\arr{7.5 4.5} {8.1 3.9} 
\arr{8.7 3.9} {9.3 4.5} 
\arr{1.5 4.8} {2.1 4.8}
\arr{2.7 4.8} {3.3 4.8}
\arr{3.9 4.8} {4.5 4.8}
\arr{5.1 4.8} {5.7 4.8}
\arr{6.3 4.8} {6.9 4.8}
\arr{1.5 5.7} {2.1 5.1} 
\arr{2.7 5.1} {3.3 5.7} 
\arr{3.9 5.7} {4.5 5.1} 
\arr{5.1 5.1} {5.7 5.7} 
\arr{6.3 5.7} {6.9 5.1} 
\arr{7.5 5.1} {8.1 5.7} 
\arr{8.7 5.7} {9.3 5.1} 
\arr{6.3 6.3} {6.9 6.9} 
\arr{7.5 6.9} {8.1 6.3} 
\setdots<2pt>
\plot 1.2 0  1.7 0 /
\plot 3.1 0  4.1 0 /
\plot 5.5 0  6.5 0 /
\plot 7.9 0  8.9 0 /
\plot 1.9 6  2.9 6 /
\plot 4.3 6  5.3 6 /
\plot 9.1 6  9.6 6 /
\setdashes <2cm>
\plot 9.6  .4  9.6 0.7 /
\plot 9.6 1.7  9.6 1.8 /
\plot 9.6 3    9.6 4.1 /
\plot 9.6 5.5  9.6 6.3 /
\plot 1.2 -.3  1.2 .5  /
\plot 1.2 1.9  1.2 3   /
\plot 1.2 4.2  1.2 4.3 /
\plot 1.2 5.3  1.2 5.6 /
\setshadegrid span <1.5mm>
\vshade   1.2   0 6  <,z,,> 
          6   0 6  <z,z,,>
          7.2   0 7.2  <z,z,,>
          8.4   0 6    <z,z,,>
          9.6   0 6 /
\endpicture} 
$$

$$
%
%
\hbox{\beginpicture
\setcoordinatesystem units <.8cm,.8cm>
\put{$\Cal S_2(5):$} at -2 2
\multiput{$\smallsq2$} at -.1 -.4  -.1 -.2  -.1 0  -.1 .2  -.1 .4 /
\put{$\sssize\bullet$} at 0 -.4
\put{$\smallsq2$} at 1.9 0
\put{$\sssize\bullet$} at 2 0
\put{$\smallsq2 $} at 3.9 0
\multiput{$\smallsq2$} at 5.9 -.2  5.9 0  5.9 .2 /
\put{$\sssize\bullet$} at 6 0
\multiput{$\smallsq2$} at .8 .6  .8 .8  .8 1  .8 1.2  .8 1.4  1 .8 /
\multiput{$\sssize\bullet$} at .9 .8  1.1 .8 /
\plot  .9 .8  1.1 .8 /
\multiput{$\smallsq2$} at  2.9 .9  2.9 1.1 /
\put{$\sssize\bullet$} at  3 .9 
\multiput{$\smallsq2$} at  3.9 .9  3.9 1.1 /
\put{$\sssize\bullet$} at 4 1.1 
\multiput{$\smallsq2$} at  4.8 .8  4.8 1  4.8 1.2  5 1 /
\multiput{$\sssize\bullet$} at  4.9 1  5.1 1 /
\plot  4.9 1  5.1 1 /
\multiput{$\smallsq2$} at 5.9 .9  5.9 1.1 /
\multiput{$\smallsq2$} at  -.2 1.7  -.2 1.9  -.2 2.1  -.2 2.3  0 1.9 /
\multiput{$\sssize\bullet$} at  -.1 1.9  .1 1.9 /
\plot  -.1 1.9  .1 1.9 /
\multiput{$\smallsq2$} at 1.8 1.6  1.8 1.8  1.8 2  1.8 2.2  1.8 2.4  2 1.8  2 2 /
\multiput{$\sssize\bullet$} at 1.9 1.8  2.1 1.8 /
\plot 1.9 1.8  2.1 1.8 /
\multiput{$\smallsq2$} at 3.9 1.8  3.9 2  3.9 2.2 /
\put{$\sssize\bullet$} at  4 1.8 
\multiput{$\smallsq2$} at  5.8 1.7  5.8 1.9  5.8 2.1  5.8 2.3  6 1.9 /
\multiput{$\sssize\bullet$} at  5.9 1.9  6.1 1.9 /
\plot 5.9 1.9  6.1 1.9 /
\multiput{$\smallsq2$} at -.1 2.8  -.1 3  -.1 3.2 /
\put{$\sssize\bullet$} at  0 3
\multiput{$\smallsq2$} at .8 2.7  .8 2.9  .8 3.1  .8 3.3  1 2.9  1 3.1 /
\multiput{$\sssize\bullet$} at  .9 2.9  1.1 2.9 /
\plot  .9 2.9  1.1 2.9 /
\multiput{$\smallsq2$} at 1.9 2.8  1.9 3  1.9 3.2 /
\multiput{$\smallsq2$} at 2.8 2.6  2.8 2.8  2.8 3  2.8 3.2  2.8 3.4  3 2.8  3 3  3 3.2 /
\multiput{$\sssize\bullet$} at 2.9 2.8  3.1 2.8 /
\plot 2.9 2.8  3.1 2.8 /
\multiput{$\smallsq2$} at 3.9 2.6  3.9 2.8  3.9 3  3.9 3.2  3.9 3.4 /
\put{$\sssize\bullet$} at 4 2.8
\multiput{$\smallsq2$} at 4.9 2.7  4.9 2.9  4.9 3.1  4.9 3.3 /
\put{$\sssize\bullet$} at 5 2.7
\multiput{$\smallsq2$} at -.1 3.9  -.1 4.1 /
\multiput{$\smallsq2$} at 1.9 3.7  1.9 3.9  1.9 4.1  1.9 4.3 /
\put{$\sssize\bullet$} at 2 3.9
\multiput{$\smallsq2$} at 3.9 3.7  3.9 3.9  3.9 4.1  3.9 4.3 /
\multiput{$\smallsq2$} at 5.9 3.6  5.9 3.8  5.9 4  5.9 4.2  5.9 4.4 /
\put{$\sssize\bullet$} at 6 3.6
\multiput{$\smallsq2$} at  4.9 4.6  4.9 4.8  4.9 5  4.9 5.2  4.9 5.4 /

\arr{0.3 0.3} {0.7 0.7} 
\arr{1.3 0.7} {1.7 0.3}
\arr{2.3 0.3} {2.7 0.7}
\arr{3.3 0.7} {3.7 0.3}
\arr{4.3 0.3} {4.7 0.7}
\arr{5.3 0.7} {5.7 0.3}

\arr{3.3 1  } {3.7 1  }
\arr{4.3 1  } {4.7 1  }
\arr{5.3 1  } {5.7 1  }

\arr{0.3 1.7} {0.7 1.3}
\arr{1.3 1.3} {1.7 1.7}
\arr{2.3 1.7} {2.7 1.3}
\arr{3.3 1.3} {3.7 1.7}
\arr{4.3 1.7} {4.7 1.3}
\arr{5.3 1.3} {5.7 1.7}

\arr{0.3 2.3} {0.7 2.7}
\arr{1.3 2.7} {1.7 2.3}
\arr{2.3 2.3} {2.7 2.7}
\arr{3.3 2.7} {3.7 2.3}
\arr{4.3 2.3} {4.7 2.7}
\arr{5.3 2.7} {5.7 2.3}

\arr{0.3 3  } {0.7 3  }
\arr{1.3 3  } {1.7 3  }
\arr{2.3 3  } {2.7 3  }
\arr{3.3 3  } {3.7 3  }
\arr{4.3 3  } {4.7 3  }

\arr{0.3 3.7} {0.7 3.3}
\arr{1.3 3.3} {1.7 3.7}
\arr{2.3 3.7} {2.7 3.3}
\arr{3.3 3.3} {3.7 3.7}
\arr{4.3 3.7} {4.7 3.3}
\arr{5.3 3.3} {5.7 3.7}

\arr{4.3 4.3} {4.7 4.7}
\arr{5.3 4.7} {5.7 4.3}
\setdots<2pt>
\plot 0.7 0  1.3 0 /
\plot 2.7 0  3.3 0 /
\plot 4.7 0  5.3 0 /
\plot 0.7 4  1.3 4 /
\plot 2.7 4  3.3 4 /
\setdashes <10mm>
\plot 0 0.6  0 1.5 /
\plot 0 2.5  0 2.6 /
\plot 0 3.4  0 3.7 /
\plot 6 0.4  6 0.7 /
\plot 6 1.3  6 1.5 /
\plot 6 2.5  6 3.4 /
\setshadegrid span <1.5mm>
\vshade   0   0 4  <,z,,> 
          4   0 4  <z,z,,>
          5   0 5  <z,z,,>
          6   0 4 /
\endpicture} 
$$

\bigskip
\centerline{\sc The Case $m=n-1$}

\medskip
There are seven 
more pairs $m<n$ left for which the category $\Cal S_m(k[T]/T^n)$
has finite or tame type.  
In this section we are dealing with the three cases where $m=n-1$.
We have observed above that in the category 
$\Cal S(k[T]/T^n)$ the modules $P$ and $I$ are isomorphic objects.  
This projective injective module is the only indecomposable object
which is not in the category $\Cal S_{n-1}(k[T]/T^n)$. 

\medskip\noindent{\bf Lemma 10.} {\it
Let $m=n-1$.  The category $\Cal S_m(n)$ is the full subcategory of 
$\Cal S(n)$ consisting of all indecomposables in $\Cal S(n)$ with 
the exception of the projective injective modules $I[i]$, $i\in\Bbb Z$.
The Auslander-Reiten quiver for $\Cal S_m(n)$ is obtained 
as the full subquiver of the 
Auslander-Reiten quiver for $\Cal S(n)$ given by deleting the points 
corresponding to the $I[i]$. }

\smallskip \noindent {\it Proof:\/} 
Any representation $M\in\Cal S(n)$ for which the subspace is not
annihilated by $\alpha^{n-1}$ admits an embedding of the projective-injective
representation $I[i]$ in $\Cal S(n)$, for some $i$; this embedding is split
exact.
For the proof of the second assertion, consider an Auslander-Reiten sequence
in $\Cal S(n)$ which contains one of the  modules $I[i]_{\Cal S(n)}$ 
as a summand of the middle term.  Here we write 
$I[i]_{\Cal S(n)}=(I_i\sub I_i)$ for 
$I_i$ the injective $k\Bbb A_\infty^\infty/\alpha^n$-module corresponding
to the point $i$.  
$$\ssize 0\;\longrightarrow \;(\rad I_i\sub I_i)\;\lto{\big({a_i\atop c_i}\big)}\;
        \beginpicture\setcoordinatesystem units <4mm,4mm>
        \put{$\ssize(I_i\sub I_i)$} at 0 1
        \put{$\ssize\oplus$} at 0 0
        \put{$\ssize(\rad I_i/\soc I_i\sub I_i/\soc I_i)$} at 0 -1
        \endpicture
   \;\lto{\ssize (b_i\;d_i)}\; (I_i/\soc I_i\sub I_i/\soc I_i)\; \longrightarrow\; 0 $$
Note that the module $(\rad I_i\sub I_i)$ is just the object $I[i]$ in
the category $\Cal S_m(n)$ and the map $c_i$ is its source map.
Similarly, $(I_i/\soc I_i\sub I_i/\soc I_i)$ is the object $P[i+1]$
in $\Cal S_m(n)$ and $d_i$ is its sink map.  
In conclusion, the Auslander-Reiten quiver for $\Cal S_m(n)$ is obtained
from the Auslander-Reiten quiver for $\Cal S(n)$ by deleting 
the points corresponding to the projective injective modules and the 
arrows representing the above maps $a_i$ and $b_i$. \qed

\smallskip\noindent{\it Example.} Assume $m=3$ and $n=4$.
The Auslander-Reiten quiver for $\Cal S_3(4)$ is obtained from the 
Auslander-Reiten quiver for $\Cal S(4)$ by deleting the points 
corresponding to the projective and injective modules.  (We are not
deleting the modules of type $Y[i]$ 
at the bottom which are projective and relatively
injective.)

$$
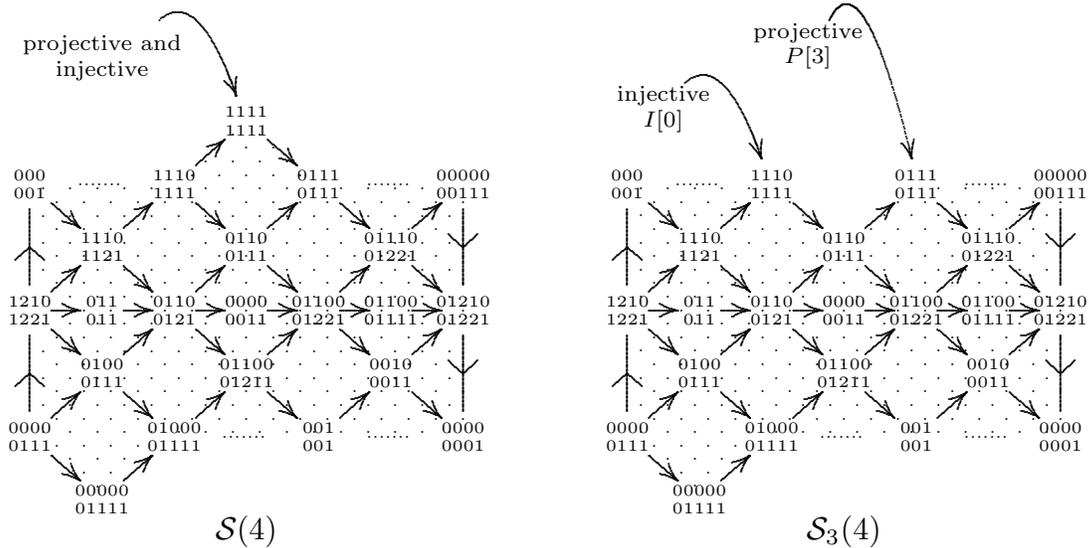
\hbox{\beginpicture
\setcoordinatesystem units <0.8cm,0.7cm>
\put{} at 0 0
\put{} at 7 5
\put{$\Cal S(4)$} at 3 -.5
\put{$\ssize {0000\atop 0111}$} at 0 1
\put{$\ssize {1210\atop 1221}$} at 0 3
\put{$\ssize {000\atop 001}$} at 0 5

\put{$\ssize {00000\atop 01111}$} at 1 0
\put{$\ssize {0100\atop 0111}$} at 1 2
\put{$\ssize {011\atop 011}$} at 1 3
\put{$\ssize {1110\atop 1121}$} at 1 4

\put{$\ssize {01000\atop 01111}$} at 2 1
\put{$\ssize {0110\atop 0121}$} at 2 3
\put{$\ssize {1110\atop 1111}$} at 2 5

\put{$\ssize {01100\atop 01211}$} at 3 2
\put{$\ssize {0000 \atop 0011}$} at 3 3
\put{$\ssize {0110 \atop 0111}$} at 3 4
\put{$\ssize {1111 \atop 1111}$} at 3 6 

\put{$\ssize {001 \atop 001}$} at 4 1
\put{$\ssize {01100 \atop 01221}$} at  4 3
\put{$\ssize {0111 \atop 0111}$} at  4 5

\put{$\ssize {01110\atop 01221}$} at 5 4
\put{$\ssize {01100\atop 01111}$} at 5 3
\put{$\ssize {0010\atop 0011}$} at 5 2

\put{$\ssize {00000\atop 00111}$} at 6 5
\put{$\ssize {01210\atop 01221}$} at 6 3
\put{$\ssize {0000\atop 0001}$} at 6 1

\arr{0.3 0.7} {0.7 0.3} 
\arr{1.3 0.3} {1.7 0.7} 
\arr{0.3 1.3} {0.7 1.7} 
\arr{1.3 1.7} {1.7 1.3} 
\arr{2.3 1.3} {2.7 1.7} 
\arr{3.3 1.7} {3.7 1.3} 
\arr{4.3 1.3} {4.7 1.7} 
\arr{5.3 1.7} {5.7 1.3} 
\arr{0.3 2.7} {0.7 2.3} 
\arr{1.3 2.3} {1.7 2.7} 
\arr{2.3 2.7} {2.7 2.3} 
\arr{3.3 2.3} {3.7 2.7} 
\arr{4.3 2.7} {4.7 2.3} 
\arr{5.3 2.3} {5.7 2.7} 
\arr{0.3 3.3} {0.7 3.7} 
\arr{1.3 3.7} {1.7 3.3} 
\arr{2.3 3.3} {2.7 3.7} 
\arr{3.3 3.7} {3.7 3.3} 
\arr{4.3 3.3} {4.7 3.7} 
\arr{5.3 3.7} {5.7 3.3} 
\arr{0.3 4.7} {0.7 4.3} 
\arr{1.3 4.3} {1.7 4.7} 
\arr{2.3 4.7} {2.7 4.3} 
\arr{3.3 4.3} {3.7 4.7} 
\arr{4.3 4.7} {4.7 4.3} 
\arr{5.3 4.3} {5.7 4.7} 
\arr{2.3 5.3} {2.7 5.7} 
\arr{3.3 5.7} {3.7 5.3} 
\arr{0.3 3}{0.7 3}
\arr{1.3 3}{1.7 3}
\arr{2.3 3}{2.7 3}
\arr{3.3 3}{3.7 3}
\arr{4.3 3}{4.7 3}
\arr{5.3 3}{5.7 3}
\setdots<2pt>
\plot 2.7 1  3.3 1 /
\plot 4.7 1  5.3 1 /
\plot 0.7 5  1.3 5 /
\plot 4.7 5  5.3 5 /
\setsolid
\plot 0 1.4  0 2.6 /
\plot 0 3.4  0 4.6 /
\plot -.2 1.8  0 2  .2 1.8 /
\plot -.2 3.8  0 4  .2 3.8 /
\plot 6 1.4  6 2.6 /
\plot 6 3.4  6 4.6 /
\plot 5.8 2.2  6 2  6.2 2.2 /
\plot 5.8 4.2  6 4  6.2 4.2 /
\put{$\ssize\text{projective and}$} at 1 7.2
\put{$\ssize\text{injective}$} at 1 6.8
\setshadegrid span <1.5mm>
\vshade   0   1 5  <,z,,> 
          1   0 5  <z,z,,>
          2   1 5  <z,z,,> 
          3   1 6  <z,z,,> 
          4   1 5  <z,z,,>
          6   1 5 /
\setquadratic
\plot 1.8 7.6  2.3 7.6  2.8 6.6 /
\arr{2.8 6.6} {2.86 6.4}
\endpicture} 
\qquad
\hbox{\beginpicture
\setcoordinatesystem units <0.8cm,0.7cm>
\put{} at 0 0
\put{} at 7 5
\put{$\Cal S_3(4)$} at 3 -.5
\put{$\ssize {0000\atop 0111}$} at 0 1
\put{$\ssize {1210\atop 1221}$} at 0 3
\put{$\ssize {000\atop 001}$} at 0 5

\put{$\ssize {00000\atop 01111}$} at 1 0
\put{$\ssize {0100\atop 0111}$} at 1 2
\put{$\ssize {011\atop 011}$} at 1 3
\put{$\ssize {1110\atop 1121}$} at 1 4

\put{$\ssize {01000\atop 01111}$} at 2 1
\put{$\ssize {0110\atop 0121}$} at 2 3
\put{$\ssize {1110\atop 1111}$} at 2 5

\put{$\ssize {01100\atop 01211}$} at 3 2
\put{$\ssize {0000 \atop 0011}$} at 3 3
\put{$\ssize {0110 \atop 0111}$} at 3 4

\put{$\ssize {001 \atop 001}$} at 4 1
\put{$\ssize {01100 \atop 01221}$} at  4 3
\put{$\ssize {0111 \atop 0111}$} at  4 5

\put{$\ssize {01110\atop 01221}$} at 5 4
\put{$\ssize {01100\atop 01111}$} at 5 3
\put{$\ssize {0010\atop 0011}$} at 5 2

\put{$\ssize {00000\atop 00111}$} at 6 5
\put{$\ssize {01210\atop 01221}$} at 6 3
\put{$\ssize {0000\atop 0001}$} at 6 1

\arr{0.3 0.7} {0.7 0.3} 
\arr{1.3 0.3} {1.7 0.7} 
\arr{0.3 1.3} {0.7 1.7} 
\arr{1.3 1.7} {1.7 1.3} 
\arr{2.3 1.3} {2.7 1.7} 
\arr{3.3 1.7} {3.7 1.3} 
\arr{4.3 1.3} {4.7 1.7} 
\arr{5.3 1.7} {5.7 1.3} 
\arr{0.3 2.7} {0.7 2.3} 
\arr{1.3 2.3} {1.7 2.7} 
\arr{2.3 2.7} {2.7 2.3} 
\arr{3.3 2.3} {3.7 2.7} 
\arr{4.3 2.7} {4.7 2.3} 
\arr{5.3 2.3} {5.7 2.7} 
\arr{0.3 3.3} {0.7 3.7} 
\arr{1.3 3.7} {1.7 3.3} 
\arr{2.3 3.3} {2.7 3.7} 
\arr{3.3 3.7} {3.7 3.3} 
\arr{4.3 3.3} {4.7 3.7} 
\arr{5.3 3.7} {5.7 3.3} 
\arr{0.3 4.7} {0.7 4.3} 
\arr{1.3 4.3} {1.7 4.7} 
\arr{2.3 4.7} {2.7 4.3} 
\arr{3.3 4.3} {3.7 4.7} 
\arr{4.3 4.7} {4.7 4.3} 
\arr{5.3 4.3} {5.7 4.7} 
\arr{0.3 3}{0.7 3}
\arr{1.3 3}{1.7 3}
\arr{2.3 3}{2.7 3}
\arr{3.3 3}{3.7 3}
\arr{4.3 3}{4.7 3}
\arr{5.3 3}{5.7 3}
\setdots<2pt>
\plot 2.7 1  3.3 1 /
\plot 4.7 1  5.3 1 /
\plot 0.7 5  1.3 5 /
\plot 4.7 5  5.3 5 /
\setsolid
\plot 0 1.4  0 2.6 /
\plot 0 3.4  0 4.6 /
\plot -.2 1.8  0 2  .2 1.8 /
\plot -.2 3.8  0 4  .2 3.8 /
\plot 6 1.4  6 2.6 /
\plot 6 3.4  6 4.6 /
\plot 5.8 2.2  6 2  6.2 2.2 /
\plot 5.8 4.2  6 4  6.2 4.2 /
\put{$\ssize\text{injective}$} at 0.5 6.4
\put{$\ssize I[0]$} at 0.5 6
\setshadegrid span <1.5mm>
\vshade   0   1 5  <,z,,> 
          1   0 5  <z,z,,>
          2   1 5  <z,z,,> 
          3   1 4  <z,z,,> 
          4   1 5  <z,z,,>
          6   1 5 /
\setquadratic
\plot .8 6.6  1.3 6.6  1.8 5.6 /
\arr{1.8 5.6} {1.86 5.4}
\put{$\ssize\text{projective}$} at 2.5 7.4
\put{$\ssize P[3]$} at 2.5 7
\plot 2.8 7.6  3.3 7.6  3.9 5.6 /
\arr{3.9 5.6} {3.94 5.4}
\endpicture} 
$$

Using covering theory, one obtains the analog of Lemma~10 for categories
of type $\Cal S_{n-1}(k[T]/T^n)$.  

\smallskip\noindent{\it Example.} The Auslander-Reiten quiver for 
$\Cal S_3(k[T]/T^4)$ is obtained by identifying the left border with
the right along the arrows in the above diagram. 
By deleting the projective injective module in $\Cal S(k[T]/T^4)$,
the stable type of the category is reduced from $\Bbb Z\Bbb E_6/\sigma\tau^3$
to $\Bbb Z\Bbb D_4/\sigma'\tau^3$ where $\sigma$ and 
$\sigma'$ are automorphisms of order two for
the Dynkin diagrams $\Bbb E_6$ and $\Bbb D_4$, respectively.
Thus, in the Auslander-Reiten quiver for $\Cal S_3(k[T]/T^4)$,
there is a stable orbit of period three in the middle
with a stable orbit of period six above and below it.  On the outside of this
second orbit is a nonstable orbit, also of length six, which has 
the projective injective module $Y$ attached to it.
This part of the quiver can be realized as a Moebius band in three dimensions.
However, there is a second stable orbit of period three, also attached to
the orbit in the middle.  This ``balcony'' cannot be realized in three
dimensions. 

$$
\hbox{\beginpicture
\setcoordinatesystem units <.8cm,.8cm>
\put{} at 2 0
\put{} at 8 6
\put{$\Cal S_3(k[T]/T^4):$} at 0 2.75
\multiput{$\smallsq2$} at 1.9 .8  1.9 1  1.9 1.2 /
\multiput{$\smallsq2$} at 1.8 2.7  1.8 2.9  1.8 3.1  1.8 3.3  2 2.9  2 3.1 /
\multiput{$\sssize\bullet$} at 1.9 3.1  2.1 3.1  2.1 2.9 /
\plot 1.9 3.1  2.1 3.1 /
\put{$\smallsq2$} at 1.9 5

\multiput{$\smallsq2$} at 2.9 -.3  2.9 -.1  2.9 .1  2.9 .3 /
\multiput{$\smallsq2$} at 2.9 1.8  2.9 2  2.9 2.2 /
\put{$\sssize\bullet$} at 3 1.8
\multiput{$\smallsq2$} at 2.9 2.65  2.9 2.85 /
\put{$\sssize\bullet$} at 3 2.85 
\multiput{$\smallsq2$} at 2.8 3.7  2.8 3.9  2.8 4.1  2.8 4.3  3 4.1 /
\multiput{$\sssize\bullet$} at 2.9 4.1  3.1 4.1 /
\plot 2.9 4.1  3.1 4.1 /

\multiput{$\smallsq2$} at 3.9 .7  3.9 .9  3.9 1.1  3.9 1.3 /
\put{$\sssize\bullet$} at 4 .7
\multiput{$\smallsq2$} at 3.8 2.8  3.8 3  3.8 3.2  4 3 /
\multiput{$\sssize\bullet$} at 3.9 3  4.1 3 /
\plot 3.9 3  4.1 3 /
\multiput{$\smallsq2$} at 3.9 4.7  3.9 4.9  3.9 5.1  3.9 5.3 /
\put{$\sssize\bullet$} at 4 5.1 

\multiput{$\smallsq2$} at 4.8 1.7  4.8 1.9  4.8 2.1  4.8 2.3  5 1.9 /
\multiput{$\sssize \bullet$} at 4.9 1.9  5.1 1.9 /
\plot 4.9 1.9  5.1 1.9 /
\multiput{$\smallsq2$} at 4.9 2.65  4.9 2.85 /
\multiput{$\smallsq2$} at 4.9 3.8  4.9 4  4.9 4.2 /
\put{$\sssize\bullet$} at 5 4

\put{$\smallsq2$} at 5.9 1
\put{$\sssize \bullet$} at 6 1
\multiput{$\smallsq2$} at 5.8 2.7  5.8 2.9  5.8 3.1  5.8 3.3  6 2.9  6 3.1 /
\multiput{$\sssize \bullet$} at 5.9 2.9  6.1 2.9 /
\plot 5.9 2.9  6.1 2.9 /
\multiput{$\smallsq2$} at 5.9 4.8  5.9 5  5.9 5.2 /
\put{$\sssize\bullet$} at 6 5.2

\multiput{$\smallsq2$} at 6.8 3.7  6.8 3.9  6.8 4.1  6.8 4.3  7 3.9  7 4.1 /
\multiput{$\sssize\bullet$} at 6.9 4.1  7.1 4.1 /
\plot 6.9 4.1  7.1 4.1 /
\multiput{$\smallsq2$} at 6.9 2.45  6.9 2.65  6.9 2.85  6.9 3.05 /
\put{$\sssize\bullet$} at 7 2.65
\multiput{$\smallsq2$} at 6.9 1.9  6.9 2.1 /
\put{$\sssize\bullet$} at 7 1.9

\multiput{$\smallsq2$} at 7.9 4.8  7.9 5  7.9 5.2 /
\multiput{$\smallsq2$} at 7.8 2.7  7.8 2.9  7.8 3.1  7.8 3.3  8 2.9  8 3.1 /
\multiput{$\sssize\bullet$} at 7.9 3.1  8.1 3.1  8.1 2.9 /
\plot 7.9 3.1  8.1 3.1 /
\put{$\smallsq2$} at 7.9 1

\arr{2.3 0.7} {2.7 0.3} 
\arr{3.3 0.3} {3.7 0.7} 
\arr{6.3 1.3} {6.7 1.7} 
\arr{7.3 1.7} {7.7 1.3} 
\arr{2.3 1.3} {2.7 1.7} 
\arr{3.3 1.7} {3.7 1.3} 
\arr{4.3 1.3} {4.7 1.7} 
\arr{5.3 1.7} {5.7 1.3} 
\arr{6.3 2.7} {6.7 2.3} 
\arr{7.3 2.3} {7.7 2.7} 
\arr{2.3 2.7} {2.7 2.3} 
\arr{3.3 2.3} {3.7 2.7} 
\arr{4.3 2.7} {4.7 2.3} 
\arr{5.3 2.3} {5.7 2.7} 
\arr{6.3 3.3} {6.7 3.7} 
\arr{7.3 3.7} {7.7 3.3} 
\arr{2.3 3.3} {2.7 3.7} 
\arr{3.3 3.7} {3.7 3.3} 
\arr{4.3 3.3} {4.7 3.7} 
\arr{5.3 3.7} {5.7 3.3} 
\arr{6.3 4.7} {6.7 4.3} 
\arr{7.3 4.3} {7.7 4.7} 
\arr{2.3 4.7} {2.7 4.3} 
\arr{3.3 4.3} {3.7 4.7} 
\arr{4.3 4.7} {4.7 4.3} 
\arr{5.3 4.3} {5.7 4.7} 
\arr{6.3 2.925}{6.7 2.825}
\arr{7.3 2.825}{7.7 2.925}
\arr{2.3 2.925}{2.7 2.825}
\arr{3.3 2.825}{3.7 2.925}
\arr{4.3 2.925}{4.7 2.825}
\arr{5.3 2.825}{5.7 2.925}
\setdots<2pt>
\plot 6.3 1  7.7 1 /
\plot 4.3 1  5.7 1 /
\plot 6.3 5  7.7 5 /
\plot 2.3 5  3.7 5 /
\setsolid
\plot 2 1.4  2 2.5 /
\plot 2 3.5  2 4.8 /
\plot 1.8 1.8  2 2  2.2 1.8 /
\plot 1.8 3.8  2 4  2.2 3.8 /
\plot 8 1.2  8 2.5 /
\plot 8 3.5  8 4.6 /
\plot 7.8 2.2  8 2  8.2 2.2 /
\plot 7.8 4.2  8 4  8.2 4.2 /
\setshadegrid span <1.5mm>
\vshade   2   1 5  <,z,,> 
          3   0 5  <z,z,,>
          4   1 5  <z,z,,> 
          5   1 4  <z,z,,> 
          6   1 5  <z,z,,>
          8   1 5 /
\endpicture} 
$$

\bigskip
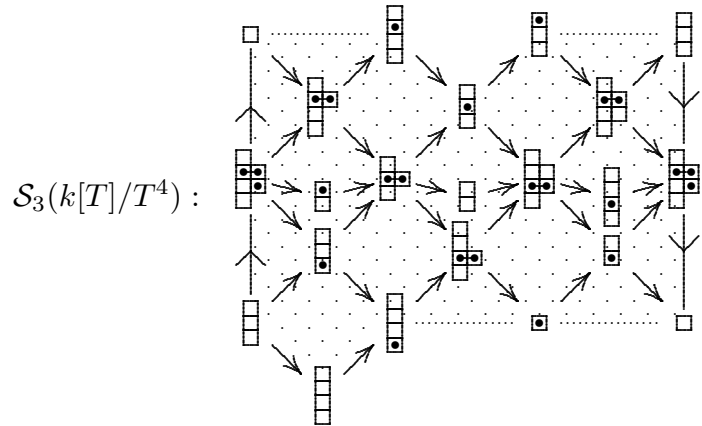
\centerline {\sc The Case $m=3$, $n=5$.}

\medskip
In this example, we obtain the Auslander-Reiten quiver for $\Cal S_3(5)$ 
from the Auslander-Reiten quiver of the path algebra of the following
quiver with commutativity and nilpotence relations as indicated.
We omit the  details as the process appears similar to the cases where
$m=2$, and even easier than the case $m=3$ and $n=7$ below.
$$
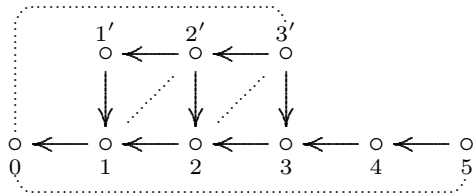
\hbox{\beginpicture
\setcoordinatesystem units <0.5cm,0.5cm>
\put{} at 0 0
\put{} at 5 2
\put{$\circ$} at 0 0
\put{$\circ$} at 2 0
\put{$\circ$} at 4 0
\put{$\circ$} at 6 0
\put{$\circ$} at 8 0
\put{$\circ$} at 10 0

\put{$\circ$} at 2 2
\put{$\circ$} at 4 2
\put{$\circ$} at 6 2
\arr{1.6 0}{0.4 0}
\arr{3.6 0}{2.4 0}
\arr{5.6 0}{4.4 0}
\arr{7.6 0}{6.4 0}
\arr{9.6 0}{8.4 0}

\arr{3.6 2}{2.4 2}
\arr{5.6 2}{4.4 2}

\arr{2 1.6}{2 0.4}
\arr{4 1.6}{4 0.4}
\arr{6 1.6}{6 0.4}


\put{$\ssize 1'$} at  2 2.5
\put{$\ssize 2'$} at  4 2.5
\put{$\ssize 3'$} at  6 2.5

\put{$\ssize 0$} at 0 -.5
\put{$\ssize 1$} at 2 -.5
\put{$\ssize 2$} at 4 -.5
\put{$\ssize 3$} at 6 -.5
\put{$\ssize 4$} at 8 -.5
\put{$\ssize 5$} at 10 -.5

\setdots<2pt>
\plot 3.5 1.5  2.5 0.5 /
\plot 5.5 1.5  4.5 0.5 /
\plot         0 -0.7  0.1 -0.9  0.2 -1   0.3 -1.05  9.7 -1.05  9.8 -1   9.9 -0.9  10 -0.7 /
\plot 0 0.3   0 2.7   0.1  2.9  0.2 3    0.3 3.05   5.7  3.05  5.8  3   5.9  2.9   6  2.7 /
\endpicture}
$$

\medskip\noindent
{\bf Proposition 11.} {\it
The category $\Cal S_3(k[T]/T^5)$ has 37 indecomposable objects, up
to isomorphism.  The Auslander-Reiten quiver has five stable orbits,
four of length 5 and one of length 10; the stable type is
$\Bbb Z\Bbb D_6/\tau^5\sigma$ where $\sigma$ permutes the two short
branches of the Dynkin diagram $\Bbb D_6$ and hence creates the
long orbit in the Auslander-Reiten quiver.  Attached to this long orbit
is a nonstable orbit from $P$ to $I$ of length 6; the projective 
injective module $Y$ is attached to the short orbit corresponding to the 
endpoint of the long branch of the diagram $\Bbb D_6$.  }

$$
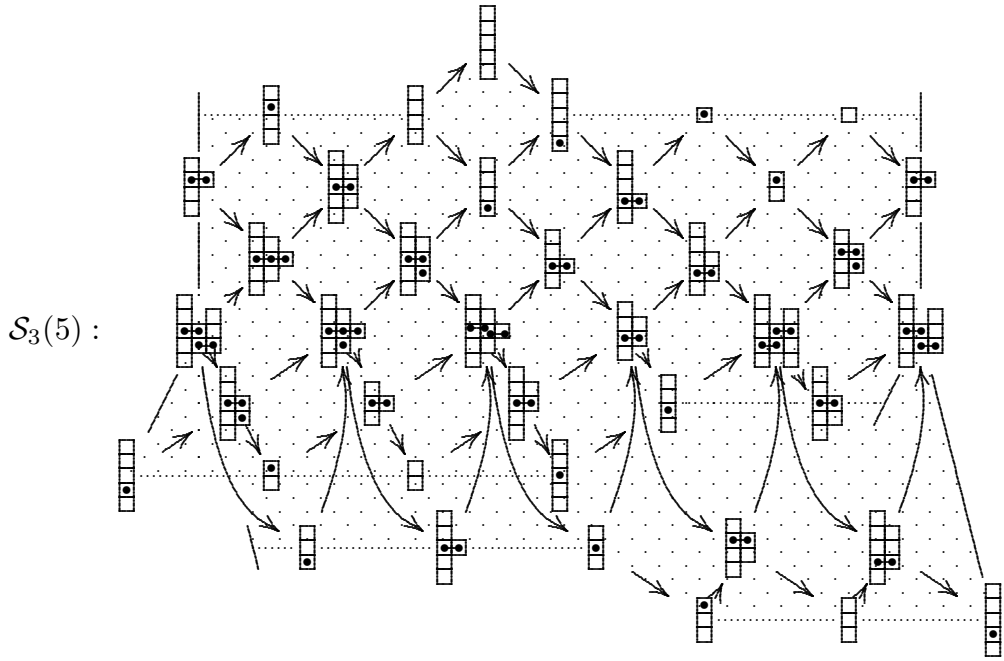
\hbox{\beginpicture
\setcoordinatesystem units <.8cm,.8cm>
\put{} at -1 0
\put{} at 11 8
\put{$\Cal S_3(5):$} at -2 4
\multiput{$\smallsq2$} at 3.9 7.6  3.9 7.8  3.9 8  3.9 8.2  3.9 8.4 /

\multiput{$\smallsq2$} at .9 6.7  .9 6.9  .9 7.1  .9 7.3 /
\put{$\sssize\bullet$} at 1 7.1
\multiput{$\smallsq2$} at 2.9 6.7  2.9 6.9  2.9 7.1  2.9 7.3 /
\multiput{$\smallsq2$} at 4.9 6.6  4.9 6.8  4.9 7  4.9 7.2  4.9 7.4 /
\put{$\sssize\bullet$} at 5 6.6
\put{$\smallsq2$} at 6.9 7
\put{$\sssize\bullet$} at 7 7 
\put{$\smallsq2$} at 8.9 7

\multiput{$\smallsq2$} at -.2 5.7  -.2 5.9  -.2 6.1  -.2 6.3  0 6.1 /
\multiput{$\sssize\bullet$} at -.1 6.1  .1 6.1 /
\plot -.1 6.1  .1 6.1 /
\multiput{$\smallsq2$} at 1.8 5.6  1.8 5.8  1.8 6  1.8 6.2  1.8 6.4  2 5.8  2 6  2 6.2 /
\multiput{$\sssize\bullet$} at 1.9 6  2.1 6 /
\plot 1.9 6 2.1 6 /
\multiput{$\smallsq2$} at 3.9 5.7  3.9 5.9  3.9 6.1  3.9 6.3 /
\put{$\sssize\bullet$} at 4 5.7
\multiput{$\smallsq2$} at 5.8 5.6  5.8 5.8  5.8 6  5.8 6.2  5.8 6.4  6 5.8 /
\multiput{$\sssize\bullet$} at 5.9 5.8  6.1 5.8 /
\plot 5.9 5.8  6.1 5.8 /
\multiput{$\smallsq2$} at 7.9 5.9  7.9 6.1 /
\put{$\sssize\bullet$} at 8 6.1
\multiput{$\smallsq2$} at 9.8 5.7  9.8 5.9  9.8 6.1  9.8 6.3  10 6.1 /
\multiput{$\sssize\bullet$} at 9.9 6.1  10.1 6.1 /
\plot 9.9 6.1  10.1 6.1 /

\multiput{$\smallsq2$} at .7 4.6  .7 4.8  .7 5  .7 5.2  .7 5.4  .9 4.8  .9 5  .9 5.2  1.1 5 /
\multiput{$\sssize\bullet$} at .8 5  1 5  1.2 5 /
\plot .8 5  1.2 5 /
\multiput{$\smallsq2$} at 2.8 4.6  2.8 4.8  2.8 5  2.8 5.2  2.8 5.4  3 4.8  3 5  3 5.2 /
\multiput{$\sssize\bullet$} at 2.9 5  3.1 5  3.1 4.8 /
\plot 2.9 5  3.1 5 /
\multiput{$\smallsq2$} at 4.8 4.7  4.8 4.9  4.8 5.1  4.8 5.3  5 4.9 /
\multiput{$\sssize\bullet$} at 4.9 4.9  5.1 4.9 /
\plot 4.9 4.9  5.1 4.9 /
\multiput{$\smallsq2$} at 6.8 4.6  6.8 4.8  6.8 5  6.8 5.2  6.8 5.4  7 4.8  7 5 /
\multiput{$\sssize\bullet$} at 6.9 4.8  7.1 4.8 /
\plot 6.9 4.8  7.1 4.8 /
\multiput{$\smallsq2$} at 8.8 4.7  8.8 4.9  8.8 5.1  8.8 5.3  9 4.9  9 5.1 /
\multiput{$\sssize\bullet$} at 8.9 5.1  9.1 5.1  9.1 4.9 /
\plot 8.9 5.1  9.1 5.1 /

\multiput{$\smallsq2$} at -.3 3.6  -.3 3.8  -.3 4  -.3 4.2  -.3 4.4  -.1 3.8  -.1 4  .1 3.8  .1 4  .1 4.2 /
\multiput{$\sssize\bullet$} at -.2 4  0 4  0 3.8  .2 3.8 /
\plot -.2 4  0 4 /
\plot 0 3.8  .2 3.8 /
\multiput{$\smallsq2$} at 1.7 3.6  1.7 3.8  1.7 4  1.7 4.2  1.7 4.4  1.9 3.8  1.9 4  1.9 4.2  2.1 4 /
\multiput{$\sssize\bullet$} at 1.8 4  2 4  2.2 4  2 3.8 /
\plot 1.8 4  2.2 4 /
\multiput{$\smallsq2$} at 3.7 3.6  3.7 3.8  3.7 4  3.7 4.2  3.7 4.4  3.9 3.8  3.9 4  3.9 4.2  4.1 4 /
\multiput{$\sssize\bullet$} at 3.76 4.04  3.96 4.04  4.04 3.96  4.24 3.96 /
\plot 3.76 4.04  3.96 4.04 /
\plot 4.04 3.96  4.24 3.96 /
\multiput{$\smallsq2$} at 5.8 3.7  5.8 3.9  5.8 4.1  5.8 4.3  6 3.9  6 4.1 /
\multiput{$\sssize\bullet$} at 5.9 3.9  6.1 3.9 /
\plot 5.9 3.9  6.1 3.9 /
\multiput{$\smallsq2$} at 7.7 3.6  7.7 3.8  7.7 4  7.7 4.2  7.7 4.4  
        7.9 3.8  7.9 4  8.1 3.6  8.1 3.8  8.1 4  8.1 4.2 /
\multiput{$\sssize\bullet$} at 7.8 3.8  8 3.8  8 4  8.2 4 /
\plot 7.8 3.8  8 3.8 /
\plot 8 4  8.2 4 /
\multiput{$\smallsq2$} at 9.7 3.6  9.7 3.8  9.7 4  9.7 4.2  9.7 4.4  9.9 3.8  9.9 4  10.1 3.8  10.1 4  10.1 4.2 /
\multiput{$\sssize\bullet$} at 9.8 4  10 4  10 3.8  10.2 3.8 /
\plot 9.8 4  10 4 /
\plot 10 3.8  10.2 3.8 /

\multiput{$\smallsq2$} at .3 2.6  .3 2.8  .3 3  .3 3.2  .3 3.4  .5 2.8  .5 3 /
\multiput{$\sssize\bullet$} at .4 3  .6 3  .6 2.8 /
\plot .4 3  .6 3 /
\multiput{$\smallsq2$} at 2.3 2.8  2.3 3  2.3 3.2  2.5 3 /
\multiput{$\sssize\bullet$} at 2.4 3  2.6 3 /
\plot 2.4 3  2.6 3 /
\multiput{$\smallsq2$} at 4.3 2.6  4.3 2.8  4.3 3  4.3 3.2  4.3 3.4  4.5 3  4.5 3.2 /
\multiput{$\sssize\bullet$} at 4.4 3  4.6 3 /
\plot 4.4 3  4.6 3 /
\multiput{$\smallsq2$} at 6.4 2.7  6.4 2.9  6.4 3.1  6.4 3.3 /
\put{$\sssize\bullet$} at 6.5 2.9
\multiput{$\smallsq2$} at 8.5 2.6  8.5 2.8  8.5 3  8.5 3.2  8.5 3.4  8.7 2.8  
        8.7 3 /
\multiput{$\sssize\bullet$} at 8.6 3  8.8 3 /
\plot 8.6 3  8.8 3 /

\multiput{$\smallsq2$} at -1.1 1.6  -1.1 1.8  -1.1 2  -1.1 2.2  -1.1 2.4 /
\put{$\sssize\bullet$} at -1 1.8
\multiput{$\smallsq2$} at .9 1.9  .9 2.1 /
\put{$\sssize\bullet$} at 1 2.1
\multiput{$\smallsq2$} at 2.9 1.9  2.9 2.1 /
\multiput{$\smallsq2$} at 4.9 1.6  4.9 1.8  4.9 2  4.9 2.2  4.9 2.4 /
\put{$\sssize\bullet$} at 5 2

\multiput{$\smallsq2$} at 1.4 .8  1.4 1  1.4 1.2 /
\put{$\sssize\bullet$} at 1.5 .8
\multiput{$\smallsq2$} at 3.3 .6  3.3 .8  3.3 1  3.3 1.2  3.3 1.4  3.5 1 /
\multiput{$\sssize\bullet$} at 3.4 1  3.6 1 /
\plot 3.4 1  3.6 1 /
\multiput{$\smallsq2$} at 5.4 .8  5.4 1  5.4 1.2 /
\put{$\sssize\bullet$} at 5.5 1
\multiput{$\smallsq2$} at 7.3 .7  7.3 .9  7.3 1.1  7.3 1.3  7.5 .9  7.5 1.1 /
\multiput{$\sssize\bullet$} at 7.4 1.1  7.6 1.1 /
\plot 7.4 1.1  7.6 1.1 /
\multiput{$\smallsq2$} at 9.3 .6  9.3 .8  9.3 1  9.3 1.2  9.3 1.4  9.5 .8  9.5 1  9.5 1.2 /
\multiput{$\sssize\bullet$} at 9.4 .8  9.6 .8 /
\plot 9.4 .8  9.6 .8 /

\multiput{$\smallsq2$} at 6.9 -.2  6.9 0  6.9 .2 /
\put{$\sssize\bullet$} at 7 .2
\multiput{$\smallsq2$} at 8.9 -.2  8.9 0  8.9 .2 /
\multiput{$\smallsq2$} at 10.9 -.4  10.9 -.2  10.9 0  10.9 .2  10.9 .4 /
\put{$\sssize\bullet$} at 11 -.2

\arr{6 .67} {6.5 .33}
\arr{8 .67} {8.5 .33}
\arr{10 .67} {10.5 .33}
\arr{7.15 .3} {7.25 .5}
\arr{9.15 .3} {9.25 .5}

\arr{-.5 2.33} {-.1 2.6}
\arr{0.65 2.7} {0.85 2.3}
\arr{1.5 2.33} {1.9 2.6}
\arr{2.65 2.7} {2.85 2.3}
\arr{3.5 2.33} {3.9 2.6}
\arr{4.65 2.7} {4.85 2.3}
\arr{5.5 2.33} {5.9 2.6}

\arr{0.15 3.7} {0.25 3.5}
\arr{1.0 3.33} {1.5 3.67}
\arr{2.15 3.7} {2.25 3.5}
\arr{3.0 3.33} {3.5 3.67}
\arr{4.15 3.7} {4.25 3.5}
\arr{5.0 3.33} {5.5 3.67}
\arr{6.15 3.7} {6.25 3.5}
\arr{7.0 3.33} {7.5 3.67}
\arr{8.3 3.4} {8.4 3.2}
\arr{9.0 3.33} {9.5 3.67}

\arr{0.4 4.4} {0.6 4.6}
\arr{1.3 4.7} {1.6 4.4}
\arr{2.3 4.3} {2.7 4.7}
\arr{3.3 4.7} {3.6 4.4}
\arr{4.3 4.3} {4.7 4.7}
\arr{5.3 4.7} {5.7 4.3}
\arr{6.3 4.3} {6.7 4.7}
\arr{7.3 4.7} {7.6 4.4}
\arr{8.4 4.4} {8.7 4.7}
\arr{9.3 4.7} {9.6 4.4}

\arr{0.3 5.7} {0.6 5.4}
\arr{1.3 5.3} {1.7 5.7}
\arr{2.3 5.7} {2.7 5.3}
\arr{3.3 5.3} {3.7 5.7}
\arr{4.3 5.7} {4.7 5.3}
\arr{5.3 5.3} {5.7 5.7}
\arr{6.3 5.7} {6.7 5.3}
\arr{7.3 5.3} {7.7 5.7}
\arr{8.3 5.7} {8.7 5.3}
\arr{9.3 5.3} {9.7 5.7}

\arr{0.3 6.3} {0.7 6.7}
\arr{1.3 6.7} {1.7 6.3}
\arr{2.3 6.3} {2.7 6.7}
\arr{3.3 6.7} {3.7 6.3}
\arr{4.3 6.3} {4.7 6.7}
\arr{5.3 6.7} {5.7 6.3}
\arr{6.3 6.3} {6.7 6.7}
\arr{7.3 6.7} {7.7 6.3}
\arr{8.3 6.3} {8.7 6.7}
\arr{9.3 6.7} {9.7 6.3}

\arr{3.3 7.3} {3.7 7.7}
\arr{4.3 7.7} {4.7 7.3}

\arr{1 1.3} {1.1 1.23}
\arr{3 1.3} {3.1 1.23}
\arr{5 1.3} {5.1 1.23}
\arr{7 1.3} {7.1 1.23}
\arr{9 1.3} {9.1 1.23}

\arr{2 3.4} {1.99 3.5}
\arr{4 3.4} {3.99 3.5}
\arr{6 3.4} {5.99 3.5}
\arr{8 3.4} {7.99 3.5}
\arr{10 3.4} {9.99 3.5}

\setdots<2pt>
\plot 0 7  .8 7 /
\plot 1.2 7  2.8 7 /
\plot 5.2 7  6.8 7 /
\plot 7.2 7  8.8 7 /
\plot 9.2 7  10 7 /

\plot -.8 2  .8 2 /
\plot 1.2 2  2.8 2 /
\plot 3.2 2  4.8 2 /

\plot 6.7 3  8.4 3 /
\plot 9 3  9.5 3 /

\plot .75 1 1.3 1 /
\plot 1.7 1  3.2 1 /
\plot 3.8 1  5.3 1 /

\plot 7.2 0  8.8 0 /
\plot 9.2 0  10.8 0 /

\setsolid
\plot 0 4.6  0 5.5 /
\plot 0 6.5  0 7.3 /
\plot -.3 3.4  -.7 2.6 /
\plot .675 1.3  .825 0.7 /

\plot 10 4.6  10 5.5 /
\plot 10 6.5  10 7.3 /
\plot 9.7 3.4  9.35 2.7 /
\plot 10.15 3.4  10.85 .6 /

\setshadegrid span <1.5mm>
\vshade  -1   2 2  <,z,,> 
         -.01 2 4  <z,z,,>
          0   2 7  <z,z,,>
          .5  2 7  <z,z,,>
          .75 1 7  <z,z,,>
          3   1 7  <z,z,,>
          4   1 8  <z,z,,>
          5   1 7  <z,z,,>
          5.5 1 7  <z,z,,>
          7   0 7  <z,z,,>
          9.99 0 7 <z,z,,>
         10   0 4  <z,z,,>
         11   0 0 /

\setquadratic
\plot .05 3.5  .4 2  1 1.3 /
\plot 2.05 3.5  2.4 2  3 1.3 /
\plot 4.05 3.5  4.4 2  5 1.3 /
\plot 6.05 3.5  6.4 2  7 1.3 /
\plot 8.05 3.5  8.4 2  9 1.3 /

\plot 1.7 1.5  2 2.6  2 3.4 /
\plot 3.7 1.5  4 2.6  4 3.4 /
\plot 5.7 1.5  6 2.6  6 3.4 /
\plot 7.7 1.5  8 2.6  8 3.4 /
\plot 9.7 1.5  10 2.6  10 3.4 /

\endpicture} 
$$

\bigskip
\centerline {\sc Case $m=3$, $n=6$.}

\medskip
One obtains the indecomposables in, and the 
Auslander-Reiten quiver for the category $\Cal S_3(6)$ 
from the Auslander-Reiten quiver for the path algebra of the 
following quiver, with the relations as indicated.  Again, we omit the proof.

$$
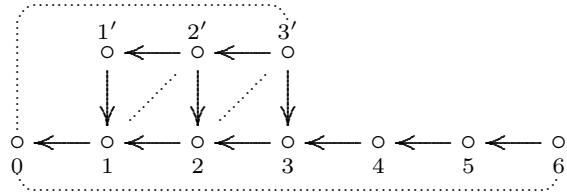
\hbox{\beginpicture
\setcoordinatesystem units <0.5cm,0.5cm>
\put{} at 0 0
\put{} at 5 2
\put{$\circ$} at 0 0
\put{$\circ$} at 2 0
\put{$\circ$} at 4 0
\put{$\circ$} at 6 0
\put{$\circ$} at 8 0
\put{$\circ$} at 10 0
\put{$\circ$} at 12 0

\put{$\circ$} at 2 2
\put{$\circ$} at 4 2
\put{$\circ$} at 6 2
\arr{1.6 0}{0.4 0}
\arr{3.6 0}{2.4 0}
\arr{5.6 0}{4.4 0}
\arr{7.6 0}{6.4 0}
\arr{9.6 0}{8.4 0}
\arr{11.6 0}{10.4 0}

\arr{3.6 2}{2.4 2}
\arr{5.6 2}{4.4 2}

\arr{2 1.6}{2 0.4}
\arr{4 1.6}{4 0.4}
\arr{6 1.6}{6 0.4}


\put{$\ssize 1'$} at  2 2.5
\put{$\ssize 2'$} at  4 2.5
\put{$\ssize 3'$} at  6 2.5

\put{$\ssize 0$} at 0 -.5
\put{$\ssize 1$} at 2 -.5
\put{$\ssize 2$} at 4 -.5
\put{$\ssize 3$} at 6 -.5
\put{$\ssize 4$} at 8 -.5
\put{$\ssize 5$} at 10 -.5
\put{$\ssize 6$} at 12 -.5

\setdots<2pt>
\plot 3.5 1.5  2.5 0.5 /
\plot 5.5 1.5  4.5 0.5 /
\plot         0 -0.7  0.1 -0.9  0.2 -1   0.3 -1.05 11.7 -1.05 11.8 -1  11.9 -0.9  12 -0.7 /
\plot 0 0.3   0 2.7   0.1  2.9  0.2 3    0.3 3.05   5.7  3.05  5.8  3   5.9  2.9   6  2.7 /
\endpicture}
$$

%
%
\def\arqthreesix{%
\beginpicture\setcoordinatesystem units <4.5mm,4mm>
\put{The Category $\Cal S_3(6)$} at 13 -4
\put{} at 42 0
\put{} at -3 0
\multiput{$\smallsq3$} at  25.85 -4.6  25.85 -4.3  25.85 -4 /
\multiput{$\boldkey.$} at 26 -4 /

\multiput{$\smallsq3$} at 29.85 -4.3  29.85 -4  29.85 -3.7 /

\multiput{$\smallsq3$} at 33.85 -4.9  33.85 -4.6  33.85 -4.3  33.85 -4
        33.85 -3.7  33.85 -3.4 /
\multiput{$\boldkey.$} at 34 -4.3 /
\multiput{$\smallsq3$} at 0 0  0 .3  0 .6  0 -.3  0 -.6  0 -.9  -.3 -.3  
       -.3 -.6 /
\multiput{$\boldkey.$} at -.15 -.3  .15 -.3 /
\plot -.15 -.3  .15 -.3 /

\multiput{$\smallsq3$} at 3.85 0  3.85 -.3  3.85 -.6 /
\multiput{$\boldkey.$} at 4 -.6 /

\multiput{$\smallsq3$} at 7.7 -.3  8 -.9  8 -.6  8 -.3  8 0  8 .3 /
\multiput{$\boldkey.$} at 7.85 -.3  8.15 -.3 /
\plot  7.85 -.3  8.15 -.3 /

\multiput{$\smallsq3$} at 11.7 -.6  11.7 -.3  11.7 0  12 -.9  12 -.6  12 -.3
        12 0  12 .3  12 .6 /
\multiput{$\boldkey.$} at 11.85 -.3  12.15 -.3 /
\plot 11.85 -.3  12.15 -.3 /

\multiput{$\smallsq3$} at 15.85 -.6  15.85 -.3  15.85 0  15.85 .3 /
\multiput{$\boldkey.$} at 16 -.6 /

\multiput{$\smallsq3$} at 19.7 -.9   19.7 -.6  19.7 -.3  19.7 0  19.7 .3
        19.7 .6  20 -.3 /
\multiput{$\boldkey.$} at 19.85 -.3  20.15 -.3 /
\plot 19.85 -.3  20.15 -.3 /

\multiput{$\smallsq3$} at 23.85 -.6  23.85 -.3  23.85 0  /
\multiput{$\boldkey.$} at 24 -.3 /

\multiput{$\smallsq3$} at 27.7 -.6  27.7 -.3  27.7 0  27.7 .3 28 -.3  28 0 /
\multiput{$\boldkey.$} at 27.85 0  28.15 0 /
\plot 27.85 0  28.15 0 /

\multiput{$\smallsq3$} at 31.7 -.3  31.7 0  31.7 .3  32 -.9  32 -.6  32 -.3 
        32 0  32 .3  32 .6 /
\multiput{$\boldkey.$} at 31.85 -.3  32.15 -.3 /
\plot 31.85 -.3  32.15 -.3 /

\multiput{$\smallsq3$} at 35.85 -.6  35.85 -.3  35.85 0  35.85 .3  35.85 .6 /
\multiput{$\boldkey.$} at 36 -.3 /

\multiput{$\smallsq3$} at 39.7 -.3  39.7 0  40 -.6  40 -.3  40 0  40 .3 
        40 .6  40 .9 /
\multiput{$\boldkey.$} at 39.85 0  40.15 0 /
\plot 39.85 0  40.15 0 /
\multiput{$\smallsq3$} at 1.55 3.4  1.55 3.7  1.85 3.1  1.85 3.4  1.85 3.7
               1.85 4  1.85 4.3  1.85 4.6  2.15 3.4  2.15 3.7  2.15 4 /
\multiput{$\boldkey.$} at 1.7 3.4  1.7 3.7  2 3.7  2.3 3.7 /
\plot 1.7 3.7  2.3 3.7 /

\multiput{$\smallsq3$} at 5.55 3.4  5.55 3.7  5.55 4  5.85 3.7  6.15 3.1 
                6.15 3.4  6.15 3.7  6.15 4  6.15 4.3 /
\multiput{$\boldkey.$} at 5.7 3.4  5.7 3.7  6 3.7  6.3 3.7 /
\plot 5.7 3.7  6.3 3.7 /

\multiput{$\smallsq3$} at  9.4 3.7  9.7 3.1  9.7 3.4  9.7 3.7  9.7 4
                9.7 4.3  10 3.4  10 3.7  10 4  10.3 3.1  10.3 3.4  10.3 3.7 
                10.3 4  10.3 4.3  10.3 4.6 /
\multiput{$\boldkey.$} at 9.5 3.65  9.8 3.65  9.9 3.75  10.2 3.75  
                10.5 3.75 / 
\plot 9.5 3.65  9.8 3.65 /
\plot 9.9 3.75  10.5 3.75 /

\multiput{$\smallsq3$} at 13.55 3.4  13.55 3.7  13.55 4  13.85 3.1  
                13.85 3.4  13.85 3.7  13.85 4  13.85 4.3  13.85 4.6
                14.15 3.4  14.15 3.7  14.15 4  14.15 4.3 /
\multiput{$\boldkey.$} at 13.7 3.4  13.7 3.7  14 3.7  14.3 3.7 /
\plot 13.7 3.7  14.3 3.7 /

\multiput{$\smallsq3$} at 17.55 3.4  17.55 3.7  17.55 4  17.55 4.3 
               17.85 3.1  17.85 3.4  17.85 3.7  17.85 4  17.85 4.3  17.85 4.6
                18.15 3.7 /
\multiput{$\boldkey.$} at 17.7 3.4  17.7 3.7  18 3.7  18.3 3.7 /
\plot  17.7 3.7  18.3 3.7 /

\multiput{$\smallsq3$} at 21.55 3.1  21.55 3.4  21.55 3.7  21.55 4  21.55 4.3
                21.55 4.6  21.85 3.7  22.15 3.4  22.15 3.7  22.15 4 /
\multiput{$\boldkey.$} at 21.65 3.65  21.95 3.65  22.05 3.75  22.35 3.75 /
\plot 21.65 3.65  21.95 3.65 /
\plot  22.05 3.75  22.35 3.75 /

\multiput{$\smallsq3$} at 25.7 3.4  25.7 3.7  25.7 4  25.7 4.3  26 3.7  26 4 /
\multiput{$\boldkey.$} at 25.85 3.7  26.15 3.7 /
\plot  25.85 3.7  26.15 3.7 /

\multiput{$\smallsq3$} at 29.55 3.4  29.55 3.7  29.55 4  29.55 4.3  29.85 3.7
                29.85 4  30.15 3.1  30.15 3.4  30.15 3.7  30.15 4  30.15 4.3
                30.15 4.6 /
\multiput{$\boldkey.$} at 29.7 4  30 4  30 3.7  30.3 3.7 /
\plot 29.7 4  30 4 /
\plot  30 3.7  30.3 3.7 /

\multiput{$\smallsq3$} at 33.7 3.7  33.7 4  33.7 4.3  34 3.4  34 3.7  34 4
                34 4.3  34 4.6 /
\multiput{$\boldkey.$} at 33.85 3.7  34.15 3.7 /
\plot 33.85 3.7  34.15 3.7 /

\multiput{$\smallsq3$} at 37.55 3.4  37.55 3.7  37.55 4  37.55 4.3  37.55 4.6
                37.85 3.7  37.85 4  38.15 3.4  38.15 3.7  38.15 4  38.15 4.3
                38.15 4.6  38.15 4.9 /
\multiput{$\boldkey.$} at 37.7 4  38 4  38 3.7  38.3 3.7 /
\plot 37.7 4  38 4 /
\plot  38 3.7  38.3 3.7 /

\multiput{$\smallsq3$} at -.6 7.4  -.6 7.7  -.3 7.4  -.3 7.7  -.3 8
        0 7.1  0 7.4  0 7.7  0 8  0 8.3  .3 7.1  .3 7.4  .3 7.7  .3 8
        .3 8.3  .3 8.6 /
\multiput{$\boldkey.$} at  -.45 7.4  -.45 7.7  -.15 7.7  .15 7.7  .15 7.4
                .45 7.4 /
\plot -.45 7.7  .15 7.7 /
\plot .15 7.4   .45 7.4 /

\multiput{$\smallsq3$} at 3.25 7.7  3.55 7.4  3.55 7.7  3.55 8  
                3.85 7.1  3.85 7.4  3.85 7.7  3.85 8  3.85 8.3
                4.15 7.1  4.15 7.4  4.15 7.7  4.15 8  4.15 8.3  4.15 8.6
                4.45 7.4  4.45 7.7 /
\multiput{$\boldkey.$} at 3.35 7.65  3.65 7.65  3.95 7.65  
                4.05 7.75  4.35 7.75  4.65 7.75  4.6 7.4 /
\plot 3.35 7.65  3.95 7.65 /
\plot 4.05 7.75  4.65 7.75 /

\multiput{$\smallsq3$} at 7.25 7.4  7.25 7.7  7.25 8  7.55 7.1  7.55 7.4
                7.55 7.7  7.55 8  7.55 8.3  7.85 7.1  7.85 7.4  7.85 7.7
                7.85 8  7.85 8.3  7.85 8.6  8.15 7.7  8.45 7.4  8.45 7.7
                8.45 8 /
\multiput{$\boldkey.$} at 7.35 7.65  7.65 7.65  7.75 7.75  8.05 7.75
                8.35 7.75  8.65 7.75  8.6 7.4 /
\plot 7.35 7.65  7.65 7.65 /
\plot 7.75 7.75  8.65 7.75 /

\multiput{$\smallsq3$} at 11.25 7.4  11.25 7.7  11.25 8  11.25 8.3
        11.55 7.1  11.55 7.4  11.55 7.7  11.55 8  11.55 8.3  11.55 8.6
        11.85 7.7  12.15 7.4  12.15 7.7  12.15 8  12.45 7.1  12.45 7.4
        12.45 7.7  12.45 8  12.45 8.3 / 
\multiput{$\boldkey.$} at 11.35 7.65  11.65 7.65  11.95 7.65  12.25 7.65
                12.35 7.75  12.65 7.75  12.3  7.4 /
\plot 11.35 7.65  12.25 7.65 /
\plot 12.35 7.75  12.65 7.75 /

\multiput{$\smallsq3$} at 15.25 7.1  15.25 7.4  15.25 7.7  15.25 8  15.25 8.3
                15.25 8.6  15.55 7.7  15.85 7.4  15.85 7.7  15.85 8  
                16.15 7.1  16.15 7.4  16.15 7.7  16.15 8  16.15 8.3  16.15 8.6
                16.45 7.4  16.45 7.7  16.45 8  16.45 8.3  /
\multiput{$\boldkey.$} at 15.35 7.65  15.65 7.65  15.95 7.65 
                16.05 7.75  16.35 7.75  16.65 7.75  16 7.4 /
\plot 15.35 7.65  15.95 7.65 /
\plot 16.05 7.75  16.65 7.75 /

\multiput{$\smallsq3$} at 19.4 7.7  19.7 7.4  19.7 7.7  19.7 8  20 7.1
                20 7.4  20 7.7  20 8  20 8.3  20 8.6  20.3 7.4  20.3 7.7
                20.3 8  20.3 8.3 /
\multiput{$\boldkey.$} at 19.5 7.65  19.8 7.65  19.9 7.75  20.2 7.75  20.5 7.75
                20.45 7.4 /
\plot 19.5 7.65  19.8 7.65 /
\plot 19.9 7.75  20.5 7.75 /

\multiput{$\smallsq3$} at 23.4 7.7  23.4 8  23.7 7.1  23.7 7.4  23.7 7.7
                23.7 8  23.7 8.3  23.7 8.6  24 7.4  24 7.7  24 8  24 8.3 
                24.3 7.7 /
\multiput{$\boldkey.$} at 23.5 7.65  23.8 7.65  24.1 7.65  24.2 7.75  
                24.5 7.75 /
\plot 23.5 7.65  24.1 7.65 /
\plot 24.2 7.75  24.5 7.75 /

\multiput{$\smallsq3$} at 27.55 7.1  27.55 7.4  27.55 7.7  27.55 8  27.55 8.3
        27.55 8.6  27.85 7.4  27.85 7.7  27.85 8  27.85 8.3  28.15 7.7  28.15 8         /
\multiput{$\boldkey.$} at 27.65 7.65  27.95 7.65  28.05 7.75  28.35 7.75 /
\plot 27.65 7.65  27.95 7.65 /
\plot   28.05 7.75  28.35 7.75 /

\multiput{$\smallsq3$} at 31.55 7.4  31.55 7.7  31.55 8  31.55 8.3  31.55 8.6
                31.85 7.7  31.85 8  32.15 7.4  32.15 7.7  32.15 8  32.15 8.3 
                /
\multiput{$\boldkey.$} at 31.7 7.7  32 7.7  32 8 32.3 8 /
\plot 31.7 7.7  32 7.7 /
\plot  32 8 32.3 8 /

\multiput{$\smallsq3$} at   35.4 7.4  35.4 7.7  35.4 8  35.4 8.3  35.4 8.6  
                35.4 8.9  35.7 7.7  35.7 8  36 7.7  36 8  36 8.3
                36.3 7.4  36.3 7.7  36.3 8  36.3 8.3  36.3 8.6 /
\multiput{$\boldkey.$} at 35.55 7.7  35.85 7.7  35.85 8  36.15 8  36.45 8 /
\plot 35.55 7.7  35.85 7.7  /
\plot 35.85 8  36.45 8 /

\multiput{$\smallsq3$} at  39.4 7.7  39.4 8  39.7 7.7  39.7 8  39.7 8.3 
                40 7.4  40 7.7  40 8  40 8.3  40 8.6  
                40.3 7.4  40.3 7.7  40.3 8  40.3 8.3  40.3 8.6  40.3 8.9 /
\multiput{$\boldkey.$} at 39.55 7.7  39.55 8  39.85 8  40.15 8  
                40.15 7.7  40.45 7.7 /
\plot 39.55 8  40.15 8 /
\plot 40.15 7.7  40.45 7.7 /

\multiput{$\smallsq3$} at 1.7 6.1  1.7 6.4  1.7 6.7  1.7 7  1.7 7.3 
        2 6.4  2 6.7 /
\multiput{$\boldkey.$} at 1.85 6.7  2.15 6.7  2.15 6.4 /
\plot 1.85 6.7  2.15 6.7 /

\multiput{$\smallsq3$} at 5.55 6.1  5.55 6.4  5.55 6.7  5.55 7  5.55 7.3
        5.55 7.6  5.85 6.7  6.15 6.4  6.15 6.7  6.15 7 /
\multiput{$\boldkey.$} at 5.7 6.7  6 6.7  6.3 6.7 /
\plot 5.7 6.7  6.3 6.7 /

\multiput{$\smallsq3$} at 9.7 6.4  9.7 6.7  9.7 7  10 6.1  10 6.4  10 6.7
                10 7 10 7.3 /
\multiput{$\boldkey.$} at 9.85 6.7  10.15 6.7  9.85 6.4 /
\plot 9.85 6.7  10.15 6.7 /

\multiput{$\smallsq3$} at 13.55 6.1  13.55 6.4  13.55 6.7  13.55 7  13.55 7.3 
                13.55 7.6  13.85 6.4  13.85 6.7  13.85 7  13.85 7.3  14.15 6.7
                /
\multiput{$\boldkey.$} at 13.7 6.7  14 6.7  14.3 6.7 /
\plot 13.7 6.7  14.3 6.7 /

\multiput{$\smallsq3$} at 17.7 6.1  17.7 6.4  17.7 6.7  17.7 7  17.7 7.3 
                17.7 7.6  18 6.4  18 6.7  18 7 /
\multiput{$\boldkey.$} at 17.85 6.7  18.15 6.7  18.15 6.4 /
\plot 17.85 6.7  18.15 6.7 /

\multiput{$\smallsq3$} at 21.7 6.4  21.7 6.7  21.7 7  21.7 7.3  22 6.7 /
\multiput{$\boldkey.$} at 21.85 6.7  22.15 6.7 /
\plot 21.85 6.7  22.15 6.7 /

\multiput{$\smallsq3$} at 25.7 6.1  25.7 6.4  25.7 6.7  25.7 7  25.7 7.3 
        25.7 7.6  26 6.7  26 7 /
\multiput{$\boldkey.$} at 25.85 6.7  26.15 6.7 /
\plot  25.85 6.7  26.15 6.7 /

\multiput{$\smallsq3$} at 29.85 6.4  29.85 6.7  29.85 7  29.85 7.3 /
\multiput{$\boldkey.$} at 30 6.7 /

\multiput{$\smallsq3$} at 33.7 6.4  33.7 6.7  33.7 7  33.7 7.3  33.7 7.6 
                34 6.7  34 7 /
\multiput{$\boldkey.$} at 33.85 7  34.15 7 /
\plot 33.85 7  34.15 7 /

\multiput{$\smallsq3$} at 37.7 6.7  37.7 7  37.7 7.3  38 6.4  38 6.7
                38 7  38 7.3  38 7.6  38 7.9 /
\multiput{$\boldkey.$} at 37.85 6.7  38.15 6.7 /
\plot 37.85 6.7  38.15 6.7 /

\multiput{$\smallsq3$} at 1.4 10.7  1.7 10.4  1.7 10.7  1.7 11  2 10.1
                2 10.4  2 10.7  2 11  2 11.3  2.3 10.1  2.3 10.4  2.3 10.7
                2.3 11  2.3 11.3  2.3 11.6 /
\multiput{$\boldkey.$} at 1.55 10.7  1.85 10.7  2.15 10.7  2.15 10.4  2.45 10.4
                /
\plot 1.55 10.7  2.15 10.7 /
\plot 2.15 10.4  2.45 10.4 /

\multiput{$\smallsq3$} at 5.4 10.4  5.4 10.7  5.4 11  5.7 10.1  5.7 10.4
                5.7 10.7  5.7 11  5.7 11.3  6 10.1  6 10.4  6 10.7  6 11
                6 11.3  6 11.6  6.3 10.4  6.3 10.7 /
\multiput{$\boldkey.$} at 5.5 10.65  5.8 10.65  5.9 10.75  6.2 10.75
                6.5 10.75  6.45 10.4 /
\plot 5.5 10.65  5.8 10.65 /
\plot 5.9 10.75  6.5 10.75 /

\multiput{$\smallsq3$} at 9.4 10.4  9.4 10.7  9.4 11  9.4 11.3  9.7 10.1
                9.7 10.4  9.7 10.7  9.7 11 9.7 11.3  9.7 11.6  10 10.7
                10.3 10.4  10.3 10.7  10.3 11 /
\multiput{$\boldkey.$} at 9.55 10.7  9.85 10.7  10.15 10.7  10.45 10.7
                10.45 10.4 /
\plot 9.55 10.7  10.45 10.7 /

\multiput{$\smallsq3$} at 13.4 10.1  13.4 10.4  13.4 10.7  13.4 11  13.4 11.3
                13.4 11.6  13.7 10.7  14 10.4  14 10.7  14 11  14.3  10.1
                14.3 10.4  14.3 10.7  14.3 11  14.3 11.3 /
\multiput{$\boldkey.$} at 13.5 10.65  13.8 10.65  14.1 10.65  14.2 10.75
                14.5 10.75  14.15 10.4 /
\plot 13.5 10.65  14.1 10.65 /
\plot 14.2 10.75  14.5 10.75 /

\multiput{$\smallsq3$} at 17.4 10.7  17.7 10.4  17.7 10.7  17.7 11  
        18 10.1  18 10.4  18 10.7  18 11  18 11.3  18 11.6  
        18.3 10.4  18.3 10.7  18.3 11  18.3 11.3 /
\multiput{$\boldkey.$} at 17.5 10.65  17.8 10.65  17.9 10.75  18.2 10.75
        18.5 10.75 /
\plot 17.5 10.65  17.8 10.65 /
\plot 17.9 10.75  18.5 10.75 /

\multiput{$\smallsq3$} at 21.55 10.7  21.55 11  21.85 10.1  21.85 10.4
        21.85 10.7  21.85 11  21.85 11.3  21.85 11.6  
        22.15 10.4  22.15 10.7  22.15 11  22.15 11.3 /
\multiput{$\boldkey.$} at 21.7 10.7  22 10.7  22.3 10.7  22.3 10.4 /
\plot 21.7 10.7  22.3 10.7 /

\multiput{$\smallsq3$} at 25.55 10.1  25.55 10.4  25.55 10.7  25.55 11 
                25.55 11.3  25.55 11.6  25.85 10.4  25.85 10.7  25.85 11
                25.85 11.3  26.15 10.7 / 
\multiput{$\boldkey.$} at 25.65 10.65  25.95 10.65  26.05 10.75  26.35 10.75 /
\plot 26.05 10.75  26.35 10.75 /
\plot 25.65 10.65  25.95 10.65 /

\multiput{$\smallsq3$} at 29.7 10.4  29.7 10.7  29.7 11 29.7 11.3  29.7 11.6 
                30 10.7  30 11 /
\multiput{$\boldkey.$} at 29.85 10.7  30.15 10.7 /
\plot 29.85 10.7  30.15 10.7 /

\multiput{$\smallsq3$} at 33.55 10.4  33.55 10.7  33.55 11  33.55 11.3  
        33.55 11.6  33.55 11.9  33.85 10.7  33.85 11  34.15 10.4  34.15 10.7
        34.15 11  34.15 11.3 /
\multiput{$\boldkey.$} at 33.7 10.7  34 10.7  34 11  34.3 11 /
\plot 34 11  34.3 11 /
\plot  33.7 10.7  34 10.7 /

\multiput{$\smallsq3$} at 37.55 10.7  37.55 11  37.85 10.7  37.85 11  
        37.85 11.3  38.15 10.4  38.15 10.7  38.15 11  38.15 11.3  38.15 11.6 /
\multiput{$\boldkey.$} at 37.7 10.7  37.7 11  38 11  38.3 11 /
\plot 37.7 11  38.3 11 /

\multiput{$\smallsq3$} at -.45 13.7  -.15 13.4  -.15 13.7
        -.15 14  .15 13.1  .15 13.4  .15 13.7  .15 14  .15 14.3 /
\multiput{$\boldkey.$} at -.3 13.7  0 13.7  .3 13.7 /
\plot -.3 13.7  .3 13.7 /

\multiput{$\smallsq3$} at 3.55 13.4  3.55 13.7  3.55 14  3.85 13.1  
        3.85 13.4  3.85 13.7  3.85 14  3.85 14.3  4.15 13.1  4.15 13.4 
        4.15 13.7  4.15 14  4.15 14.3  4.15 14.6 /
\multiput{$\boldkey.$} at 3.7 13.7  4 13.7  4 13.4  4.3 13.4 /
\plot 4 13.4  4.3 13.4 /
\plot 3.7 13.7  4 13.7 /

\multiput{$\smallsq3$} at 7.55 13.4  7.55 13.7  7.55 14  7.55 14.3  7.85 13.1
        7.85 13.4  7.85 13.7  7.85 14  7.85 14.3  7.85 14.6 
        8.15 13.4  8.15 13.7 /
\multiput{$\boldkey.$} at 7.7 13.7  8 13.7  8.3 13.7  8.3 13.4 /
\plot  7.7 13.7  8.3 13.7 /

\multiput{$\smallsq3$} at 11.55 13.1  11.55 13.4  11.55 13.7  11.55 14 
        11.55 14.3  11.55 14.6  11.85 13.7  12.15 13.4  12.15 13.7  12.15 14 /
\multiput{$\boldkey.$} at 11.7 13.7  12 13.7  12.3 13.7  12.3 13.4 /
\plot 11.7 13.7  12.3 13.7 /

\multiput{$\smallsq3$} at 15.55 13.7  15.85 13.4  15.85 13.7  15.85 14
                16.15 13.1  16.15 13.4  16.15 13.7  16.15 14  16.15 14.3 /
\multiput{$\boldkey.$} at 15.65 13.65  15.95 13.65  16.05 13.75  16.35 13.75 /
\plot  16.05 13.75  16.35 13.75 /
\plot  15.65 13.65  15.95 13.65 /

\multiput{$\smallsq3$} at 19.55 13.7  19.55 14  19.85 13.1  19.85 13.4  
        19.85 13.7  19.85 14  19.85 14.3  19.85 14.6  20.15 13.4  20.15 13.7 
        20.15 14  20.15 14.3 /
\multiput{$\boldkey.$} at 19.7 13.7  20 13.7  20.3 13.7 /
\plot 19.7 13.7  20.3 13.7 /

\multiput{$\smallsq3$} at 23.7 13.1  23.7 13.4  23.7 13.7  23.7 14 
        23.7 14.3  23.7 14.6  24 13.4  24 13.7  24 14  24 14.3 /
\multiput{$\boldkey.$} at 23.85 13.7  24.15 13.7  24.15 13.4 /
\plot 23.85 13.7  24.15 13.7 /

\multiput{$\smallsq3$} at 27.7 13.4  27.7 13.7  27.7 14  27.7 14.3  27.7 14.6
        28 13.7 /
\multiput{$\boldkey.$} at 27.85 13.7  28.15 13.7 /
\plot 27.85 13.7  28.15 13.7 /

\multiput{$\smallsq3$} at 31.7 13.4  31.7 13.7  31.7 14  31.7 14.3  31.7 14.6
        31.7 14.9  32 13.7  32 14 /
\multiput{$\boldkey.$} at 31.85 13.7  32.15 13.7 /
\plot 31.85 13.7  32.15 13.7 /

\multiput{$\smallsq3$} at 35.7 13.7  35.7 14  36 13.4  36 13.7  36 14  36 14.3 
                        /
\multiput{$\boldkey.$} at 35.85 13.7  35.85 14  36.15 14 /
\plot 35.85 14  36.15 14 /

\multiput{$\smallsq3$} at 39.55 14  39.85 13.7  39.85 14  
        39.85 14.3  40.15 13.4  40.15 13.7  40.15 14  40.15 14.3  40.15 14.6 /
\multiput{$\boldkey.$} at 39.7 14  40 14  40.3 14 /
\plot 39.7 14  40.3 14 /

\multiput{$\smallsq3$} at 1.7 16.4  1.7 16.7  1.7 17
                2 16.1  2 16.4  2 16.7  2 17  2 17.3 /
\multiput{$\boldkey.$} at 1.85 16.7  2.15 16.7 /
\plot  1.85 16.7  2.15 16.7 /

\multiput{$\smallsq3$} at 5.7 16.4  5.7 16.7  5.7 17  5.7 17.3  6 16.1
        6 16.4  6 16.7  6 17  6 17.3  6 17.6 /
\multiput{$\boldkey.$} at 5.85 16.4  6.15 16.4 /
\plot 5.85 16.4  6.15 16.4 /

\multiput{$\smallsq3$} at 9.7 16.1  9.7 16.4  9.7 16.7  9.7 17  9.7 17.3
        9.7 17.6  10 16.4  10 16.7 /
\multiput{$\boldkey.$} at 9.85 16.7  10.15 16.7  10.15 16.4 /
\plot 9.85 16.7  10.15 16.7 /

\multiput{$\smallsq3$} at 13.7 16.7  14 16.4  14 16.7  14 17 /
\multiput{$\boldkey.$} at 13.85 16.7  14.15 16.7 /
\plot 13.85 16.7  14.15 16.7 /

\multiput{$\smallsq3$} at 17.7 16.7  17.7 17  18 16.1  18 16.4  18 16.7
        18 17  18 17.3 /
\multiput{$\boldkey.$} at 17.85 16.7  18.15 16.7 /
\plot 17.85 16.7  18.15 16.7 /

\multiput{$\smallsq3$} at 21.7 16.1  21.7 16.4  21.7 16.7  21.7 17  21.7 17.3
        21.7 17.6  22 16.4  22 16.7  22 17  22 17.3 /
\multiput{$\boldkey.$} at 21.85 16.7  22.15 16.7 /
\plot 21.85 16.7  22.15 16.7 /

\multiput{$\smallsq3$} at 25.85 16.4  25.85 16.7  25.85 17  25.85 17.3
        25.85 17.6 /
\multiput{$\boldkey.$} at 26 16.4 /

\multiput{$\smallsq3$} at 29.7 16.4  29.7 16.7  29.7 17  29.7 17.3  29.7 17.6
        29.7 17.9  30 16.7 /
\multiput{$\boldkey.$} at 29.85 16.7  30.15 16.7 /
\plot  29.85 16.7  30.15 16.7 /

\multiput{$\smallsq3$} at 33.85 16.7  33.85 17 /
\multiput{$\boldkey.$} at 34 16.7 /

\multiput{$\smallsq3$} at 37.7 17  38 16.4  38 16.7  38 17  38 17.3 /
\multiput{$\boldkey.$} at 37.85 17  38.15 17 /
\plot 37.85 17  38.15 17 /

\multiput{$\smallsq3$} at -.15 19.1  -.15 19.4  -.15 19.7  -.15 20 /
\multiput{$\boldkey.$} at 0 19.7 /

\multiput{$\smallsq3$} at 3.85 19.4  3.85 19.7  3.85 20  3.85 20.3 /

\multiput{$\smallsq3$} at 7.85 19.1  7.85 19.4  7.85 19.7  7.85 20 
        7.85 20.3  7.85 20.6 /
\multiput{$\boldkey.$} at 8 19.4 /

\multiput{$\smallsq3$} at 11.85 19.4  11.85 19.7 /
\multiput{$\boldkey.$} at 12 19.7 /

\multiput{$\smallsq3$} at 15.85 19.7  15.85 20 /

\multiput{$\smallsq3$} at 19.85 19.1  19.85 19.4  19.85 19.7  19.85 20
        19.85 20.3  /
\multiput{$\boldkey.$} at 20 19.7 /

\multiput{$\smallsq3$} at 23.85 19.4  23.85 19.7  23.85 20  23.85 20.3
        23.85 20.6 / 

\multiput{$\smallsq3$} at 27.85 19.4  27.85 19.7  27.85 20  27.85 20.3
        27.85 20.6  27.85 20.9 /
\multiput{$\boldkey.$} at 28 19.4 /

\multiput{$\smallsq3$} at 31.85 19.7 /
\multiput{$\boldkey.$} at 32 19.7 /
 
\multiput{$\smallsq3$} at 35.85 20 /

\multiput{$\smallsq3$} at 39.85 19.4  39.85 19.7  39.85 20  39.85 20.3 /
\multiput{$\boldkey.$} at 40 20 /
\multiput{$\smallsq3$} at 25.85 22.4  25.85 22.7  25.85 23  25.85 23.3
        25.85 23.6  25.85 23.9 /
\arr{24.5 -1} {25.5 -3}
\arr{28.5 -1} {29.5 -3}
\arr{32.5 -1} {33.5 -3}
\arr{26.5 -3} {27.5 -1}
\arr{30.5 -3} {31.5 -1}
\arr{34.5 -3} {35.5 -1}
\arr{0.5 1} {1.25 2.5}
\arr{4.5 1} {5.25 2.5}
\arr{8.5 1} {9.25 2.5}
\arr{12.5 1} {13.25 2.5}
\arr{16.5 1} {17.25 2.5}
\arr{20.5 1} {21.25 2.5}
\arr{24.5 1} {25.25 2.5}
\arr{28.5 1} {29.25 2.5}
\arr{32.5 1} {33.25 2.5}
\arr{36.5 1} {37.25 2.5}
\arr{2.75 2.5} {3.5 1}
\arr{6.75 2.5} {7.5 1}
\arr{10.75 2.5} {11.5 1}
\arr{14.75 2.5} {15.5 1}
\arr{18.75 2.5} {19.5 1}
\arr{22.75 2.5} {23.5 1}
\arr{26.75 2.5} {27.5 1}
\arr{30.75 2.5} {31.5 1}
\arr{34.75 2.5} {35.5 1}
\arr{38.75 2.5} {39.5 1}
\arr{0.75 6.5} {1.5 5}
\arr{4.75 6.5} {5.5 5}
\arr{8.75 6.5} {9.5 5}
\arr{12.75 6.5} {13.5 5}
\arr{16.75 6.5} {17.5 5}
\arr{20.75 6.5} {21.5 5}
\arr{24.75 6.5} {25.5 5}
\arr{28.75 6.5} {29.5 5}
\arr{32.75 6.5} {33.5 5}
\arr{36.75 6.5} {37.5 5}
\arr{2.5 5} {3.25 6.5}
\arr{6.5 5} {7.25 6.5}
\arr{10.5 5} {11.25 6.5}
\arr{14.5 5} {15.25 6.5}
\arr{18.5 5} {19.25 6.5}
\arr{22.5 5} {23.25 6.5}
\arr{26.5 5} {27.25 6.5}
\arr{30.5 5} {31.25 6.5}
\arr{34.5 5} {35.25 6.5}
\arr{38.6 5.2} {39.25 6.5}
\arr{1 7.5} {1.5 7.25}
\arr{5 7.5} {5.4 7.3}
\arr{9 7.5} {9.5 7.25}
\arr{13 7.5} {13.4 7.3}
\arr{17 7.5} {17.5 7.25}
\arr{21 7.5} {21.5 7.25}
\arr{25 7.5} {25.5 7.25}
\arr{29 7.5} {29.5 7.25}
\arr{33 7.5} {33.5 7.25}
\arr{37 7.5} {37.5 7.25}
\arr{2.5 7.25} {3 7.5}
\arr{6.5 7.25} {7 7.5}
\arr{10.5 7.25} {11 7.5}
\arr{14.5 7.25} {15 7.5}
\arr{18.5 7.25} {19 7.5}
\arr{22.5 7.25} {23 7.5}
\arr{26.5 7.25} {27 7.5}
\arr{30.5 7.25} {31 7.5}
\arr{34.5 7.25} {35 7.5}
\arr{38.5 7.25} {39 7.5}
\arr{0.67 9} {1.33 10}
\arr{4.67 9} {5.33 10}
\arr{8.67 9} {9.33 10}
\arr{12.67 9} {13.23 9.85}
\arr{16.67 9} {17.33 10}
\arr{20.67 9} {21.33 10}
\arr{24.67 9} {25.33 10}
\arr{28.67 9} {29.33 10}
\arr{32.67 9} {33.33 10}
\arr{36.67 9} {37.33 10}
\arr{2.77 9.85} {3.33 9}
\arr{6.67 10} {7.33 9}
\arr{10.67 10} {11.33 9}
\arr{14.77 9.85} {15.33 9}
\arr{18.67 10} {19.33 9}
\arr{22.67 10} {23.33 9}
\arr{26.67 10} {27.33 9}
\arr{30.67 10} {31.33 9}
\arr{34.67 10} {35.23 9.15}
\arr{38.67 10} {39.33 9}
\arr{0.67 13} {1.33 12}
\arr{4.67 13} {5.33 12}
\arr{8.67 13} {9.33 12}
\arr{12.67 13} {13.33 12}
\arr{16.67 13} {17.33 12}
\arr{20.67 13} {21.33 12}
\arr{24.67 13} {25.33 12}
\arr{28.67 13} {29.33 12}
\arr{32.67 13} {33.33 12}
\arr{36.67 13} {37.33 12}
\arr{2.67 12} {3.33 13}
\arr{6.67 12} {7.33 13}
\arr{10.67 12} {11.33 13}
\arr{14.67 12} {15.33 13}
\arr{18.67 12} {19.33 13}
\arr{22.67 12} {23.33 13}
\arr{26.67 12} {27.33 13}
\arr{30.67 12} {31.33 13}
\arr{34.67 12} {35.33 13}
\arr{38.67 12} {39.33 13}
\arr{0.67 15} {1.33 16}
\arr{4.67 15} {5.33 16}
\arr{8.67 15} {9.33 16}
\arr{12.67 15} {13.33 16}
\arr{16.67 15} {17.33 16}
\arr{20.67 15} {21.33 16}
\arr{24.67 15} {25.33 16}
\arr{28.67 15} {29.33 16}
\arr{32.67 15} {33.33 16}
\arr{36.67 15} {37.33 16}
\arr{2.67 16} {3.33 15}
\arr{6.67 16} {7.33 15}
\arr{10.67 16} {11.33 15}
\arr{14.67 16} {15.33 15}
\arr{18.67 16} {19.33 15}
\arr{22.67 16} {23.33 15}
\arr{26.67 16} {27.33 15}
\arr{30.67 16} {31.33 15}
\arr{34.67 16} {35.33 15}
\arr{38.67 16} {39.33 15}
\arr{0.5 19.25} {1.33 18}
\arr{4.5 19.25} {5.33 18}
\arr{8.5 19.25} {9.33 18}
\arr{12.5 19.25} {13.33 18}
\arr{16.5 19.25} {17.33 18}
\arr{20.5 19.25} {21.33 18}
\arr{24.5 19.25} {25.33 18}
\arr{28.5 19.25} {29.33 18}
\arr{32.5 19.25} {33.33 18}
\arr{36.5 19.25} {37.33 18}
\arr{2.67 18} {3.5 19.25}
\arr{6.67 18} {7.5 19.25}
\arr{10.67 18} {11.5 19.25}
\arr{14.67 18} {15.5 19.25}
\arr{18.67 18} {19.5 19.25}
\arr{22.67 18} {23.5 19.25}
\arr{26.67 18} {27.5 19.25}
\arr{30.67 18} {31.5 19.25}
\arr{34.67 18} {35.5 19.25}
\arr{38.67 18} {39.5 19.25}
\arr{24.5 20.75} {25.5 22.25}
\arr{26.5 22.25} {27.5 20.75}
\setdots<2pt>
\plot 0.5 -.15  3.5 -.15 /
\plot 4.5 -.15  7.5 -.15 /
\plot 8.5 -.15  11.5 -.15 /
\plot 12.5 -.15  15.5 -.15 /
\plot 16.5 -.15  19.5 -.15 /
\plot 20.5 -.15  23.5 -.15 /
\plot 26.5 -4.15  29.5 -4.15 /
\plot 30.5 -4.15  33.5 -4.15 /
\plot 36.5 -.15  39.5 -.15 /
\plot 0.5 19.85  3.5 19.85 /
\plot 4.5 19.85  7.5 19.85 /
\plot 8.5 19.85  11.5 19.85 /
\plot 12.5 19.85  15.5 19.85 /
\plot 16.5 19.85  19.5 19.85 /
\plot 20.5 19.85  23.5 19.85 /
\plot 28.5 19.85  31.5 19.85 /
\plot 32.5 19.85  35.5 19.85 /
\plot 36.5 19.85  39.5 19.85 /

\setsolid 
\plot 0 1.2  0 6.6 /
\plot 0 9  0 12.6 /
\plot 0 15  0 18.6 /
\plot 40 1.5  40 6.9 /
\plot 40 9.3  40 12.9 /
\plot 40 15.3  40 18.9 /
\setshadegrid span <1.5mm>
\vshade   0   -.15 19.85  <,z,,> 
          24  -.15 19.85  <z,z,,>
          26 -4    22.85  <z,z,,>
          28 -4    19.85  <z,z,,>
          34 -4    19.85  <z,z,,>
          36  -.15 19.85  <z,z,,>
          40  -.15 19.85 /
\endpicture}

\topinsert\noindent
\rotatebox{90}{\arqthreesix}
\endinsert

\medskip\noindent {\bf Proposition 12.} {\it 
The category $\Cal S_3(k[T]/T^6)$ has 84 indecomposable
objects and hence is the largest among the
categories of type $\Cal S_m(k[T]/T^n)$ which 
are representation finite.
The stable part of the Auslander-Reiten quiver
has type  $\Bbb Z\Bbb E_8/\tau^{10}$. There is a
nonstable orbit of length 3 from $P$ to $I$ attached
to the orbit for the endpoint of the branch in 
the Dynkin diagram $\Bbb E_8$ of length 2; and the
module $Y$ is attached to the orbit for the endpoint
of the branch of length 4.}

\medskip
Two indecomposable objects appear to be of particular interest 
as they are the only ones (in any of the representation finite categories
of type $\Cal S_m(k[T]/T^n)$) for which the 
submodule is a free module of rank 2. 
Let $X_1$ be the module on the intersection of the ray starting at $P$
and the coray ending at $I$, and put $X_2=\tau^5X_1=\tau^{-5}X_1$.
The two modules are
pictured as follows.
$$
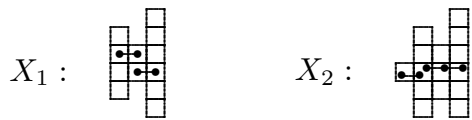
\beginpicture 
\setcoordinatesystem units <.2cm,.2cm>
\put{$X_1:$} at -4 2
\multiput{\sq} at  0 4  0 3  0 2  0 1   1 3  1 2 
                2 5  2 4  2 3  2 2  2 1  2 0 /
\multiput{$\sssize\bullet$} at 0.5 3  1.5 3  1.5 2  2.5 2 /
\plot 0.5 3  1.5 3 /
\plot 1.5 2  2.5 2 /
\endpicture \qquad\qquad
\beginpicture 
\setcoordinatesystem units <.2cm,.2cm>
\put{$X_2:$} at -5 2
\multiput{\sq} at -1  2  0 4  0 3  0 2  0 1  0 0    1 3  1 2  1 1 
               2 5   2 4  2 3  2 2  2 1  2 0 /
\multiput{$\sssize\bullet$} at -.7 1.8  .3 1.8  .7 2.2  1.7 2.2  2.7 2.2 /
\plot -.7 1.8    .3 1.8 /
\plot .7 2.2  2.7 2.2 /
\endpicture$$
Both $X_1$ and $X_2$ have as submodule two copies of $k[T]/T^3$.
In $X_1$ the copies are shifted against each other, in $X_2$ they are not.
In fact, $\Cal S_3(k[T]/T^6)$ is the only representation finite category in
which the indecomposable summands of the submodule are not all pairwise
nonisomorphic. 
Besides $X_1$ and $X_2$, there are six more such indecomposables in
$\Cal S_3(k[T]/T^6)$ where the subgroup is isomorphic to $2\cdot
k[T]/T^3\oplus k[T]/T$; note that this is the subgroup of the modules in the
radical family in $\Cal S_3(k[T]/T^7)$. 
Concerning the total spaces, 
there are three indecomposable modules for which
the total space has isomorphic indecomposable summands, and which live 
in a representation finite category of type $\Cal S_m(k[T]/T^n)$. Two of them,
$Z_1$ and $Z_2$, have two isomorphic summands 
of type $k[T]/T^3$ and $k[T]/T^6$, respectively, which 
are not shifted against each other.
Both occur on the central orbit in $\Cal S_3(k[T]/T^6)$. The third one, $Z_3$
from $\Cal S_4(k[T]/T^5)$, 
has two summands of type $k[T]/T^5$ that are shifted against each other.

$$
\beginpicture 
\setcoordinatesystem units <.2cm,.2cm>
\put{$Z_1\:$} at -5 2
\multiput{\sq} at -1 1  -1 2  -1 3  0 4  0 3  0 2  0 1  0 0   1 5  1 4 1 3  1 2
        1 1  1 0   2 2 3 1  3 2  3 3 /
\multiput{$\sssize\bullet$} at -.7 1.8  .3 1.8  .7 2.2  1.7 2.2  2.7 2.2  
                3.7 2.2  3.5 1 /
\plot -.7 1.8  .3 1.8 /
\plot .7 2.2  3.7 2.2 /
\endpicture
 \qquad\qquad
\beginpicture 
\setcoordinatesystem units <.2cm,.2cm>
\put{$Z_2\:$} at -4 2
\multiput{\sq} at  0 5  0 4  0 3  0 2  0 1  0 0   1 2
           2 3  2 2  2 1  3 5  3 4  3 3  3 2  3 1  3 0  4 4  4 3  4 2  4 1 /
\multiput{$\sssize\bullet$} at .3 1.8  1.3 1.8  2.3 1.8  2.7 2.2  3.7 2.2
                4.7 2.2  2.5 1 /
\plot .3 1.8  2.3 1.8 /
\plot 2.7 2.2  4.7 2.2 /
\endpicture
 \qquad\qquad
\beginpicture 
\setcoordinatesystem units <.2cm,.2cm>
\put{$Z_3\:$} at -4 2
\multiput{\sq} at  0 4  0 3  0 2  0 1  0 0    1 3  1 2
               2 5   2 4  2 3  2 2  2 1 /
\multiput{$\sssize\bullet$} at .5 3  1.5 3  1.5 2  2.5 2 /
\plot .5 3  1.5 3 /
\plot 1.5 2  2.5 2 /
\endpicture$$
\bigskip
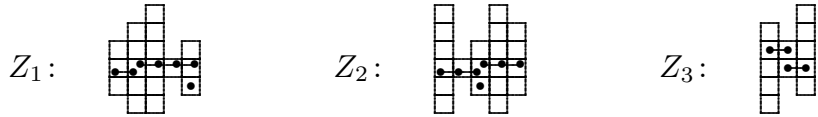
\centerline{\sc The Case $m=3$, $n=7$.}

\medskip
The category $\Cal S_3(k[T]/T^7)$ is the first of two categories
of tame infinite representation type considered in detail in this 
manuscript.  The indecomposable objects occur either in
a stable family of tubes of tubular type $(5,3,2)$ or in a 
connecting component. 
We use the indecomposables of dimension type given by the 
radical vector of the tubular family
to obtain a second minimal family of subgroup embeddings, as indicated 
in the introduction.

\smallskip 
First consider the covering category $\Cal S_3(7)$. 
Let $B_0$ be the algebra given by the following quiver
and the commutativity relation as indicated.
$$
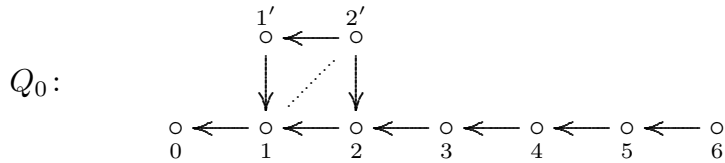
\hbox{\beginpicture
\setcoordinatesystem units <0.5cm,0.5cm>
\put{} at 0 -1.55
\put{} at 16 3.05
\put{$Q_0\:$} at -3 1
\put{$\circ$} at 0 0
\put{$\circ$} at 2 0
\put{$\circ$} at 4 0
\put{$\circ$} at 6 0
\put{$\circ$} at 8 0
\put{$\circ$} at 10 0
\put{$\circ$} at 12 0

\put{$\circ$} at 2 2
\put{$\circ$} at 4 2
\arr{1.6 0}{0.4 0}
\arr{3.6 0}{2.4 0}
\arr{5.6 0}{4.4 0}
\arr{7.6 0}{6.4 0}
\arr{9.6 0}{8.4 0}
\arr{11.6 0}{10.4 0}

\arr{3.6 2}{2.4 2}

\arr{2 1.6}{2 0.4}
\arr{4 1.6}{4 0.4}


\put{$\ssize 1'$} at  2 2.5
\put{$\ssize 2'$} at  4 2.5

\put{$\ssize 0$} at 0 -.5
\put{$\ssize 1$} at 2 -.5
\put{$\ssize 2$} at 4 -.5
\put{$\ssize 3$} at 6 -.5
\put{$\ssize 4$} at 8 -.5
\put{$\ssize 5$} at 10 -.5
\put{$\ssize 6$} at 12 -.5

\setdots<2pt>
\plot 3.5 1.5  2.5 0.5 /
\endpicture}
$$
This algebra is tame concealed of type $\tilde{\Bbb E}_8$, according to
the list of Happel and Vossiek [5, A.2], and hence the Auslander-Reiten
quiver consists of a preprojective component $\Cal P$, a family $\Cal T_0$
of tubes, and  a preinjective component $\Cal I_0$.  
With three exceptions, the tubes are homogeneous and contain representations
of dimension type a multiple of 
$\ssize{ 2320000\atop 2454321}$.
For later use 
we picture the mouth for each of the remaining three big tubes of 
circumference 2, 3 and 5.
$$
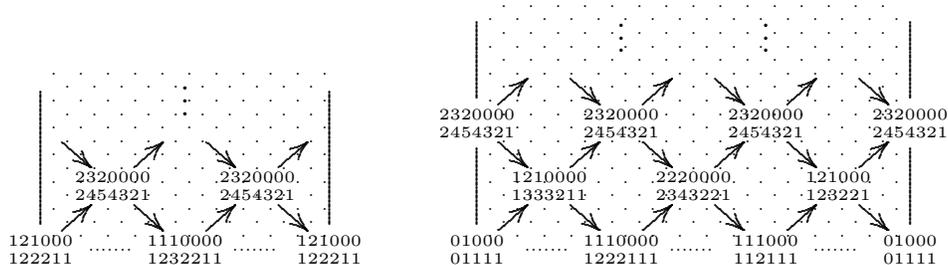
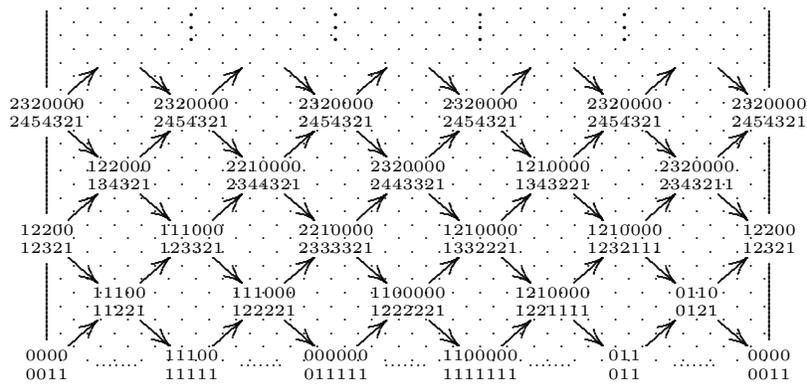
\hbox{\beginpicture
\setcoordinatesystem units <0.8cm,0.7cm>
\put{} at 0 0
\put{} at 4 2.5
\put{$\vdots$} at 2 2.5
\put{$\ssize {121000    \atop 122211}$} at 0 0
\put{$\ssize {2320000   \atop 2454321}$} at 1 1
\put{$\ssize {1110000   \atop 1232211}$} at 2 0
\put{$\ssize {2320000   \atop 2454321}$} at 3 1
\put{$\ssize {121000    \atop 122211}$} at 4 0

\arr{0.3 0.3} {0.7 0.7}
\arr{1.3 0.7} {1.7 0.3}
\arr{2.3 0.3} {2.7 0.7}
\arr{3.3 0.7} {3.7 0.3}
\arr{0.3 1.7} {0.7 1.3}
\arr{1.3 1.3} {1.7 1.7}
\arr{2.3 1.7} {2.7 1.3}
\arr{3.3 1.3} {3.7 1.7}

\setdots<2pt>
\plot 0.7 0  1.3 0 /
\plot 2.7 0  3.3 0 /

\setsolid
\plot 0 0.4  0 2.5 /
\plot 4 0.4  4 2.5 /
\setshadegrid span <1.5mm>
\vshade   0   0 3  <,,,>
          4   0 3 /
\endpicture}
\qquad
\hbox{\beginpicture
\setcoordinatesystem units <0.8cm,0.7cm>
\put{} at 0 0
\put{} at 6 3.5
\put{$\vdots$} at 2 3.5
\put{$\vdots$} at 4 3.5
\put{$\ssize {01000    \atop 01111}$} at 0 0
\put{$\ssize {1210000  \atop 1333211}$} at 1 1
\put{$\ssize {1110000  \atop 1222111}$} at 2 0
\put{$\ssize {2220000  \atop 2343221}$} at 3 1
\put{$\ssize {111000   \atop 112111}$} at 4 0
\put{$\ssize {121000   \atop 123221}$} at 5 1
\put{$\ssize {01000    \atop 01111}$} at 6 0
\put{$\ssize {2320000  \atop 2454321}$} at 0 2
\put{$\ssize {2320000  \atop 2454321}$} at 2 2
\put{$\ssize {2320000  \atop 2454321}$} at 4 2
\put{$\ssize {2320000  \atop 2454321}$} at 6 2

\arr{0.3 0.3} {0.7 0.7}
\arr{1.3 0.7} {1.7 0.3}
\arr{2.3 0.3} {2.7 0.7}
\arr{3.3 0.7} {3.7 0.3}
\arr{4.3 0.3} {4.7 0.7}
\arr{5.3 0.7} {5.7 0.3}
\arr{0.3 1.7} {0.7 1.3}
\arr{1.3 1.3} {1.7 1.7}
\arr{2.3 1.7} {2.7 1.3}
\arr{3.3 1.3} {3.7 1.7}
\arr{4.3 1.7} {4.7 1.3}
\arr{5.3 1.3} {5.7 1.7}
\arr{0.3 2.3} {0.7 2.7}
\arr{1.3 2.7} {1.7 2.3}
\arr{2.3 2.3} {2.7 2.7}
\arr{3.3 2.7} {3.7 2.3}
\arr{4.3 2.3} {4.7 2.7}
\arr{5.3 2.7} {5.7 2.3}

\setdots<2pt>
\plot 0.7 0  1.3 0 /
\plot 2.7 0  3.3 0 /
\plot 4.7 0  5.3 0 /

\setsolid
\plot 0 0.4  0 1.6 /
\plot 6 0.4  6 1.6 /
\plot 0 2.4  0 3.6 /
\plot 6 2.4  6 3.6 /
\setshadegrid span <1.5mm>
\vshade   0   0 4  <,,,>
          6   0 4 /
\endpicture}
$$
$$
\hbox{\beginpicture
\setcoordinatesystem units <0.8cm,0.7cm>
\put{} at 0 1
\put{} at 10 3.5
\put{$\vdots$} at 2 6.5
\put{$\vdots$} at 4 6.5
\put{$\vdots$} at 6 6.5
\put{$\vdots$} at 8 6.5

\put{$\ssize {0000     \atop 0011    }$} at 0 1
\put{$\ssize {11100    \atop 11111  }$} at 2 1
\put{$\ssize {000000   \atop 011111     }$} at 4 1
\put{$\ssize {1100000  \atop 1111111   }$} at 6 1
\put{$\ssize {011      \atop 011 }$} at 8 1
\put{$\ssize {0000     \atop 0011   }$} at 10 1

\put{$\ssize {11100  \atop   11221}$} at 1 2
\put{$\ssize {111000  \atop  122221 }$} at 3 2
\put{$\ssize {1100000 \atop  1222221  }$} at 5 2
\put{$\ssize {1210000  \atop 1221111}$} at 7 2
\put{$\ssize {0110     \atop 0121 }$} at 9 2

\put{$\ssize {12200     \atop 12321    }$} at 0 3
\put{$\ssize {111000    \atop 123321  }$} at 2 3
\put{$\ssize {2210000   \atop 2333321     }$} at 4 3
\put{$\ssize {1210000  \atop  1332221   }$} at 6 3
\put{$\ssize {1210000      \atop 1232111 }$} at 8 3
\put{$\ssize {12200    \atop 12321   }$} at 10 3

\put{$\ssize {122000 \atop    134321}$} at 1 4
\put{$\ssize {2210000  \atop  2344321 }$} at 3 4
\put{$\ssize {2320000 \atop   2443321 }$} at 5 4
\put{$\ssize {1210000  \atop  1343221}$} at 7 4
\put{$\ssize {2320000   \atop 2343211}$} at 9 4

\put{$\ssize {2320000     \atop 2454321}$} at 0 5
\put{$\ssize {2320000     \atop 2454321}$} at 2 5
\put{$\ssize {2320000     \atop 2454321}$} at 4 5
\put{$\ssize {2320000     \atop 2454321}$} at 6 5
\put{$\ssize {2320000     \atop 2454321}$} at 8 5
\put{$\ssize {2320000     \atop 2454321}$} at 10 5


\arr{0.3 1.3} {0.7 1.7}
\arr{1.3 1.7} {1.7 1.3}
\arr{2.3 1.3} {2.7 1.7}
\arr{3.3 1.7} {3.7 1.3}
\arr{4.3 1.3} {4.7 1.7}
\arr{5.3 1.7} {5.7 1.3}
\arr{6.3 1.3} {6.7 1.7}
\arr{7.3 1.7} {7.7 1.3}
\arr{8.3 1.3} {8.7 1.7}
\arr{9.3 1.7} {9.7 1.3}

\arr{0.3 2.7} {0.7 2.3}
\arr{1.3 2.3} {1.7 2.7}
\arr{2.3 2.7} {2.7 2.3}
\arr{3.3 2.3} {3.7 2.7}
\arr{4.3 2.7} {4.7 2.3}
\arr{5.3 2.3} {5.7 2.7}
\arr{6.3 2.7} {6.7 2.3}
\arr{7.3 2.3} {7.7 2.7}
\arr{8.3 2.7} {8.7 2.3}
\arr{9.3 2.3} {9.7 2.7}

\arr{0.3 3.3} {0.7 3.7}
\arr{1.3 3.7} {1.7 3.3}
\arr{2.3 3.3} {2.7 3.7}
\arr{3.3 3.7} {3.7 3.3}
\arr{4.3 3.3} {4.7 3.7}
\arr{5.3 3.7} {5.7 3.3}
\arr{6.3 3.3} {6.7 3.7}
\arr{7.3 3.7} {7.7 3.3}
\arr{8.3 3.3} {8.7 3.7}
\arr{9.3 3.7} {9.7 3.3}

\arr{0.3 4.7} {0.7 4.3}
\arr{1.3 4.3} {1.7 4.7}
\arr{2.3 4.7} {2.7 4.3}
\arr{3.3 4.3} {3.7 4.7}
\arr{4.3 4.7} {4.7 4.3}
\arr{5.3 4.3} {5.7 4.7}
\arr{6.3 4.7} {6.7 4.3}
\arr{7.3 4.3} {7.7 4.7}
\arr{8.3 4.7} {8.7 4.3}
\arr{9.3 4.3} {9.7 4.7}

\arr{0.3 5.3} {0.7 5.7}
\arr{1.3 5.7} {1.7 5.3}
\arr{2.3 5.3} {2.7 5.7}
\arr{3.3 5.7} {3.7 5.3}
\arr{4.3 5.3} {4.7 5.7}
\arr{5.3 5.7} {5.7 5.3}
\arr{6.3 5.3} {6.7 5.7}
\arr{7.3 5.7} {7.7 5.3}
\arr{8.3 5.3} {8.7 5.7}
\arr{9.3 5.7} {9.7 5.3}

\setdots<2pt>
\plot 0.7 1  1.3 1 /
\plot 2.7 1  3.3 1 /
\plot 4.7 1  5.3 1 /
\plot 6.7 1  7.3 1 /
\plot 8.7 1  9.3 1 /

\setsolid
\plot 0 1.4  0 2.6 /
\plot 10 1.4  10 2.6 /
\plot 0 3.4  0 4.6 /
\plot 10 3.4  10 4.6 /
\plot 0 5.4  0 6.6 /
\plot 10 5.4  10 6.6 /
\setshadegrid span <1.5mm>
\vshade   0   1 7  <,,,>
         10   1 7 /
\endpicture}
$$

It turns out that the preinjective component $\Cal I_0$ contains
the radicals of the modules $P[3]$ and $Y[1]$.  Consider the corresponding
iterated one point extension algebra $B_1$ given by the following quiver $Q_1$
and the indicated relations.

$$
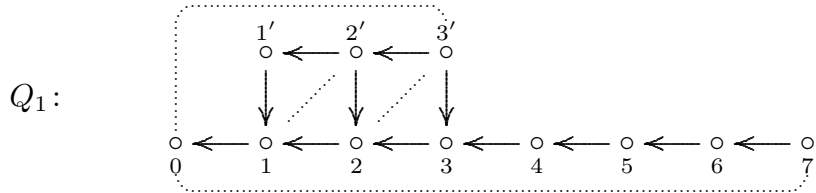
\hbox{\beginpicture
\setcoordinatesystem units <0.5cm,0.5cm>
\put{} at 0 -1.55
\put{} at 16 3.05
\put{$Q_1\:$} at -3 1
\put{$\circ$} at 0 0
\put{$\circ$} at 2 0
\put{$\circ$} at 4 0
\put{$\circ$} at 6 0
\put{$\circ$} at 8 0
\put{$\circ$} at 10 0
\put{$\circ$} at 12 0
\put{$\circ$} at 14 0

\put{$\circ$} at 2 2
\put{$\circ$} at 4 2
\put{$\circ$} at 6 2
\arr{1.6 0}{0.4 0}
\arr{3.6 0}{2.4 0}
\arr{5.6 0}{4.4 0}
\arr{7.6 0}{6.4 0}
\arr{9.6 0}{8.4 0}
\arr{11.6 0}{10.4 0}
\arr{13.6 0}{12.4 0}

\arr{3.6 2}{2.4 2}
\arr{5.6 2}{4.4 2}

\arr{2 1.6}{2 0.4}
\arr{4 1.6}{4 0.4}
\arr{6 1.6}{6 0.4}


\put{$\ssize 1'$} at  2 2.5
\put{$\ssize 2'$} at  4 2.5
\put{$\ssize 3'$} at  6 2.5

\put{$\ssize 0$} at 0 -.5
\put{$\ssize 1$} at 2 -.5
\put{$\ssize 2$} at 4 -.5
\put{$\ssize 3$} at 6 -.5
\put{$\ssize 4$} at 8 -.5
\put{$\ssize 5$} at 10 -.5
\put{$\ssize 6$} at 12 -.5
\put{$\ssize 7$} at 14 -.5

\setdots<2pt>
\plot 3.5 1.5  2.5 0.5 /
\plot 5.5 1.5  4.5 0.5 /
\plot         0 -0.7  0.1 -0.9  0.2 -1   0.3 -1.05 13.7 -1.05 13.8 -1  13.9 -0.9  14 -0.7 /
\plot 0 0.3   0 2.7   0.1  2.9  0.2 3    0.3 3.05   5.7  3.05  5.8  3   5.9  2.9   6  2.7 /
\endpicture}
$$

We picture a part of the component $\Cal I_1$
obtained from $\Cal I_0$ by inserting the two rays for $P[3]$ and $Y[1]$.
From this component we are going to obtain the connecting component $\Cal U$
for $\Cal S_3(7)$. 
$$
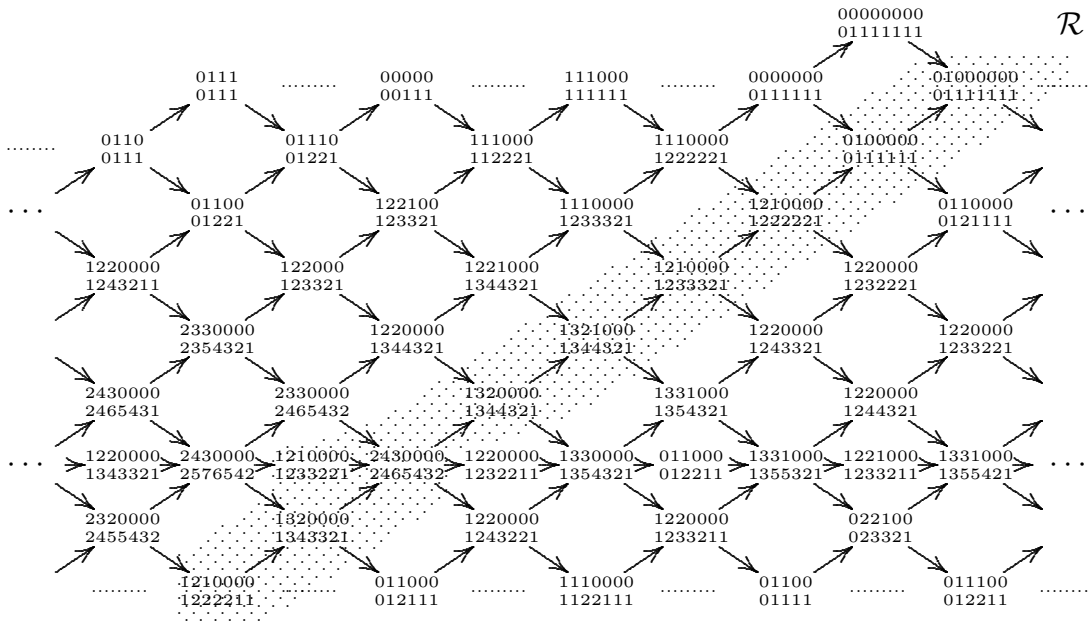
\hbox{\beginpicture
\setcoordinatesystem units <1.05cm,0.7cm>

\put{$\ssize {1210000 \atop 1222211}$} at 2 0
\put{$\ssize {011000  \atop 012111 }$} at 4 0
\put{$\ssize {1110000 \atop 1122111}$} at 6 0
\put{$\ssize {01100   \atop   01111}$} at 8 0
\put{$\ssize {011100  \atop  012211}$} at 10 0

\put{$\ssize {2320000 \atop 2455432}$} at 1 1
\put{$\ssize {1320000 \atop 1343321}$} at 3 1
\put{$\ssize {1220000 \atop 1243221}$} at 5 1
\put{$\ssize {1220000 \atop 1233211}$} at 7 1
\put{$\ssize {022100  \atop  023321}$} at 9 1

\put{$\cdots$} at 0 2
\put{$\ssize {1220000 \atop 1343321}$} at 1 2
\put{$\ssize {2430000 \atop 2576542}$} at 2 2
\put{$\ssize {1210000 \atop 1233221}$} at 3 2
\put{$\ssize {2430000 \atop 2465432}$} at 4 2
\put{$\ssize {1220000 \atop 1232211}$} at 5 2
\put{$\ssize {1330000 \atop 1354321}$} at 6 2
\put{$\ssize {011000  \atop  012211}$} at 7 2
\put{$\ssize {1331000 \atop 1355321}$} at 8 2
\put{$\ssize {1221000 \atop 1233211}$} at 9 2
\put{$\ssize {1331000 \atop 1355421}$} at 10 2
\put{$\cdots$} at 11 2

\put{$\ssize {2430000 \atop 2465431}$} at 1 3
\put{$\ssize {2330000 \atop 2465432}$} at 3 3
\put{$\ssize {1320000 \atop 1344321}$} at 5 3
\put{$\ssize {1331000 \atop 1354321}$} at 7 3
\put{$\ssize {1220000 \atop 1244321}$} at 9 3

\put{$\ssize {2330000 \atop 2354321}$} at 2 4
\put{$\ssize {1220000 \atop 1344321}$} at 4 4
\put{$\ssize {1321000 \atop 1344321}$} at 6 4
\put{$\ssize {1220000 \atop 1243321}$} at 8 4
\put{$\ssize {1220000 \atop 1233221}$} at 10 4

\put{$\ssize {1220000 \atop 1243211}$} at 1 5
\put{$\ssize {122000  \atop  123321}$} at 3 5
\put{$\ssize {1221000 \atop 1344321}$} at 5 5
\put{$\ssize {1210000 \atop 1233321}$} at 7 5
\put{$\ssize {1220000 \atop 1232221}$} at 9 5

\put{$\cdots$} at 0 6
\put{$\ssize {01100   \atop   01221}$} at 2 6
\put{$\ssize {122100  \atop  123321}$} at 4 6
\put{$\ssize {1110000 \atop 1233321}$} at 6 6
\put{$\ssize {1210000 \atop 1222221}$} at 8 6
\put{$\ssize {0110000 \atop 0121111}$} at 10 6
\put{$\cdots$} at 11 6

\put{$\ssize {0110    \atop    0111}$} at 1 7
\put{$\ssize {01110   \atop   01221}$} at 3 7
\put{$\ssize {111000  \atop  112221}$} at 5 7
\put{$\ssize {1110000 \atop 1222221}$} at 7 7
\put{$\ssize {0100000 \atop 0111111}$} at 9 7

\put{$\ssize {0111    \atop    0111}$} at 2 8
\put{$\ssize {00000   \atop   00111}$} at 4 8
\put{$\ssize {111000  \atop  111111}$} at 6 8
\put{$\ssize {0000000 \atop 0111111}$} at 8 8
\put{$\ssize {01000000\atop01111111}$} at 10 8

\put{$\ssize {00000000\atop01111111}$} at 9 9
\arr{0.3 0.3} {0.7 0.7}
\arr{2.3 0.3} {2.7 0.7}
\arr{4.3 0.3} {4.7 0.7}
\arr{6.3 0.3} {6.7 0.7}
\arr{8.3 0.3} {8.7 0.7}
\arr{10.3 0.3} {10.7 0.7}

\arr{1.3 0.7} {1.7 0.3}
\arr{3.3 0.7} {3.7 0.3}
\arr{5.3 0.7} {5.7 0.3}
\arr{7.3 0.7} {7.7 0.3}
\arr{9.3 0.7} {9.7 0.3}

\arr{0.3 1.7} {0.7 1.3} 
\arr{2.3 1.7} {2.7 1.3} 
\arr{4.3 1.7} {4.7 1.3} 
\arr{6.3 1.7} {6.7 1.3} 
\arr{8.3 1.7} {8.7 1.3} 
\arr{10.3 1.7} {10.7 1.3} 

\arr{1.3 1.3} {1.7 1.7} 
\arr{3.3 1.3} {3.7 1.7} 
\arr{5.3 1.3} {5.7 1.7} 
\arr{7.3 1.3} {7.7 1.7} 
\arr{9.3 1.3} {9.7 1.7} 

\arr{0.4 2.0} {0.6 2.0}
\arr{1.4 2.0} {1.6 2.0}
\arr{2.4 2.0} {2.6 2.0}
\arr{3.4 2.0} {3.6 2.0}
\arr{4.4 2.0} {4.6 2.0}
\arr{5.4 2.0} {5.6 2.0}
\arr{6.4 2.0} {6.6 2.0}
\arr{7.4 2.0} {7.6 2.0}
\arr{8.4 2.0} {8.6 2.0}
\arr{9.4 2.0} {9.6 2.0}
\arr{10.4 2.0} {10.6 2.0}

\arr{0.3 2.3} {0.7 2.7} 
\arr{2.3 2.3} {2.7 2.7} 
\arr{4.3 2.3} {4.7 2.7} 
\arr{6.3 2.3} {6.7 2.7} 
\arr{8.3 2.3} {8.7 2.7} 
\arr{10.3 2.3} {10.7 2.7} 

\arr{1.3 2.7} {1.7 2.3} 
\arr{3.3 2.7} {3.7 2.3} 
\arr{5.3 2.7} {5.7 2.3} 
\arr{7.3 2.7} {7.7 2.3} 
\arr{9.3 2.7} {9.7 2.3} 

\arr{0.3 3.7} {0.7 3.3} 
\arr{2.3 3.7} {2.7 3.3} 
\arr{4.3 3.7} {4.7 3.3} 
\arr{6.3 3.7} {6.7 3.3} 
\arr{8.3 3.7} {8.7 3.3} 
\arr{10.3 3.7} {10.7 3.3} 

\arr{1.3 3.3} {1.7 3.7} 
\arr{3.3 3.3} {3.7 3.7} 
\arr{5.3 3.3} {5.7 3.7} 
\arr{7.3 3.3} {7.7 3.7} 
\arr{9.3 3.3} {9.7 3.7} 

\arr{0.3 4.3} {0.7 4.7} 
\arr{2.3 4.3} {2.7 4.7} 
\arr{4.3 4.3} {4.7 4.7} 
\arr{6.3 4.3} {6.7 4.7} 
\arr{8.3 4.3} {8.7 4.7} 
\arr{10.3 4.3} {10.7 4.7} 

\arr{1.3 4.7} {1.7 4.3} 
\arr{3.3 4.7} {3.7 4.3} 
\arr{5.3 4.7} {5.7 4.3} 
\arr{7.3 4.7} {7.7 4.3} 
\arr{9.3 4.7} {9.7 4.3} 

\arr{0.3 5.7} {0.7 5.3} 
\arr{2.3 5.7} {2.7 5.3} 
\arr{4.3 5.7} {4.7 5.3} 
\arr{6.3 5.7} {6.7 5.3} 
\arr{8.3 5.7} {8.7 5.3} 
\arr{10.3 5.7} {10.7 5.3} 

\arr{1.3 5.3} {1.7 5.7} 
\arr{3.3 5.3} {3.7 5.7} 
\arr{5.3 5.3} {5.7 5.7} 
\arr{7.3 5.3} {7.7 5.7} 
\arr{9.3 5.3} {9.7 5.7} 

\arr{0.3 6.3} {0.7 6.7} 
\arr{2.3 6.3} {2.7 6.7} 
\arr{4.3 6.3} {4.7 6.7} 
\arr{6.3 6.3} {6.7 6.7} 
\arr{8.3 6.3} {8.7 6.7} 
\arr{10.3 6.3} {10.7 6.7} 

\arr{1.3 6.7} {1.7 6.3} 
\arr{3.3 6.7} {3.7 6.3} 
\arr{5.3 6.7} {5.7 6.3} 
\arr{7.3 6.7} {7.7 6.3} 
\arr{9.3 6.7} {9.7 6.3} 

\arr{2.3 7.7} {2.7 7.3} 
\arr{4.3 7.7} {4.7 7.3} 
\arr{6.3 7.7} {6.7 7.3} 
\arr{8.3 7.7} {8.7 7.3} 
\arr{10.3 7.7} {10.7 7.3} 

\arr{1.3 7.3} {1.7 7.7} 
\arr{3.3 7.3} {3.7 7.7} 
\arr{5.3 7.3} {5.7 7.7} 
\arr{7.3 7.3} {7.7 7.7} 
\arr{9.3 7.3} {9.7 7.7} 

\arr{8.3 8.3} {8.7 8.7}
\arr{9.3 8.7} {9.7 8.3}
\setdots<2pt>
\plot 0.7 0  1.3 0 /
\plot 2.7 0  3.3 0 /
\plot 4.7 0  5.3 0 /
\plot 6.7 0  7.3 0 /
\plot 8.7 0  9.3 0 /
\plot 10.7 0  11.3 0 /

\plot -0.2 7  0.3 7 /

\plot 2.7 8  3.3 8 /
\plot 4.7 8  5.3 8 /
\plot 6.7 8  7.3 8 /
\plot 10.7 8 11.3 8 /

\put{$\Cal R$} at 11 9
\setshadegrid span <0.8mm>
\hshade -0.5 1.5 2.5 <,,,z> 0.5 1.5 3.5 <,,z,z> 1.5 2.5 4.5  <,,z,z> 2.2 2.5 4.5 <,,z,z> 2.7 3.7 5.7 <,,z,z> 8 9 11 <,,z,> 8.5 9.5 11 / 
\endpicture}
$$
The diagram actually pictures part of the connecting component $\Cal U$
as we have already deleted the
module $S_2'$ which occurs as $\tau^{-1}Y[1]$ (in $\Cal I_1$)
in the upper right corner. 
Due to the presence of $S_2'$, not all modules in $\Cal I_1$ are 
in $\Cal S_3(7)$, but certainly those in the shaded region $\Cal R$
(which forms a slice) and those on the left of $\Cal R$ are.
Each module $M$ in $\Cal S_3(7)$ 
which lies in $\Cal R$ or which admits a path from $\Cal R$ has
no nonzero morphism into the module $Y[1]$ and hence is a module over
the algebra $B_2$ given by the quiver $Q_2$ and the relations as indicated.
Namely, if $\Hom(M,Y[1])=0$, then also 
$\Hom((0\sub (M_i/\Im M(\iota_i))),Y[1])=0$ and hence $M(\iota_1)$ is an
isomorphism and $M$ a $B_2$-module, as seen in Lemma~8.
$$
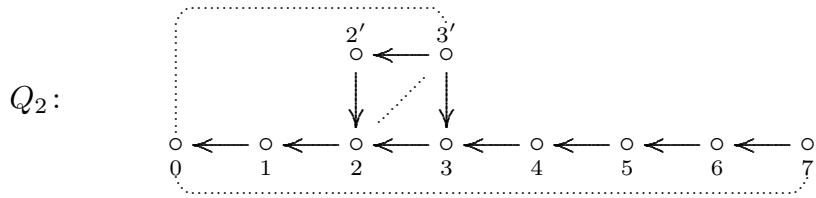
\hbox{\beginpicture
\setcoordinatesystem units <0.5cm,0.5cm>
\put{} at 0 -1.55
\put{} at 16 3.05
\put{$Q_2\:$} at -3 1
\put{$\circ$} at 0 0
\put{$\circ$} at 2 0
\put{$\circ$} at 4 0
\put{$\circ$} at 6 0
\put{$\circ$} at 8 0
\put{$\circ$} at 10 0
\put{$\circ$} at 12 0
\put{$\circ$} at 14 0

\put{$\circ$} at 4 2
\put{$\circ$} at 6 2
\arr{1.6 0}{0.4 0}
\arr{3.6 0}{2.4 0}
\arr{5.6 0}{4.4 0}
\arr{7.6 0}{6.4 0}
\arr{9.6 0}{8.4 0}
\arr{11.6 0}{10.4 0}
\arr{13.6 0}{12.4 0}

\arr{5.6 2}{4.4 2}

\arr{4 1.6}{4 0.4}
\arr{6 1.6}{6 0.4}


\put{$\ssize 2'$} at  4 2.5
\put{$\ssize 3'$} at  6 2.5

\put{$\ssize 0$} at 0 -.5
\put{$\ssize 1$} at 2 -.5
\put{$\ssize 2$} at 4 -.5
\put{$\ssize 3$} at 6 -.5
\put{$\ssize 4$} at 8 -.5
\put{$\ssize 5$} at 10 -.5
\put{$\ssize 6$} at 12 -.5
\put{$\ssize 7$} at 14 -.5

\setdots<2pt>
\plot 5.5 1.5  4.5 0.5 /
\plot         0 -0.7  0.1 -0.9  0.2 -1   0.3 -1.05 13.7 -1.05 13.8 -1  13.9 -0.9  14 -0.7 /
\plot 0 0.3   0 2.7   0.1  2.9  0.2 3    0.3 3.05   5.7  3.05  5.8  3   5.9  2.9   6  2.7 /
\endpicture}
$$

Note that by deleting the point 0 in $Q_2$ we get back the quiver $Q_0$,
up to a shift. Let $\overline I_0$ be the socle factor of the $B_2$-module
$I_0$.  Then $\overline I_0[-1]$ occurs in the preprojective component
$\Cal P$ for $B_0$. The component $\Cal P[1]$ with the corresponding 
coextension becomes the preprojective component $\Cal P_2$ for $B_2$.
Thus, the Auslander-Reiten quiver for $B_2$ consists of 
this preprojective component
$\Cal P_2$, the tubular family $\Cal T=\Cal T_0[1]$, 
and the preinjective component
$\Cal I_0[1]$.
Since the modules $S_2'$ and $S_3'$ occur in $\Cal I_0[1]$,
all modules in $\Cal P_2$ and in $\Cal T$ are objects in $\Cal S_3(7)$. 
The shaded region $\Cal R$ above occurs as a slice
in $\Cal P_2$. Let $\Cal U$ consist of the left hand part of $\Cal I_0$
and the right hand part of $\Cal P_2$, with the two copies of $\Cal R$
identified. 
We observe that the injective module $I[1]$
occurs as $\tau^{-5}Z$ where $Z$ is the irreducible successor of $Y[1]$.
Hence $\Cal U$ has stable type 
$\Bbb Z\widetilde{\Bbb E}_8$ and there are two nonstable $\tau$-orbits 
of length ten and one attached, as sketched above. 

\smallskip
We claim that $\Cal D=\Cal U\vee\Cal T$ 
forms a fundamental domain for $\Cal S_3(7)$, modulo the shift $M\mapsto M[1]$.
Clearly, $\Cal D$ consists of pairwise nonisomorphic modules 
since they are pairwise nonisomorphic as $B_1$-modules.
As $\Cal D'=\Cal P\vee\Cal T_0\vee\Cal U$ contains the injective $I[0]$
and the relatively injective $Y[1]$ but 
no unextended radicals of projective $\Cal S_3(7)$
representations, it follows that every 
indecomposable in $\Cal S_3(7)$ has a translate
in $\Cal D'$. 
Moreover, each module in $\Cal P$ occurs in some
category $\Cal T[i]$ or $\Cal U[i]$, as indicated below.  For example, the three wings
in $\Cal P\cap\Cal T[-1]$ occur in the three big tubes in $\Cal T[-1]$. 
$$
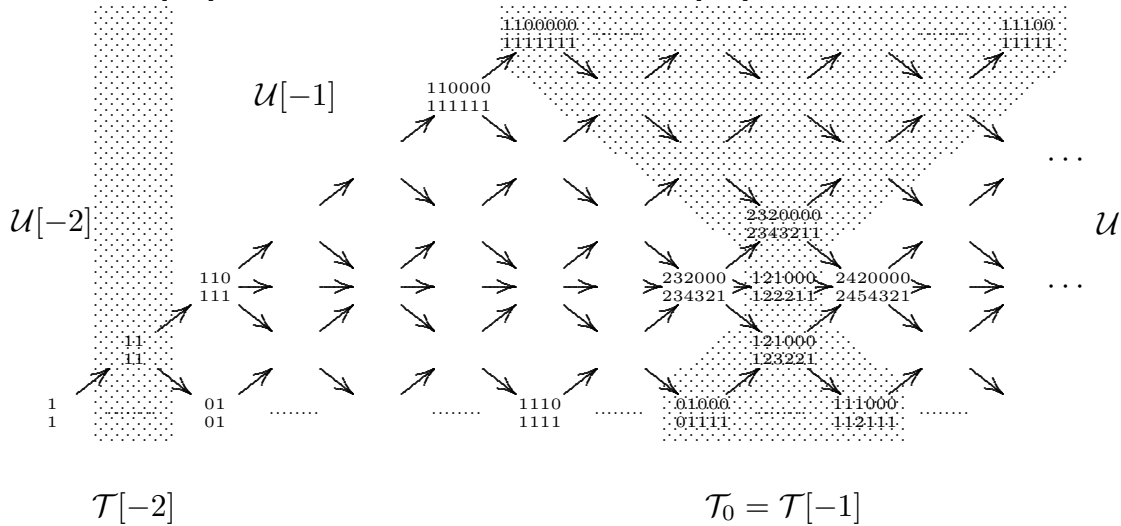
\hbox{\beginpicture
\setcoordinatesystem units <0.9cm,0.7cm>
\put{} at 0 0
\put{} at 12 6

\put{$\ssize {1       \atop       1}$} at 0 0
\put{$\ssize {01 \atop           01}$} at 2 0
\put{$\ssize {1110    \atop    1111}$} at 6 0
\put{$\ssize {01000   \atop   01111}$} at 8 0
\put{$\ssize {111000  \atop  112111}$} at 10 0

\put{$\ssize {11      \atop      11}$} at 1 1
\put{$\ssize {121000  \atop  123221}$} at 9 1

\put{$\ssize {110     \atop     111}$} at 2 2
\put{$\ssize {232000  \atop  234321}$} at 7.9 2
\put{$\ssize {121000  \atop  122211}$} at 9 2
\put{$\ssize {2420000 \atop 2454321}$} at 10.1 2

\put{$\ssize {2320000 \atop 2343211}$} at 9 3


\put{$\ssize {110000  \atop  111111}$} at 5 5

\put{$\ssize {1100000 \atop 1111111}$} at 6 6
\put{$\ssize {11100   \atop   11111}$} at 12 6
\arr{0.3 0.3} {0.7 0.7}
\arr{2.3 0.3} {2.7 0.7}
\arr{4.3 0.3} {4.7 0.7}
\arr{6.3 0.3} {6.7 0.7}
\arr{8.3 0.3} {8.7 0.7}
\arr{10.3 0.3} {10.7 0.7}

\arr{1.3 0.7} {1.7 0.3}
\arr{3.3 0.7} {3.7 0.3}
\arr{5.3 0.7} {5.7 0.3}
\arr{7.3 0.7} {7.7 0.3}
\arr{9.3 0.7} {9.7 0.3}
\arr{11.3 0.7} {11.7 0.3}

\arr{2.3 1.7} {2.7 1.3} 
\arr{4.3 1.7} {4.7 1.3} 
\arr{6.3 1.7} {6.7 1.3} 
\arr{8.3 1.7} {8.7 1.3} 
\arr{10.3 1.7} {10.7 1.3} 

\arr{1.3 1.3} {1.7 1.7} 
\arr{3.3 1.3} {3.7 1.7} 
\arr{5.3 1.3} {5.7 1.7} 
\arr{7.3 1.3} {7.7 1.7} 
\arr{9.3 1.3} {9.7 1.7} 
\arr{11.3 1.3} {11.7 1.7} 

\arr{2.3 2.0} {2.7 2.0}
\arr{3.3 2.0} {3.7 2.0}
\arr{4.3 2.0} {4.7 2.0}
\arr{5.3 2.0} {5.7 2.0}
\arr{6.3 2.0} {6.7 2.0}
\arr{7.3 2.0} {7.5 2.0}
\arr{8.4 2.0} {8.6 2.0}
\arr{9.4 2.0} {9.6 2.0}
\arr{10.6 2.0} {10.8 2.0}
\arr{11.3 2.0} {11.7 2.0}

\arr{2.3 2.3} {2.7 2.7} 
\arr{4.3 2.3} {4.7 2.7} 
\arr{6.3 2.3} {6.7 2.7} 
\arr{8.3 2.3} {8.7 2.7} 
\arr{10.3 2.3} {10.7 2.7} 

\arr{3.3 2.7} {3.7 2.3} 
\arr{5.3 2.7} {5.7 2.3} 
\arr{7.3 2.7} {7.7 2.3} 
\arr{9.3 2.7} {9.7 2.3} 
\arr{11.3 2.7} {11.7 2.3} 

\arr{4.3 3.7} {4.7 3.3} 
\arr{6.3 3.7} {6.7 3.3} 
\arr{8.3 3.7} {8.7 3.3} 
\arr{10.3 3.7} {10.7 3.3} 

\arr{3.3 3.3} {3.7 3.7} 
\arr{5.3 3.3} {5.7 3.7} 
\arr{7.3 3.3} {7.7 3.7} 
\arr{9.3 3.3} {9.7 3.7} 
\arr{11.3 3.3} {11.7 3.7} 

\arr{4.3 4.3} {4.7 4.7} 
\arr{6.3 4.3} {6.7 4.7} 
\arr{8.3 4.3} {8.7 4.7} 
\arr{10.3 4.3} {10.7 4.7} 

\arr{5.3 4.7} {5.7 4.3} 
\arr{7.3 4.7} {7.7 4.3} 
\arr{9.3 4.7} {9.7 4.3} 
\arr{11.3 4.7} {11.7 4.3} 

\arr{6.3 5.7} {6.7 5.3} 
\arr{8.3 5.7} {8.7 5.3} 
\arr{10.3 5.7} {10.7 5.3} 

\arr{5.3 5.3} {5.7 5.7} 
\arr{7.3 5.3} {7.7 5.7} 
\arr{9.3 5.3} {9.7 5.7} 
\arr{11.3 5.3} {11.7 5.7} 
\setdots<2pt>
\plot 0.7 0  1.3 0 /
\plot 2.7 0  3.3 0 /
\plot 4.7 0  5.3 0 /
\plot 6.7 0  7.3 0 /
\plot 8.7 0  9.3 0 /
\plot 10.7 0  11.3 0 /

\plot 6.7 6  7.3 6 /
\plot 8.7 6  9.3 6 /
\plot 10.7 6 11.3 6 /

\put{$\cdots$} at 12.5 4
\put{$\cdots$} at 12.5 2
\setshadegrid span <0.6mm>
\hshade -0.5  7.5 10.5 <,,,z>  0.5 7.5 10.5 <,,z,z> 1.5 8.5 9.5 <,,z,> 2.5 8.5 9.5 <,,z,> 5.5 5.5 12.5 <,,z,> 6.5 5.5 12.5 /
\hshade -0.5  0.5 1.5  6.5 0.5 1.5 /

\put{$\Cal U[-2]$} at 0 3
\put{$\Cal U[-1]$} at 3 5
\put{$\Cal U$} at 13 3
\put{$\Cal T[-2]$} at 1 -1.5
\put{$\Cal T_0=\Cal T[-1]$} at 9 -1.5
\endpicture}
$$

\smallskip
We show that $\Cal D$ consists of components of the Auslander-Reiten
quiver for $\Cal S_3(7)$. The source maps in $\Cal D$ are as follows:
If $M\in\Cal D$ has a path into $\Cal R$ then the source map $M\to N$
is a source map in the category $B_1$-mod, otherwise it is a source map
in $B_2$-mod.  Let $B_*=B_1$ or $B_*=B_2$ for the first and second case,
respectively. In each case, $M\to N$ is a source map in the category
$\Cal S_3(7)$: Let $T\in \Cal D[i]$ be an indecomposable object and
$t\:M\to T$ a nonisomorphism.  If $i=0$ and $T\in B_*$-mod, then $t$
factors through $N$ (this includes the case that $B_*=B_1$ and $T$ is 
a $B_2$-module). Otherwise, $T$ has a path into $\Cal R$ and $M$ has
not, so $\Hom_{B_1}(M,T)=0$. 
If $i>0$ then $T$ is a module over an iterated one point extension algebra
$B^+$ for $B_*$ and the component for $M$ in the Auslander-Reiten quiver
for $B_*$-mod is also a component for the Auslander-Reiten quiver
for $B^+$ (since all radicals of projective indecomposable objects have
been extended). Again, $t$ factors through $N$. 
Finally, if $i<0$ then $\Hom_{\Cal S}(M,T)=\Hom_{B^+}(M[-i],T[-i])=0$
since $T[-i]$ is a $B_1$-module and $M$ a module over an iterated 
one point extension algebra $B^+$. 
Hence the Auslander-Reiten quiver 
for $\Cal S_3(7)$ is as pictured below.

$$
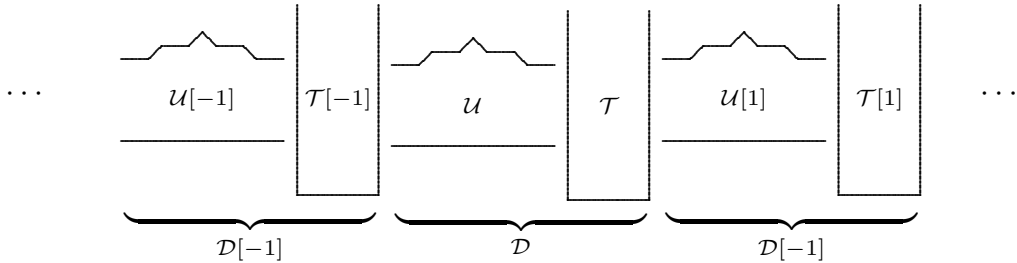
\hbox{\beginpicture
\setcoordinatesystem units <0.6cm,0.6cm>
\put{\beginpicture
\setcoordinatesystem units <0.3cm,0.3cm>
\put{} at 0 0
\put{} at 10 7 
\put{$\underbrace{\phantom{mmmmmmmmm}}_{\Cal D[-1]}$} at 5.25 -1.3

\plot 0.5 2  6.5 2 /
\plot 0.5 5  1.5 5  2 5.5  3 5.5  3.5 6  4 5.5  5 5.5  5.5 5  6.5 5 /
\plot 7 7  7 0  10 0  10 7 /

\put{$\ssize \Cal T[-1]$} at 8.5 3.5
\put{$\ssize \Cal U[-1]$} at 3.5 3.5
\endpicture} at 0 0 
\put{\beginpicture
\setcoordinatesystem units <0.3cm,0.3cm>
\put{} at 0 0
\put{} at 10 7 
\put{$\underbrace{\phantom{mmmmmmmmm}}_{\Cal D}$} at 5.25 -1

\plot 0.5 2  6.5 2 /
\plot 0.5 5  1.5 5  2 5.5  3 5.5  3.5 6  4 5.5  5 5.5  5.5 5  6.5 5 /
\plot 7 7  7 0  10 0  10 7 /

\put{$\ssize \Cal T$} at 8.5 3.5
\put{$\ssize \Cal U$} at 3.5 3.5
\endpicture} at 5 0 
\put{\beginpicture
\setcoordinatesystem units <0.3cm,0.3cm>
\put{} at 0 0
\put{} at 10 7 
\put{$\underbrace{\phantom{mmmmmmmmm}}_{\Cal D[-1]}$} at 5.25 -1.3

\plot 0.5 2  6.5 2 /
\plot 0.5 5  1.5 5  2 5.5  3 5.5  3.5 6  4 5.5  5 5.5  5.5 5  6.5 5 /
\plot 7 7  7 0  10 0  10 7 /

\put{$\ssize \Cal T[1]$} at 8.5 3.5
\put{$\ssize \Cal U[1]$} at 3.5 3.5
\endpicture} at 10 0 
\put{$\cdots$} at -4 0.75
\put{$\cdots$} at 14 0.75
\endpicture}
$$

\noindent{\bf Proposition 13.} (Separation properties of $\Cal S_3(7)$)

\smallskip
\item{1.} {\it Let 
        $\Cal P_T[i]=\bigvee_{j<i}\big(\Cal T[j]\vee\Cal U[j+1]\big)$ and
        $\Cal I_T[i]=\bigvee_{j>i}\big(\Cal U[j]\vee\Cal T[j]\big)$. 
        Then
        $\Cal T[i]$ separates $\Cal P_T[i]$ from $\Cal I_T[i]$ 
        in the sense that}
        $$\Hom(\Cal T[i],\Cal P_T[i])=
        \Hom(\Cal I_T[i],\Cal T[i])=
        \Hom(\Cal I_T[i],\Cal P_T[i])=0$$
        {\it and that for every map $f\in\Hom(M,N)$ where
        $M\in \Cal P_T[i]$ and $N\in \Cal I_T[i]$, and for every
        tube $\Cal T_\gamma$ in $\Cal T[i]$ the map 
        $f$ factors through a sum of objects in 
        $\Cal T_\gamma$. }
\item{2.} {\it If $\Cal P_U[i]=\bigvee_{j<i}\big(\Cal U[j]\vee\Cal T[j]\big)$ 
        and $\Cal I_U[i]=\bigvee_{j>i}\big(\Cal T[j-1]\vee\Cal U[j]\big)$ then
        $\Cal U[i]$ separates $\Cal P_U[i]$ from $\Cal I_U[i]$ 
        in the sense that }
        $$\Hom(\Cal U[i],\Cal P_U[i])=
        \Hom(\Cal I_U[i],\Cal U[i])=
        \Hom(\Cal I_U[i],\Cal P_U[i])=0$$
        {\it holds and that for every map $f\in\Hom(M,N)$ where
         $M\in \Cal P_U[i]$ and $N\in \Cal I_U[i]$ and for 
        every slice $\Cal R$ in $\Cal U[i]$ the map
        $f$ factors through a sum of objects in 
        $\Cal R$. }

\smallskip\noindent{\it Proof:}
First we determine the minimal left and right approximations for certain
objects in $\Cal S_3(7)$ in $B_2$-mod.
Suppose that $M\in\Cal S_3(7)$ has support in the set
$\{i|i\leq 7\}\cup\{i'|i\leq 3\}$. Then the map 
$$\xymatrix@1@=5mm{\ssize M_{-1}' \ar[d]_{j_{-1}}
        & \ssize M_0'\ar[d]^{j_0} \ar[l]_{a_0} 
        & \ssize M_1' \ar[d]^{j_1} \ar[l]_{a_1}
                & \ssize M_2' \ar[l]_{a_2} \ar[d]^{j_2} & \cdots \ar[l] \\
       \ssize M_{-1} & \ssize M_0 \ar[l]^{b_0} & \ssize M_1 \ar[l]^{b_1} 
                & \ssize M_2 \ar[l]^{b_2} & \cdots \ar[l]}
\quad\beginpicture\setcoordinatesystem units <1cm,1cm> 
        \arr{-.5 -.4}{.5 -.4} 
        \put{} at -.5 -1
        \put{} at .5 0 \endpicture
\quad \xymatrix@1@=5mm{\ssize 0 
                & \ssize M_0\ar[d]^1\ar[l] 
                & \ssize M_1 \ar[d]^1 \ar[l]_{b_1} 
                & \ssize M_2' \ar[l]_{j_2b_2} \ar[d]^{j_2} & \cdots \ar[l] \\
        \ssize 0 & \ssize M_0 \ar[l] & \ssize M_1 \ar[l]^{b_1} 
                & \ssize M_2 \ar[l]^{b_2} & \cdots \ar[l]}$$
is a minimal left approximation for $M$ in $B_2$-mod.
Similarly, if $M$ has support in $\{i|i\geq0\}\cup\{i'|i\geq0\}$ and if
$j_0$ and $j_1$ are isomorphisms, then the 
inclusion
$$\xymatrix@1@=5mm{\cdots & \ssize M_3' \ar[l] \ar[d]
        & \ssize 0 \ar[d] \ar[l]
        & \cdots  \ar[l]
                & \ssize 0 \ar[l] \ar[d] & \ssize 0 \ar[l] \ar[d] \\
       \cdots & \ssize M_3 \ar[l] & \ssize M_4 \ar[l] & \cdots \ar[l]
                & \ssize M_7 \ar[l] & \ssize 0 \ar[l]}
\quad\beginpicture\setcoordinatesystem units <1cm,1cm> 
        \arr{-.5 -.4}{.5 -.4} 
        \put{} at -.5 -1
        \put{} at .5 0 \endpicture
\quad \xymatrix@1@=5mm{\cdots & \ssize M_3' \ar[d] \ar[l]
        & \ssize M_4' \ar[d] \ar[l]
        & \cdots  \ar[l]
                & \ssize M_7' \ar[l] \ar[d] & \ssize M_8' \ar[l]\ar[d] \\
       \cdots & \ssize M_3 \ar[l] & \ssize M_4 \ar[l] & \cdots \ar[l]
                & \ssize M_7 \ar[l] & \ssize M_8 \ar[l]} $$
is a minimal right approximation for $M$ in $B_2$-mod.

\smallskip For the proof of the first assertion, we assume that $i=0$.
Then $\Cal T[0]$ is the tubular family for $B_2$. 
The statement about the nonexistence 
of homomorphisms follows from the corresponding statement about 
preprojective, regular, and preinjective $B_2$-modules and the fact
that the minimal left approximation of an object in $\Cal P_T[0]$ 
is a preprojective $B_2$-module and the minimal right approximation of
an object in $\Cal I_T[0]$ is a preinjective $B_2$-module. 
If $f\in \Hom(M,N)$ and $\Cal T_\gamma$ a tube in $\Cal T[0]$,
then $f$ factors through both the minimal left approximation for $M$
and the minimal right approximation for $N$ and hence through $\Cal T_\gamma$.

\smallskip For the second assertion we use the fact that the minimal left approximation for an 
object in $\Cal P_U[1]$ is a preprojective or regular $B_2$-module,
that the minimal right approximation for an object in $\Cal I_U[0]$
is regular or preinjective, that the minimal right approximation for an 
object in $\Cal U[1]$ is preinjective as a $B_2$-module and that
the minimal left approximation for a module in $\Cal U[0]$ is preprojective.
Each claim  about the nonexistence of homomorphisms follows for either
$i=0$ or $i=1$.  To show the last assertion, let $f\in\Hom(M,N)$
where $M\in \Cal P_U[0]$ and $N\in\Cal I_U[0]$, 
and let $\Cal R$ be a slice in $\Cal U[0]$. 
Since the minimal left approximation
for $M$ is a preprojective $B_2$-module, $f$ factors through a preprojective
$B_2$-module. Using the property of source maps in $B_2$-mod, we 
obtain that $f$ factors through some slice in $\Cal U[0]$ which we may assume
to be on the right hand side of $\Cal R$ (where all objects are
preprojective $B_2$-modules); and using the property of sink
maps in $\Cal U$ it follows that $f$ actually factors through $\Cal R$. \qed

\smallskip
Note that all modules in $\Cal U$ and in $\Cal T$ are  modules over some 
finite dimensional algebra, and hence $\Cal S_3(7)$ is locally support
finite.  Theorem~5 above yields the following information
about the corresponding category $\Cal S_3(k[T]/T^7)$.

\medskip\noindent{\bf Theorem 14.} {\it
The category $\Cal S_3(k[T]/T^7)$ consists of a $\Bbb P_1(k)$-family
of stable tubes of type $(5,3,2)$ and of a connecting component 
$\Cal U$ of stable orbit type $\Bbb Z\widetilde{\Bbb E}_8$ with two nonstable 
orbits of length 10 and 1 attached.  Homomorphisms in the infinite
radical of $\Cal T$ factor through any slice in $\Cal U$, and 
maps in the infinite radical of $\Cal U$ factor through any of the
tubes in $\Cal T$. } \qed

\def\boundarythreeseven{\beginpicture\setcoordinatesystem units <6.3mm,6.3mm>
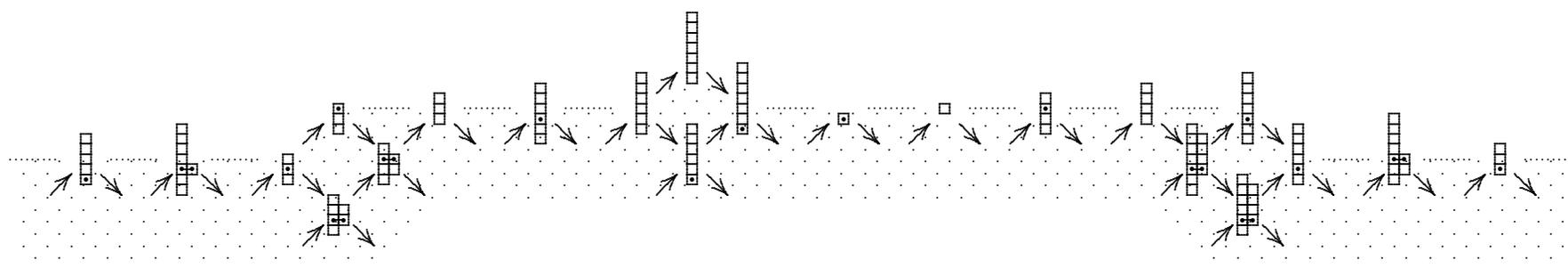
\put{The Nonstable Modules in the Connecting Component in  
$\Cal S_3(k[T]/T^7)$} at 15 -2
\multiput{$\smallsq2$} at .9 .6 *4 0 .2 /
\multiput{$\boldkey.$} at 1 .6 /

\multiput{$\smallsq2$} at 2.8 .4 *6 0 .2 /
\multiput{$\smallsq2$} at 3 .8  /
\multiput{$\boldkey.$} at 2.9 .8  3.1 .8 /
\plot 2.9 .8  3.1 .8 /

\multiput{$\smallsq2$} at 4.9 .6 *2 0 .2 /
\multiput{$\boldkey.$} at 5 .8 /

\multiput{$\smallsq2$} at 5.9 1.6 *2 0 .2 /
\multiput{$\boldkey.$} at 6 2 /

\multiput{$\smallsq2$} at 5.8 -.4 *3 0 .2 /
\multiput{$\smallsq2$} at 6 -.2 *1 0 .2 /
\multiput{$\boldkey.$} at 5.9 -.2  6.1 -.2 /
\plot  5.9 -.2  6.1 -.2 /

\multiput{$\smallsq2$} at 6.8 .6 *3  0 .2 /
\multiput{$\smallsq2$} at 7 .8 *1  0 .2 /
\multiput{$\boldkey.$} at 6.9 1  7.1 1 /
\plot  6.9 1  7.1 1 /

\multiput{$\smallsq2$} at 7.9 1.8 *2 0 .2 /

\multiput{$\smallsq2$} at 9.9 1.4 *5 0 .2 /
\multiput{$\boldkey.$} at 10 1.8 /

\multiput{$\smallsq2$} at 11.9 1.6 *5 0 .2 /

\multiput{$\smallsq2$} at 12.9 2.6 *6  0 .2 /

\multiput{$\smallsq2$} at 12.9 .6 *5 0 .2 /
\multiput{$\boldkey.$} at 13 .6  /

\multiput{$\smallsq2$} at 13.9 1.6 *6 0 .2 /
\multiput{$\boldkey.$} at 14 1.6 /

\multiput{$\smallsq2$} at 15.9 1.8 /
\multiput{$\boldkey.$} at 16 1.8 /

\multiput{$\smallsq2$} at 17.9 2 /

\multiput{$\smallsq2$} at 19.9 1.6 *3 0 .2 /
\multiput{$\boldkey.$} at 20 2 /

\multiput{$\smallsq2$} at 21.9 1.8 *3  0 .2 /

\multiput{$\smallsq2$} at 22.8 .4 *6 0 .2 /
\multiput{$\smallsq2$} at 23 .8 *3 0 .2 /
\multiput{$\boldkey.$} at 22.9 .8  23.1 .8  /
\plot 22.9 .8  23.1 .8  /

\multiput{$\smallsq2$} at 23.9 1.4 *6 0 .2 /
\multiput{$\boldkey.$} at 24 1.8 /

\multiput{$\smallsq2$} at 23.8 -.4 *5 0 .2 /
\multiput{$\smallsq2$} at 24 -.2 *3 0 .2 /
\multiput{$\boldkey.$} at 23.9 -.2  24.1 -.2  /
\plot 23.9 -.2  24.1 -.2  /

\multiput{$\smallsq2$} at 24.9 .6 *5  0 .2 /
\multiput{$\boldkey.$} at 25 .8 /

\multiput{$\smallsq2$} at 26.8 .6 *6 0 .2 /
\multiput{$\smallsq2$} at 27 .8 *1 0 .2 /
\multiput{$\boldkey.$} at 26.9 1  27.1 1 /
\plot  26.9 1  27.1 1 /

\multiput{$\smallsq2$} at 28.9 .8 *2 0 .2 /
\multiput{$\boldkey.$} at 29 .8 /

\arr{0.3 0.3} {0.7 0.7}
\arr{1.3 0.7} {1.7 0.3}
\arr{2.3 0.3} {2.7 0.7}
\arr{3.3 0.7} {3.7 0.3}
\arr{4.3 0.3} {4.7 0.7}
\arr{5.3 0.7} {5.7 0.3}
\arr{6.3 0.3} {6.7 0.7}
\arr{7.3 0.7} {7.7 0.3}
\arr{5.3 -.7} {5.7 -.3}
\arr{6.3 -.3} {6.7 -.7}
\arr{5.3 1.3} {5.7 1.7}
\arr{6.3 1.7} {6.7 1.3}
\arr{7.3 1.3} {7.7 1.7}
\arr{8.3 1.7} {8.7 1.3}
\arr{9.3 1.3} {9.7 1.7}
\arr{10.3 1.7} {10.7 1.3}
\arr{11.3 1.3} {11.7 1.7}
\arr{12.3 1.7} {12.7 1.3}
\arr{13.3 1.3} {13.7 1.7}
\arr{14.3 1.7} {14.7 1.3}
\arr{15.3 1.3} {15.7 1.7}
\arr{16.3 1.7} {16.7 1.3}
\arr{17.3 1.3} {17.7 1.7}
\arr{18.3 1.7} {18.7 1.3}
\arr{19.3 1.3} {19.7 1.7}
\arr{20.3 1.7} {20.7 1.3}
\arr{21.3 1.3} {21.7 1.7}
\arr{22.3 1.7} {22.7 1.3}
\arr{23.3 1.3} {23.7 1.7}
\arr{24.3 1.7} {24.7 1.3}
\arr{12.3 2.3} {12.7 2.7}
\arr{13.3 2.7} {13.7 2.3}
\arr{12.3  .3} {12.7  .7}
\arr{13.3  .7} {13.7  .3}
\arr{22.3 0.3} {22.7 0.7}
\arr{23.3 0.7} {23.7 0.3}
\arr{24.3 0.3} {24.7 0.7}
\arr{25.3 0.7} {25.7 0.3}
\arr{26.3 0.3} {26.7 0.7}
\arr{27.3 0.7} {27.7 0.3}
\arr{28.3 0.3} {28.7 0.7}
\arr{29.3 0.7} {29.7 0.3}
\arr{23.3 -.7} {23.7 -.3}
\arr{24.3 -.3} {24.7 -.7}

\setdots<2pt>
\plot -.5 1  0.5 1 /
\plot 1.5 1  2.5 1 /
\plot 3.5 1  4.5 1 /
\plot 6.5 2  7.5 2 /
\plot 8.5 2  9.5 2 /
\plot 10.5 2  11.5 2 /
\plot 14.5 2  15.5 2 /
\plot 16.5 2  17.5 2 /
\plot 18.5 2  19.5 2 /
\plot 20.5 2  21.5 2 /
\plot 22.5 2  23.5 2 /
\plot 25.5 1  26.5 1 /
\plot 27.5 1  28.5 1 /
\plot 29.5 1  30.5 1 /
\setshadegrid span <1.5mm>
\vshade   -.5 -1 1 <,z,,> 
          5   -1 1 <z,z,,>
          6   -1 2 <z,z,,>
          7   -1 2 <z,z,,>
          8    0 2 <z,z,,>
          12   0 2 <z,z,,>
          13   0 3 <z,z,,>
          14   0 2 <z,z,,>
          22   0 2 <z,z,,>
          23  -1 2 <z,z,,>
          24  -1 2 <z,z,,>
          25  -1 1 <z,,,>
         30.5 -1 1 /
\endpicture} 

\def\fivetubethreeseven{\beginpicture\setcoordinatesystem units <6.3mm,6.3mm>
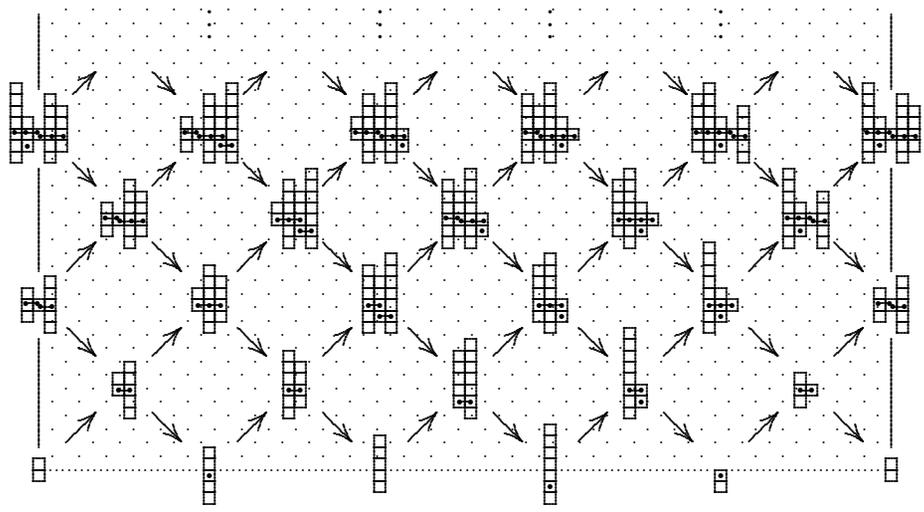
\put{The Tube of Circumference Five in $\Cal S_3(k[T]/T^7)$} at 7.5 -2
\put{} at 15.5 8.5
\multiput{$\smallsq2$} at -.1 -.1 *1 0 .2 /

\multiput{$\smallsq2$} at 2.9 -.5 *4  0 .2 /
\multiput{$\boldkey.$} at 3 -.1 /

\multiput{$\smallsq2$} at 5.9 -.3 *4 0 .2 /

\multiput{$\smallsq2$} at 8.9 -.5 *6  0 .2 /
\multiput{$\boldkey.$} at 9 -.3 /

\multiput{$\smallsq2$} at 11.9 -.3 *1  0 .2 /
\multiput{$\boldkey.$} at 12 -.1 /

\multiput{$\smallsq2$} at 14.9 -.1 *1  0 .2 /

\multiput{$\smallsq2$} at 1.3 1.4 *1  0 .2 /
\multiput{$\smallsq2$} at 1.5 1 *4  0 .2 /
\multiput{$\boldkey.$} at 1.4 1.4  1.6 1.4 /
\plot 1.4 1.4  1.6 1.4 /

\multiput{$\smallsq2$} at 4.3 1 *5  0 .2 /
\multiput{$\smallsq2$} at 4.5 1.2 *3  0 .2 /
\multiput{$\boldkey.$} at 4.4 1.4  4.6 1.4 /
\plot 4.4 1.4  4.6 1.4 /

\multiput{$\smallsq2$} at 7.3 1.2 *4  0 .2 /
\multiput{$\smallsq2$} at 7.5 1 *6  0 .2 /
\multiput{$\boldkey.$} at 7.4 1.2  7.6 1.2 /
\plot 7.4 1.2  7.6 1.2 /

\multiput{$\smallsq2$} at 10.3 1 *7  0 .2 /
\multiput{$\smallsq2$} at 10.5 1.2 *1  0 .2 /
\multiput{$\boldkey.$} at 10.4 1.4  10.6 1.4  10.6 1.2 /
\plot 10.4 1.4  10.6 1.4 /

\multiput{$\smallsq2$} at 13.3 1.2 *2  0 .2 /
\multiput{$\smallsq2$} at 13.5 1.4 /
\multiput{$\boldkey.$} at 13.4 1.4  13.6 1.4 /
\plot 13.4 1.4  13.6 1.4 /

\multiput{$\smallsq2$} at -.3 2.7 *2  0 .2 /
\multiput{$\smallsq2$} at -.1 2.9 /
\multiput{$\smallsq2$} at .1 2.5 *4  0 .2 /
\multiput{$\boldkey.$} at -.23 2.93  -.03 2.93  .03 2.87  .23 2.87 /
\plot -.23 2.93  -.03 2.93 /
\plot  .03 2.87  .23 2.87 /

\multiput{$\smallsq2$} at 2.7 2.9 *1  0 .2 /
\multiput{$\smallsq2$} at 2.9 2.5 *5  0 .2 /
\multiput{$\smallsq2$} at 3.1 2.7 *3  0 .2 /
\multiput{$\boldkey.$} at 2.8 2.9  3 2.9  3.2 2.9 /
\plot 2.8 2.9  3 2.9  3.2 2.9 /

\multiput{$\smallsq2$} at 5.7 2.5 *5  0 .2 /
\multiput{$\smallsq2$} at 5.9 2.7 *3  0 .2 /
\multiput{$\smallsq2$} at 6.1 2.5 *6  0 .2 /
\multiput{$\boldkey.$} at 5.8 2.9  6 2.9  6 2.7  6.2 2.7 /
\plot 5.8 2.9  6 2.9 /
\plot 6 2.7  6.2 2.7 /

\multiput{$\smallsq2$} at 8.7 2.7 *4  0 .2 /
\multiput{$\smallsq2$} at 8.9 2.5 *6  0 .2 /
\multiput{$\smallsq2$} at 9.1 2.7 *1  0 .2 /
\multiput{$\boldkey.$} at 8.8 2.9  9 2.9  9.2 2.9  9.2 2.7 /
\plot 8.8 2.9  9 2.9  9.2 2.9 /

\multiput{$\smallsq2$} at 11.7 2.5 *7  0 .2 /
\multiput{$\smallsq2$} at 11.9 2.7 *2  0 .2 /
\multiput{$\smallsq2$} at 12.1 2.9 /
\multiput{$\boldkey.$} at 11.8 2.9  12 2.9  12.2 2.9  12 2.7 /
\plot 11.8 2.9  12 2.9  12.2 2.9 /

\multiput{$\smallsq2$} at 14.7 2.7 *2  0 .2 /
\multiput{$\smallsq2$} at 14.9 2.9 /
\multiput{$\smallsq2$} at 15.1 2.5 *4  0 .2 /
\multiput{$\boldkey.$} at 14.77 2.93  14.97 2.93 15.03 2.87  15.23 2.87 /
\plot 14.77 2.93  14.97 2.93 /
\plot 15.03 2.87  15.23 2.87 /

\multiput{$\smallsq2$} at 1.1 4.2 *2  0 .2 /
\multiput{$\smallsq2$} at 1.3 4.4 /
\multiput{$\smallsq2$} at 1.5 4 *5  0 .2 /
\multiput{$\smallsq2$} at 1.7 4.2 *3  0 .2 /
\multiput{$\boldkey.$} at 1.17 4.43  1.37 4.43  
                    1.43 4.37  1.63 4.37  1.83 4.37 /
\plot 1.17 4.43  1.37 4.43  /
\plot 1.43 4.37  1.63 4.37  1.83 4.37 /

\multiput{$\smallsq2$} at 4.1 4.4 *1  0 .2 /
\multiput{$\smallsq2$} at 4.3 4 *5  0 .2 /
\multiput{$\smallsq2$} at 4.5 4.2 *3  0 .2 /
\multiput{$\smallsq2$} at 4.7 4 *6  0 .2 /
\multiput{$\boldkey.$} at 4.2 4.4  4.4 4.4  4.6 4.4  4.6 4.2  4.8 4.2  /
\plot 4.2 4.4  4.4 4.4  4.6 4.4 /
\plot 4.6 4.2  4.8 4.2  /

\multiput{$\smallsq2$} at 7.1 4 *5  0 .2 /
\multiput{$\smallsq2$} at 7.3 4.2 *3  0 .2 /
\multiput{$\smallsq2$} at 7.5 4 *6  0 .2 /
\multiput{$\smallsq2$} at 7.7 4.2 *1  0 .2 /
\multiput{$\boldkey.$} at 7.17 4.43  7.37 4.43  7.43 4.37  7.63 4.37  7.83 4.37  7.8 4.2 /
\plot 7.17 4.43  7.37 4.43 /
\plot 7.43 4.37  7.63 4.37  7.83 4.37 /

\multiput{$\smallsq2$} at 10.1 4.2 *4  0 .2 /
\multiput{$\smallsq2$} at 10.3 4 *6  0 .2 /
\multiput{$\smallsq2$} at 10.5 4.2 *2  0 .2 /
\multiput{$\smallsq2$} at 10.7 4.4 /
\multiput{$\boldkey.$} at 10.2 4.4  10.4 4.4  10.6 4.4  10.8 4.4  10.6 4.2 /
\plot 10.2 4.4  10.8 4.4 /

\multiput{$\smallsq2$} at 13.1 4 *6  0 .2 /
\multiput{$\smallsq2$} at 13.3 4.2 *2  0 .2 /
\multiput{$\smallsq2$} at 13.5 4.4 /
\multiput{$\smallsq2$} at 13.7 4 *4  0 .2 /
\multiput{$\boldkey.$} at 13.17 4.43  13.37 4.43  13.57 4.43  
                        13.63 4.37  13.83 4.37  13.4 4.2 /
\plot 13.17 4.43  13.37 4.43  13.57 4.43  /
\plot   13.63 4.37  13.83 4.37 /

\multiput{$\smallsq2$} at -.5 5.5 *6  0 .2 /
\multiput{$\smallsq2$} at -.3 5.7 *2  0 .2 /
\multiput{$\smallsq2$} at -.1 5.9 /
\multiput{$\smallsq2$} at .1 5.5 *5  0 .2 /
\multiput{$\smallsq2$} at .3 5.7 *3  0 .2 /
\multiput{$\boldkey.$} at -.43 5.93  -.23 5.93  -.03 5.93  .03 5.87 
                        .23 5.87  .43 5.87  -.2 5.7 /
\plot -.43 5.93  -.03 5.93 /
\plot  .03 5.87  .43 5.87  /

\multiput{$\smallsq2$} at 2.5 5.7 *2  0 .2 /
\multiput{$\smallsq2$} at 2.7 5.9 /
\multiput{$\smallsq2$} at 2.9 5.5 *5  0 .2 /
\multiput{$\smallsq2$} at 3.1 5.7 *3  0 .2 /
\multiput{$\smallsq2$} at 3.3 5.5 *6  0 .2 /
\multiput{$\boldkey.$} at 2.57 5.93  2.77 5.93  
                        2.83 5.87  3.03 5.87  3.23 5.87  3.2 5.7  3.4 5.7 /
\plot 2.57 5.93  2.77 5.93 /
\plot  2.83 5.87  3.23 5.87 /
\plot    3.2 5.7  3.4 5.7 /

\multiput{$\smallsq2$} at 5.5 5.9 *1  0 .2 /
\multiput{$\smallsq2$} at 5.7 5.5 *5  0 .2 /
\multiput{$\smallsq2$} at 5.9 5.7 *3  0 .2 /
\multiput{$\smallsq2$} at 6.1 5.5 *6  0 .2 /
\multiput{$\smallsq2$} at 6.3 5.7 *1  0 .2 /
\multiput{$\boldkey.$} at 5.57 5.93  5.77 5.93  5.97 5.93 
                        6.03 5.87  6.23 5.87  6.43 5.87  6.4 5.7 /
\plot  5.57 5.93  5.77 5.93  5.97 5.93 /
\plot  6.03 5.87  6.23 5.87  6.43 5.87 /

\multiput{$\smallsq2$} at 8.5 5.5 *5  0 .2 /
\multiput{$\smallsq2$} at 8.7 5.7 *3  0 .2 /
\multiput{$\smallsq2$} at 8.9 5.5 *6  0 .2 /
\multiput{$\smallsq2$} at 9.1 5.7 *2  0 .2 /
\multiput{$\smallsq2$} at 9.3 5.9 /
\multiput{$\boldkey.$} at 8.57 5.93  8.77 5.93
                        8.83 5.87  9.03 5.87  9.23 5.87  9.43 5.87  9.2 5.7 /
\plot  8.57 5.93  8.77 5.93 /
\plot  8.83 5.87  9.43 5.87 /

\multiput{$\smallsq2$} at 11.5 5.7 *4  0 .2 /
\multiput{$\smallsq2$} at 11.7 5.5 *6  0 .2 /
\multiput{$\smallsq2$} at 11.9 5.7 *2  0 .2 /
\multiput{$\smallsq2$} at 12.1 5.9 /
\multiput{$\smallsq2$} at 12.3 5.5 *4  0 .2 /
\multiput{$\boldkey.$} at 11.57 5.93  11.77 5.93  11.97 5.93  12.17 5.93
                        12.23 5.87  12.43 5.87  12 5.7 /
\plot 11.57 5.93  12.17 5.93 /
\plot   12.23 5.87  12.43 5.87 /

\multiput{$\smallsq2$} at 14.5 5.5 *6  0 .2 /
\multiput{$\smallsq2$} at 14.7 5.7 *2  0 .2 /
\multiput{$\smallsq2$} at 14.9 5.9 /
\multiput{$\smallsq2$} at 15.1 5.5 *5  0 .2 /
\multiput{$\smallsq2$} at 15.3 5.7 *3  0 .2 /
\multiput{$\boldkey.$} at 14.57 5.93  14.77 5.93  14.97 5.93
                        15.03 5.87  15.23 5.87  15.43 5.87  14.8 5.7 /
\plot  14.57 5.93  14.97 5.93 /
\plot  15.03 5.87  15.43 5.87 /
\arr{0.5 0.5}{1 1}
\arr{2 1} {2.5 0.5}
\arr{3.5 0.5}{4 1}
\arr{5 1} {5.5 0.5}
\arr{6.5 0.5}{7 1}
\arr{8 1} {8.5 0.5}
\arr{9.5 0.5}{10 1}
\arr{11 1} {11.5 0.5}
\arr{12.5 0.5}{13 1}
\arr{14 1} {14.5 0.5}
\arr{0.5 2.5}{1 2}
\arr{2 2} {2.5 2.5}
\arr{3.5 2.5}{4 2}
\arr{5 2} {5.5 2.5}
\arr{6.5 2.5}{7 2}
\arr{8 2} {8.5 2.5}
\arr{9.5 2.5}{10 2}
\arr{11 2} {11.5 2.5}
\arr{12.5 2.5}{13 2}
\arr{14 2} {14.5 2.5}
\arr{0.5 3.5}{1 4}
\arr{2 4} {2.5 3.5}
\arr{3.5 3.5}{4 4}
\arr{5 4} {5.5 3.5}
\arr{6.5 3.5}{7 4}
\arr{8 4} {8.5 3.5}
\arr{9.5 3.5}{10 4}
\arr{11 4} {11.5 3.5}
\arr{12.5 3.5}{13 4}
\arr{14 4} {14.5 3.5}
\arr{0.6 5.4}{1 5}
\arr{2 5} {2.4 5.4}
\arr{3.6 5.4}{4 5}
\arr{5 5} {5.4 5.4}
\arr{6.6 5.4}{7 5}
\arr{8 5} {8.4 5.4}
\arr{9.6 5.4}{10 5}
\arr{11 5} {11.4 5.4}
\arr{12.6 5.4}{13 5}
\arr{14 5} {14.4 5.4}
\arr{0.6 6.6}{1 7}
\arr{2 7} {2.4 6.6}
\arr{3.6 6.6}{4 7}
\arr{5 7} {5.4 6.6}
\arr{6.6 6.6}{7 7}
\arr{8 7} {8.4 6.6}
\arr{9.6 6.6}{10 7}
\arr{11 7} {11.4 6.6}
\arr{12.6 6.6}{13 7}
\arr{14 7} {14.4 6.6}
\plot 0 0.4  0 2.3 /
\plot 0 3.5 0 5.3 /
\plot 0 6.7 0 8 /
\plot 15 0.4  15 2.3 /
\plot 15 3.5 15 5.3 /
\plot 15 6.7 15 8 /
\setdots <2pt>
\plot 0.2 0  2.8 0 /
\plot 3.2 0  5.8 0 /
\plot 6.2 0  8.8 0 /
\plot 9.2 0  11.8 0 /
\plot 12.2 0  14.8 0 /

\multiput{$\vdots$} at 3 8  6 8  9 8  12 8 /
\setshadegrid span <1.5mm>
\vshade    0  0 8.5  <,,,> 
          15  0 8.5  /
\endpicture}

\def\threetubethreeseven{\beginpicture
\setcoordinatesystem units <6.3mm,6.3mm>
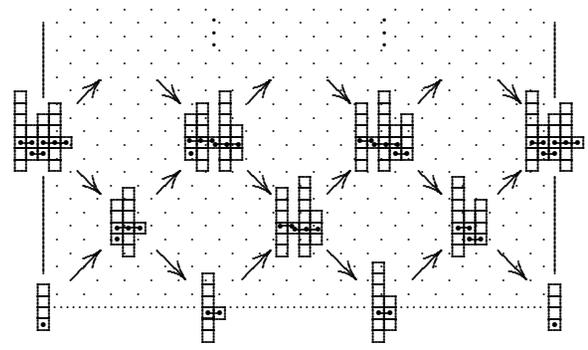
\put{The Tube of Circumference Three in $\Cal S_3(k[T]/T^7)$} at 4.5 -2
\put{} at 9.5 5.5
\multiput{$\smallsq2$} at -.1 -.3 *3   0 .2 /
\multiput{$\boldkey.$} at 0 -.3   /

\multiput{$\smallsq2$} at 2.8 -.5 *5   0 .2 /
\multiput{$\smallsq2$} at 3 -.1 /
\multiput{$\boldkey.$} at 2.9 -.1  3.1 -.1   /
\plot  2.9 -.1  3.1 -.1   /

\multiput{$\smallsq2$} at  5.8 -.5 *6  0 .2 /
\multiput{$\smallsq2$} at  6 -.3 *2  0 .2 /
\multiput{$\boldkey.$} at  5.9 -.1  6.1 -.1  /
\plot  5.9 -.1  6.1 -.1  /

\multiput{$\smallsq2$} at  8.9 -.3 *3  0 .2 /
\multiput{$\boldkey.$} at  9 -.3  /

\multiput{$\smallsq2$} at 1.2 1.2 *3   0 .2 /
\multiput{$\smallsq2$} at 1.4 1 *5   0 .2 /
\multiput{$\smallsq2$} at 1.6 1.4 /
\multiput{$\boldkey.$} at 1.3 1.4  1.5 1.4  1.7 1.4  1.3 1.2   /
\plot 1.3 1.4  1.7 1.4 /

\multiput{$\smallsq2$} at 4.1 1 *5   0 .2 /
\multiput{$\smallsq2$} at 4.3 1.4 /
\multiput{$\smallsq2$} at 4.5 1 *6   0 .2 /
\multiput{$\smallsq2$} at 4.7 1.2 *2   0 .2 /
\multiput{$\boldkey.$} at 4.17 1.43  4.37 1.43  4.43 1.37  4.63 1.37  4.83 1.37   /
\plot  4.17 1.43  4.37 1.43 /
\plot  4.43 1.37  4.83 1.37 /

\multiput{$\smallsq2$} at 7.2 1 *6   0 .2 /
\multiput{$\smallsq2$} at 7.4 1.2 *2   0 .2 /
\multiput{$\smallsq2$} at 7.6 1.2 *3   0 .2 /
\multiput{$\boldkey.$} at 7.3 1.4  7.5 1.4  7.5 1.2  7.7 1.2    /
\plot  7.3 1.4  7.5 1.4 /
\plot  7.5 1.2  7.7 1.2 /

\multiput{$\smallsq2$} at -.5 2.5 *6   0 .2 /
\multiput{$\smallsq2$} at -.3 2.7 *2   0 .2 /
\multiput{$\smallsq2$} at -.1 2.7 *3   0 .2 /
\multiput{$\smallsq2$} at .1 2.5 *5   0 .2 /
\multiput{$\smallsq2$} at .3 2.9 /
\multiput{$\boldkey.$} at -.4 2.9  -.2 2.9  -.2 2.7  0 2.7  0 2.9  .2 2.9  
                        .4 2.9   /
\plot  -.4 2.9  -.2 2.9 /
\plot  -.2 2.7  0 2.7 /
\plot  0 2.9  .4 2.9   /

\multiput{$\smallsq2$} at  2.5 2.7 *3  0 .2 /
\multiput{$\smallsq2$} at  2.7 2.5 *5  0 .2 /
\multiput{$\smallsq2$} at  2.9 2.9 /
\multiput{$\smallsq2$} at  3.1 2.5 *6  0 .2 /
\multiput{$\smallsq2$} at  3.3 2.7 *2  0 .2 /
\multiput{$\boldkey.$} at  2.57 2.93  2.77 2.93  2.97 2.93
                3.03 2.87  3.23 2.87 3.43 2.87  2.6 2.7  /
\plot   2.57 2.93  2.77 2.93  2.97 2.93 /
\plot  3.03 2.87  3.23 2.87 3.43 2.87 /

\multiput{$\smallsq2$} at  5.5 2.5 *5  0 .2 /
\multiput{$\smallsq2$} at  5.7 2.9 /
\multiput{$\smallsq2$} at  5.9 2.5 *6  0 .2 /
\multiput{$\smallsq2$} at  6.1 2.7 *2  0 .2 /
\multiput{$\smallsq2$} at  6.3 2.7 *3  0 .2 /
\multiput{$\boldkey.$} at  5.57 2.93  5.77 2.93
                5.83 2.87  6.03 2.87  6.23 2.87  6.2 2.7  6.4 2.7  /
\plot  5.57 2.93  5.77 2.93 /
\plot  5.83 2.87  6.03 2.87  6.23 2.87 /
\plot  6.2 2.7  6.4 2.7  /

\multiput{$\smallsq2$} at 8.5 2.5 *6   0 .2 /
\multiput{$\smallsq2$} at 8.7 2.7 *2   0 .2 /
\multiput{$\smallsq2$} at 8.9 2.7 *3   0 .2 /
\multiput{$\smallsq2$} at 9.1 2.5 *5   0 .2 /
\multiput{$\smallsq2$} at 9.3 2.9 /
\multiput{$\boldkey.$} at 8.6 2.9  8.8 2.9  8.8 2.7  9 2.7 
                        9 2.9  9.2 2.9  9.4 2.9    /
\plot  8.6 2.9  8.8 2.9 /
\plot  8.8 2.7  9 2.7 /
\plot   9 2.9  9.2 2.9  9.4 2.9    /
\arr{0.5 0.5}{1 1}
\arr{2 1} {2.5 0.5}
\arr{3.5 0.5}{4 1}
\arr{5 1} {5.5 0.5}
\arr{6.5 0.5}{7 1}
\arr{8 1} {8.5 0.5}
\arr{0.6 2.4}{1 2}
\arr{2 2} {2.4 2.4}
\arr{3.6 2.4}{4 2}
\arr{5 2} {5.4 2.4}
\arr{6.6 2.4}{7 2}
\arr{8 2} {8.4 2.4}
\arr{0.6 3.6}{1 4}
\arr{2 4} {2.4 3.6}
\arr{3.6 3.6}{4 4}
\arr{5 4} {5.4 3.6}
\arr{6.6 3.6}{7 4}
\arr{8 4} {8.4 3.6}
\plot 0 .6  0 2.3 /
\plot 0 3.7 0 5 /
\plot 9 .6  9 2.3 /
\plot 9 3.7 9 5 /
\setdots <2pt>
\plot 0.2 0  2.7 0 /
\plot 3.3 0  5.7 0 /
\plot 6.3 0  8.8 0 /

\multiput{$\vdots$} at 3 5  6 5 /
\setshadegrid span <1.5mm>
\vshade    0  0 5.5  <,,,> 
           9  0 5.5  /
\endpicture}

\def\twotubethreeseven{\beginpicture
\setcoordinatesystem units <6.3mm,6.3mm>
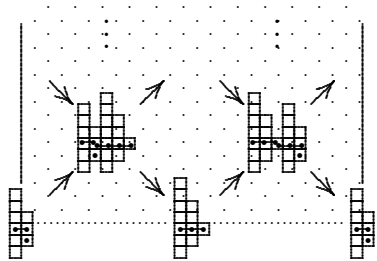
\put{The Tube of Circumference Two in $\Cal S_3(k[T]/T^7)$} at 3 -1.5

\put{} at 6.5 4
\multiput{$\smallsq2$} at -.2 -.5 *5  0 .2 /
\multiput{$\smallsq2$} at 0 -.3 *2  0 .2 /
\multiput{$\boldkey.$} at -.1 -.1  .1 -.1  .1 -.3  /
\plot -.1 -.1  .1 -.1 /

\multiput{$\smallsq2$} at 2.7 -.5  *6  0 .2 /
\multiput{$\smallsq2$} at 2.9 -.3 *3  0 .2 /
\multiput{$\smallsq2$} at 3.1 -.1 /
\multiput{$\boldkey.$} at 2.8 -.1  3 -.1  3.2 -.1  /
\plot 2.8 -.1  3 -.1  3.2 -.1  /

\multiput{$\smallsq2$} at 5.8 -.5 *5  0 .2 /
\multiput{$\smallsq2$} at 6 -.3 *2  0 .2 /
\multiput{$\boldkey.$} at 5.9 -.1  6.1 -.1  6.1 -.3   /
\plot 5.9 -.1  6.1 -.1 /

\multiput{$\smallsq2$} at 1 1 *5  0 .2 /
\multiput{$\smallsq2$} at 1.2 1.2 *2  0 .2 /
\multiput{$\smallsq2$} at 1.4 1 *6  0 .2 /
\multiput{$\smallsq2$} at 1.6 1.2 *3  0 .2 /
\multiput{$\smallsq2$} at 1.8 1.4 /
\multiput{$\boldkey.$} at 1.07 1.43  1.27 1.43  
                1.33 1.37  1.53 1.37  1.73 1.37  1.93 1.37  1.3 1.2  /
\plot 1.07 1.43  1.27 1.43 /
\plot  1.33 1.37   1.93 1.37 /

\multiput{$\smallsq2$} at 4 1 *6  0 .2 /
\multiput{$\smallsq2$} at 4.2 1.2 *3  0 .2 /
\multiput{$\smallsq2$} at 4.4 1.4 /
\multiput{$\smallsq2$} at 4.6 1 *5  0 .2 /
\multiput{$\smallsq2$} at 4.8 1.2 *2  0 .2 /
\multiput{$\boldkey.$} at 4.07 1.43  4.27 1.43  4.47 1.43
                        4.53 1.37  4.73 1.37  4.93 1.37  4.9 1.2  /
\plot 4.07 1.43  4.27 1.43  4.47 1.43 /
\plot   4.53 1.37  4.73 1.37  4.93 1.37 /
\arr{0.5 0.5}{.9 .9}
\arr{2.1 .9} {2.5 0.5}
\arr{3.5 0.5}{3.9 .9}
\arr{5.1 .9} {5.5 0.5}
\arr{0.5 2.5}{.9 2.1}
\arr{2.1 2.1} {2.5 2.5}
\arr{3.5 2.5}{3.9 2.1}
\arr{5.1 2.1} {5.5 2.5}

\plot 0 0.7  0 3.5 /
\plot 6 0.7  6 3.5 /
\setdots <2pt>
\plot .3 0  2.6 0 /
\plot 3.4 0  5.7 0 /
\multiput{$\vdots$} at 1.5 3.5  4.5 3.5 /
\setshadegrid span <1.5mm>
\vshade    0  0 4  <,,,> 
           6  0 4  /
\endpicture}

\topinsert\noindent
\rotatebox{90}{
$$\beginpicture
\setcoordinatesystem units <1cm,1cm>
\put{} at -10 0
\put{} at 10 0
\put{\boundarythreeseven} at 0 7
\put{\fivetubethreeseven} at -4.5 0
\put{\threetubethreeseven} at 6 -.94
\endpicture$$}
\endinsert

\smallskip
The two nonstable orbits in $\Cal U$ and the mouths of the big tubes 
are pictured below. 

\smallskip The indecomposables on the mouth of the homogeneous tubes
which do not involve any field extensions have the following type
$$
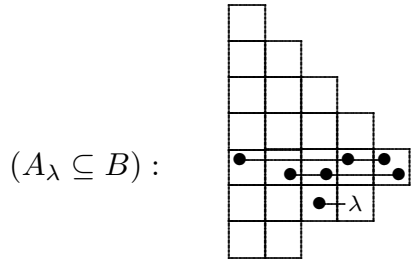
\beginpicture 
\setcoordinatesystem units <.4cm,.4cm>
\put{$(A_\lambda\sub B):$} at -4 2
\multiput{\sq} at 0 6  0 5  0 4  0 3  0 2  0 1  0 0  1 5  1 4  1 3  1 2  1 1  1 0  
                  2 4  2 3  2 2  2 1  3 3  3 2  3 1  4 2 /
\put{$\bullet$} at 0.3 2.2
\put{$\bullet$} at 3.3 2.2
\put{$\bullet$} at 4.3 2.2
\put{$\bullet$} at 1.7 1.8
\put{$\bullet$} at 2.7 1.8 
\put{$\bullet$} at 4.7 1.8 
\put{$\bullet$} at 2.5 1 
\put{$\ssize \lambda$} at 3.5 1 
\plot 0.3 2.2  4.3 2.2 /
\plot 1.7 1.8  4.7 1.8 /
\plot 2.5 1  3.2 1 /
\endpicture
$$
where $\lambda\in k\backslash\{0,-1\}$.  
For $\lambda=\infty$ we obtain the
last module in the fifth row in the tube of circumference five;
for $\lambda=0$, $(A_\lambda\sub B)$ is the second module in the
third row in the tube of circumference three; and for $\lambda=-1$
one obtains the first module in the second row in the tube of 
circumference two below.  The isomorphism is as follows.
$$
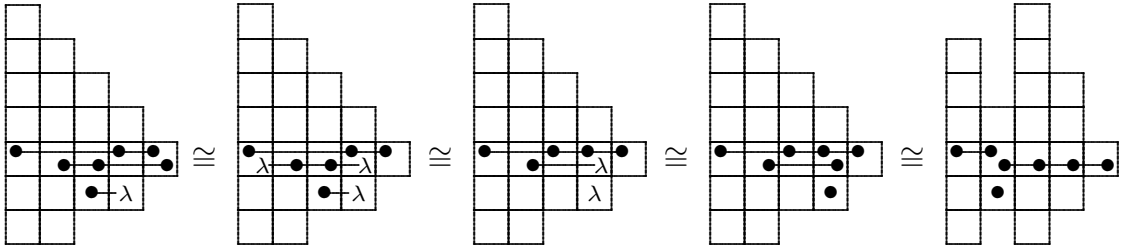
\beginpicture 
\setcoordinatesystem units <.38cm,.38cm>
\multiput{\sq} at 0 6  0 5  0 4  0 3  0 2  0 1  0 0  1 5  1 4  1 3  1 2  1 1  1 0  
                  2 4  2 3  2 2  2 1  3 3  3 2  3 1  4 2 /
\put{$\bullet$} at 0.3 2.2
\put{$\bullet$} at 3.3 2.2
\put{$\bullet$} at 4.3 2.2
\put{$\bullet$} at 1.7 1.8
\put{$\bullet$} at 2.7 1.8 
\put{$\bullet$} at 4.7 1.8 
\put{$\bullet$} at 2.5 1 
\put{$\ssize \lambda$} at 3.5 1 
\plot 0.3 2.2  4.3 2.2 /
\plot 1.7 1.8  4.7 1.8 /
\plot 2.5 1  3.2 1 /
\endpicture
\;
\beginpicture 
\setcoordinatesystem units <.38cm,.38cm>
\put{$\cong$} at -1 2
\put{} at -1.5 0
\multiput{\sq} at 0 6  0 5  0 4  0 3  0 2  0 1  0 0  1 5  1 4  1 3  1 2  1 1  1 0  
                  2 4  2 3  2 2  2 1  3 3  3 2  3 1  4 2 /
\put{$\bullet$} at 0.3 2.2
\put{$\bullet$} at 3.3 2.2
\put{$\bullet$} at 4.3 2.2
\put{$\ssize\lambda$} at 0.7 1.8
\put{$\bullet$} at 1.7 1.8
\put{$\bullet$} at 2.7 1.8 
\put{$\ssize\lambda$} at 3.7 1.8
\put{$\bullet$} at 2.5 1 
\put{$\ssize \lambda$} at 3.5 1 
\plot 0.3 2.2  4.3 2.2 /
\plot 0.9 1.8  3.5  1.8  /
\plot 2.5 1  3.2 1 /
\endpicture
\;
\beginpicture 
\setcoordinatesystem units <.38cm,.38cm>
\put{$\cong$} at -1 2
\put{} at -2 0
\multiput{\sq} at 0 6  0 5  0 4  0 3  0 2  0 1  0 0  1 5  1 4  1 3  1 2  1 1  1 0  
                  2 4  2 3  2 2  2 1  3 3  3 2  3 1  4 2 /
\put{$\bullet$} at 0.3 2.2
\put{$\bullet$} at 2.3 2.2
\put{$\bullet$} at 3.3 2.2
\put{$\bullet$} at 4.3 2.2
\put{$\bullet$} at 1.7 1.8
\put{$\ssize\lambda$} at 3.7 1.8
\put{$\ssize \lambda$} at 3.5 1 
\plot 0.3 2.2  4.3 2.2 /
\plot 1.7 1.8  3.5  1.8  /
\endpicture
\;
\beginpicture 
\setcoordinatesystem units <.38cm,.38cm>
\put{$\cong$} at -1 2
\put{} at -2 0
\multiput{\sq} at 0 6  0 5  0 4  0 3  0 2  0 1  0 0  1 5  1 4  1 3  1 2  1 1  1 0  
                  2 4  2 3  2 2  2 1  3 3  3 2  3 1  4 2 /
\put{$\bullet$} at 0.3 2.2
\put{$\bullet$} at 2.3 2.2
\put{$\bullet$} at 3.3 2.2
\put{$\bullet$} at 4.3 2.2
\put{$\bullet$} at 1.7 1.8
\put{$\bullet$} at 3.7 1.8
\put{$\bullet$} at 3.5 1 
\plot 0.3 2.2  4.3 2.2 /
\plot 1.7 1.8  3.7  1.8  /
\endpicture
\;
\beginpicture 
\setcoordinatesystem units <.38cm,.38cm>
\put{$\cong$} at -1 2
\put{} at -2 0
\multiput{\sq} at 0 5  0 4  0 3  0 2  0 1  0 0  1 3  1 2  1 1   
                  2 6  2 5  2 4  2 3  2 2  2 1  2 0  3 4  3 3  3 2  3 1  4 2 /
\put{$\bullet$} at 0.3 2.2
\put{$\bullet$} at 1.3 2.2
\put{$\bullet$} at 1.7 1.8
\put{$\bullet$} at 2.7 1.8 
\put{$\bullet$} at 3.7 1.8 
\put{$\bullet$} at 4.7 1.8 
\put{$\bullet$} at 1.5 1 
\plot 0.3 2.2  1.3 2.2 /
\plot 1.7 1.8  4.7 1.8 /
\endpicture
$$
(The first isomorphism is given by subtracting the first generator
of the subgroup from the second; if the big group has generators
$x_7,x_6,x_4,x_3,x_1$, then the second isomorphism is obtained by
replacing $x_6$ by $x_6'=x_6+\lambda Tx_7$ and $x_3$ by 
$x_3'=x_3+\lambda^{-1}Tx_4$; for 
the third isomorphism we multiply the second and the third generator
of the subgroup by $\lambda^{-1}$ and replace 
$x_6'$ by $x_6''=\lambda^{-1}x_6'$;
the last isomorphism is just a rearrangement of the columns.)

$$\twotubethreeseven$$

In particular, the objects $(A_\lambda\sub B)$ pictured above 
are indecomposable and pairwise nonisomorphic for $\lambda\in k\cup\{\infty\}$.
We deduce the corresponding information about the subgroup embeddings
$(C_\lambda\sub D)$ studied in the introduction.

\medskip\noindent{\bf Corollary 15.} {\it 
For $\lambda\in\{0,\ldots,p-1\}$, the subgroup embeddings
$(C_\lambda\sub D)$ in the introduction are indecomposable and pairwise
nonisomorphic objects in $\Cal S_3(\Bbb Z/p^7)$.}

\smallskip\noindent{\it Proof:\/} Consider the layer functors
$L_i\: \Cal S_3(\Bbb Z/p^7)\to \mod\Bbb Z/p^7$ defined by
$$ \eqalign{L_2(A\sub B) &= (p^4B\cap p^{-2}0_B)+(p^2B\cap p^{-1}0_B)
\T{and} \cr
L_i(A\sub B) & =p^{2-i}(L_2(A\sub B))\qquad\qquad\qquad\T{for $i\neq 2$.}}$$
For each $i$, the $i$-th layer $L_i=L_i(A\sub B)$ is a subgroup
of $B$, and the sequence $(L_i)_i$ defines a filtration for $B$,
as follows.
$$
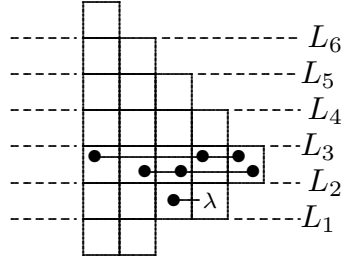
\beginpicture 
\setcoordinatesystem units <.4cm,.4cm>
\multiput{\sq} at 0 6  0 5  0 4  0 3  0 2  0 1  0 0  1 5  1 4  1 3  1 2  1 1  1 0  
                  2 4  2 3  2 2  2 1  3 3  3 2  3 1  4 2 /
\put{$\bullet$} at 0.3 2.2
\put{$\bullet$} at 3.3 2.2
\put{$\bullet$} at 4.3 2.2
\put{$\bullet$} at 1.7 1.8
\put{$\bullet$} at 2.7 1.8 
\put{$\bullet$} at 4.7 1.8 
\put{$\bullet$} at 2.5 1 
\put{$\ssize \lambda$} at 3.5 1 
\plot 0.3 2.2  4.3 2.2 /
\plot 1.7 1.8  4.7 1.8 /
\plot 2.5 1  3.2 1 /
\setdashes <2pt>
\plot -2 0.5  -.2 0.5 /
\plot 4.2 0.5  6 0.5 /
\plot -2 1.5  -.2 1.5 /
\plot 5.2 1.5  6 1.5 /
\plot -2 2.5  -.2 2.5 /
\plot 5.2 2.5  6 2.5 /
\plot -2 3.5  -.2 3.5 /
\plot 4.2 3.5  6 3.5 /
\plot -2 4.5  -.2 4.5 /
\plot 3.2 4.5  6 4.5 /
\plot -2 5.5  -.2 5.5 /
\plot 2.2 5.5  6 5.5 /
\put{$L_1$} at 6.5 0.5 
\put{$L_2$} at 6.7 1.5 
\put{$L_3$} at 6.5 2.5 
\put{$L_4$} at 6.7 3.5 
\put{$L_5$} at 6.5 4.5 
\put{$L_6$} at 6.7 5.5 
\endpicture
$$
The vectorspaces $M_i= L_i/L_{i-1}$ and 
$M_i'=A\cap L_i/A\cap L_{i-1}$ together with the linear maps 
obtained from the multiplication by $p$
define a $\Bbb Z/p$-linear representation
of $\widetilde Q$. In fact, the composition of this construction 
with the covering functor gives rise to an additive functor
$\Cal S_3(\Bbb Z/p^7)\to \Cal S_3((\Bbb Z/p)[T]/T^7)$
which maps the subgroup embeddings $(C_\lambda\sub D)$ from the introduction to
the above objects $(A_\lambda\sub B)$ in $ \Cal S_3((\Bbb Z/p)[T]/T^7)$.
Since the objects 
$(A_\lambda\sub B)$ are indecomposable and pairwise nonisomorphic
for $\lambda\in\Bbb Z/p$, so are the 
subgroup embeddings $(C_\lambda\sub D)$. \qed

\bigskip
\centerline{\sc The Cases $m=3$, $n\geq 8$.}

\medskip
For $m\geq 4$ and $n\geq 7$ or for $m\geq 3$ and $n\geq 8$, the 
categories $\Cal S_m(n)$ and $\Cal S_m(k[T]/T^n)$ have wild 
representation type.  In [8], it is shown that $\Cal S_4(k[T]/T^7)$
in fact is controlled $k$-wild with a single control object. 
Here we show the following result.

\medskip\noindent{\bf Proposition 16.} {\it
For $n\geq 8$, the categories $\Cal S_3(n)$ and $\Cal S_3(k[T]/T^n)$
have wild representation type. 
Moreover, the category $\Cal S_3(n)$ is even strictly wild 
in the sense that any finite dimensional $k$-algebra can be
realized as the endomorphism ring of an object in $\Cal S_3(n)$.  }

\smallskip\noindent{\it Proof:\/}
The algebra $B_0$ given by the quiver
$$
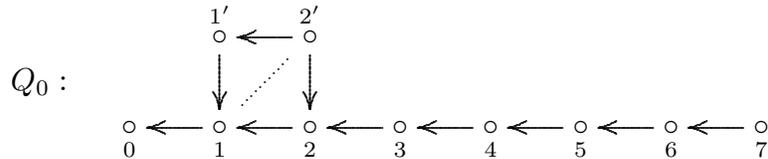
\hbox{\beginpicture
\setcoordinatesystem units <0.5cm,0.5cm>
\put{} at 0 -0.55
\put{} at 12 2.55
\put{$Q_0:$} at -2 1
\put{$\circ$} at 0 0
\put{$\circ$} at 2 0
\put{$\circ$} at 4 0
\put{$\circ$} at 6 0
\put{$\circ$} at 8 0
\put{$\circ$} at 10 0
\put{$\circ$} at 12 0
\put{$\circ$} at 14 0

\put{$\circ$} at 2 2
\put{$\circ$} at 4 2
\arr{1.6 0}{0.4 0}
\arr{3.6 0}{2.4 0}
\arr{5.6 0}{4.4 0}
\arr{7.6 0}{6.4 0}
\arr{9.6 0}{8.4 0}
\arr{11.6 0}{10.4 0}
\arr{13.6 0}{12.4 0}

\arr{3.6 2}{2.4 2}

\arr{2 1.6}{2 0.4}
\arr{4 1.6}{4 0.4}


\put{$\ssize 1'$} at  2 2.5
\put{$\ssize 2'$} at  4 2.5

\put{$\ssize 0$} at 0 -.5
\put{$\ssize 1$} at 2 -.5
\put{$\ssize 2$} at 4 -.5
\put{$\ssize 3$} at 6 -.5
\put{$\ssize 4$} at 8 -.5
\put{$\ssize 5$} at 10 -.5
\put{$\ssize 6$} at 12 -.5
\put{$\ssize 7$} at 14 -.5

\setdots<2pt>
\plot 3.5 1.5  2.5 0.5 /
\endpicture}
$$
with the commutativity relation as indicated 
is tilted from a wild
hereditary algebra of type $\widetilde{\widetilde{ \Bbb E}}_8$ 
by using a preprojective
tilting module, and hence has wild representation type.
Indeed, the tilting functor is defined by forming the pushout of the
commutative square.
The simple $B_0$-modules $S_1'$ and $S_2'$ occur 
in the preinjective component and
hence the category $\Cal S_3(n)$ contains as a full subcategory the category
of regular $B_0$-modules (which is equivalent to the category of regular
$k\widetilde{\widetilde{\Bbb E}}_8$-modules).
Since any finite dimensional $k$-algebra can be realized as the endomorphism
ring of some regular $B_0$-module, 
it can also be realized as 
the endomorphism ring of an object in $\Cal S_3(n)$. \qed

\bigskip
\centerline{\sc The Case $m=4$, $n=6$.}

\nopagebreak
\medskip
The category $\Cal S_4(k[T]/T^6)$ has tame infinite representation type,
and the Auslander-Reiten quiver consists of a family of tubes and
a connecting component, similar to the situation in  $\Cal S_3(k[T]/T^7)$.
We also consider $\Cal S_4(k[T]/T^6)$ as a subcategory of $\Cal S(k[T]/T^6)$.

\smallskip
First we investigate the covering category $\Cal S_4(6)$. 
We will see that the 
Aus\-lan\-der-\-Rei\-ten quiver has the following
global structure.

$$
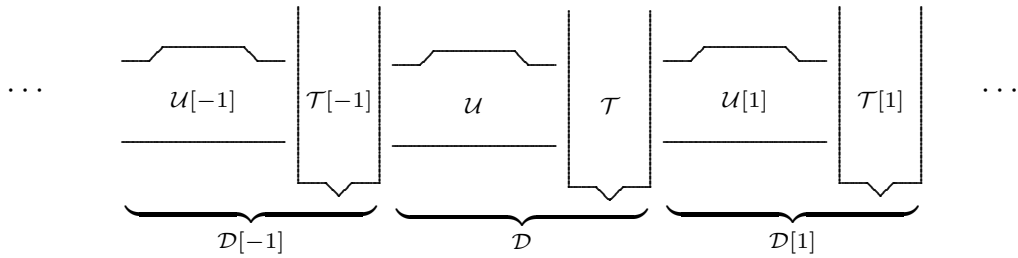
\hbox{\beginpicture
\setcoordinatesystem units <0.6cm,0.6cm>
\put{\beginpicture
\setcoordinatesystem units <0.3cm,0.3cm>
\put{} at 0 0
\put{} at 10 7 
\put{$\underbrace{\phantom{mmmmmmmmm}}_{\Cal D[-1]}$} at 5.25 -1

\plot 0.5 2  6.5 2 /
\plot 0.5 5  1.5 5  2 5.5  5 5.5  5.5 5  6.5 5 /
\plot 7 7  7 0.5  8 0.5  8.5 0  9 0.5  10 0.5  10 7 /

\put{$\ssize \Cal T[-1]$} at 8.5 3.5
\put{$\ssize \Cal U[-1]$} at 3.5 3.5
\endpicture} at 0 0 
\put{\beginpicture
\setcoordinatesystem units <0.3cm,0.3cm>
\put{} at 0 0
\put{} at 10 7 
\put{$\underbrace{\phantom{mmmmmmmmm}}_{\Cal D}$} at 5.25 -.8

\plot 0.5 2  6.5 2 /
\plot 0.5 5  1.5 5  2 5.5  5 5.5  5.5 5  6.5 5 /
\plot 7 7  7 0.5  8 0.5  8.5 0  9 0.5  10 0.5  10 7 /

\put{$\ssize \Cal T$} at 8.5 3.5
\put{$\ssize \Cal U$} at 3.5 3.5
\endpicture} at 5 0 
\put{\beginpicture
\setcoordinatesystem units <0.3cm,0.3cm>
\put{} at 0 0
\put{} at 10 7 
\put{$\underbrace{\phantom{mmmmmmmmm}}_{\Cal D[1]}$} at 5.25 -1

\plot 0.5 2  6.5 2 /
\plot 0.5 5  1.5 5  2 5.5  5 5.5  5.5 5  6.5 5 /
\plot 7 7  7 0.5  8 0.5  8.5 0  9 0.5  10 0.5  10 7 /

\put{$\ssize \Cal T[1]$} at 8.5 3.5
\put{$\ssize \Cal U[1]$} at 3.5 3.5
\endpicture} at 10 0 
\put{$\cdots$} at -4 0.75
\put{$\cdots$} at 14 0.75
\endpicture}
$$
Here, $\Cal T$ consists of a $\Bbb P_1(k)$-family of tubes $\Cal T$ of type $(5,3,2)$;
all the tubes in $\Cal T$ are stable with the exception of the  big tube 
which contains the projective injective module $Y[2]$. 
The connecting component $\Cal U$ has stable part of type 
$\Bbb Z\tilde{\Bbb E}_8$,
to which a non-stable orbit of length six is attached.

\medskip
It turns out that most indecomposables in $\Cal S_4(6)$
have a translate under the shift 
which is a module over the following algebra $B$ given by the quiver $Q$
and the relations as indicated.
$$
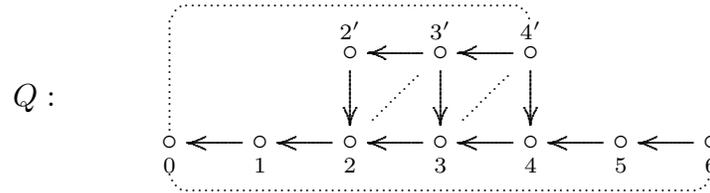
\hbox{\beginpicture
\setcoordinatesystem units <0.5cm,0.5cm>
\put{} at 0 -1.05
\put{} at 12 3.05
\put{$Q:$} at -3 1
\put{$\circ$} at 0 0
\put{$\circ$} at 2 0
\put{$\circ$} at 4 0
\put{$\circ$} at 6 0
\put{$\circ$} at 8 0
\put{$\circ$} at 10 0
\put{$\circ$} at 12 0

\put{$\circ$} at 4 2
\put{$\circ$} at 6 2
\put{$\circ$} at 8 2
\arr{1.6 0}{0.4 0}
\arr{3.6 0}{2.4 0}
\arr{5.6 0}{4.4 0}
\arr{7.6 0}{6.4 0}
\arr{9.6 0}{8.4 0}
\arr{11.6 0}{10.4 0}

\arr{7.6 2}{6.4 2}
\arr{5.6 2}{4.4 2}

\arr{8 1.6}{8 0.4}
\arr{4 1.6}{4 0.4}
\arr{6 1.6}{6 0.4}

\put{$\ssize 2'$} at  4 2.5
\put{$\ssize 3'$} at  6 2.5
\put{$\ssize 4'$} at  8 2.5

\put{$\ssize 0$} at 0 -.5
\put{$\ssize 1$} at 2 -.5
\put{$\ssize 2$} at 4 -.5
\put{$\ssize 3$} at 6 -.5
\put{$\ssize 4$} at 8 -.5
\put{$\ssize 5$} at 10 -.5
\put{$\ssize 6$} at 12 -.5

\setdots<2pt>
\plot 5.5 1.5  4.5 0.5 /
\plot 7.5 1.5  6.5 0.5 /
\plot         0 -0.7  0.1 -0.9  0.2 -1   0.3 -1.05 11.7 -1.05 11.8 -1  11.9 -0.9  12 -0.7 /
\plot 0 0.3   0 2.7   0.1  2.9  0.2 3    0.3 3.05   7.7  3.05  7.8  3   7.9  2.9   8  2.7 /
\endpicture}
$$
First we show that the Auslander-Reiten quiver of $B$ has the following structure.
$$
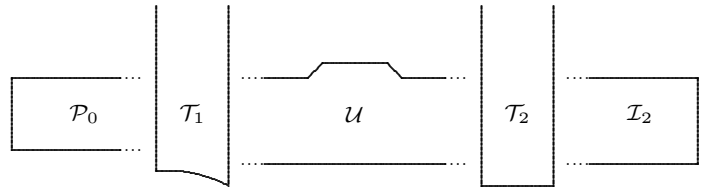
\hbox{\beginpicture
\setcoordinatesystem units <0.8cm,0.8cm>
\put{} at 0 0
\put{} at 9 2.5 
\plot 1.5 0.5  0 0.5  0 1.5  1.5 1.5 /
\setdots<2pt>
\plot 1.5 0.5  1.8 0.5 /
\plot 1.5 1.5  1.8 1.5 /
\plot 3.2 0.3  3.5 0.3 /
\plot 3.2 1.5  3.5 1.5 /
\plot 6 0.3  6.3 0.3 /
\plot 6 1.5  6.3 1.5 /
\plot 7.7 1.5  8 1.5 /
\plot 7.7 0.3  8 0.3 /
\setsolid
\plot 2 2.5  2 0.2 /
\plot 3 2.5  3 0 /
\plot 3.5 0.3  6 0.3 /
\plot 3.5 1.5  4.1 1.5  4.3 1.7  5.2 1.7  5.4 1.5  6 1.5 /
\plot 6.5 2.5  6.5 0  7.5 0  7.5 2.5 /
\plot 8 1.5  9.5 1.5  9.5 0.3  8 0.3 /
\setquadratic
\plot 2 0.2  2.5 0.18  3 0 /

\put{$\ssize \Cal P_0$} at 1 1
\put{$\ssize \Cal T_1$} at 2.5 1
\put{$\ssize \Cal U$} at 4.75 1
\put{$\ssize \Cal T_2$} at 7 1
\put{$\ssize \Cal I_2$} at 8.7 1
\endpicture}
$$
Let $Q_0$ be the full subquiver 
obtained from $Q$ by deleting the points $4'$ and $6$,
and let $B_0$ be the algebra given by $Q_0$ and the commutativity relation
as indicated. This algebra is tame concealed of tubular type $(4,3,2)$, 
and as $\rad P_6$ occurs on the mouth of the 
big tube, and $\rad P'_4$ is preinjective, it follows that 
the preprojective component $\Cal P_0$ of $B_0$ 
forms a preprojective component for $B$. We will picture this component later.
Here is the mouth of the 4-tube for $B_0$.
$$
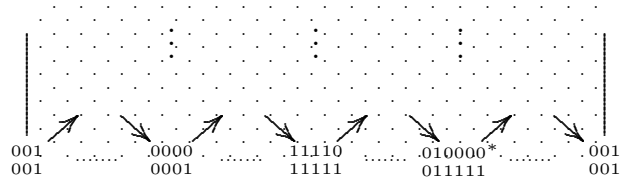
\hbox{\beginpicture
\setcoordinatesystem units <0.8cm,0.7cm>
\put{} at 0 1
\put{} at 8 2.5
\put{$\ssize {001     \atop 001    }$} at 0 1
\put{$\ssize {0000    \atop    0001}$} at 2 1
\put{$\ssize {11110   \atop   11111}$} at 4 1
\put{$\ssize {010000^*  \atop  011111^{\phantom *}}$} at 6 1
\put{$\ssize {001     \atop     001}$} at 8 1

\arr{0.3 1.3} {0.7 1.7}
\arr{1.3 1.7} {1.7 1.3}
\arr{2.3 1.3} {2.7 1.7}
\arr{3.3 1.7} {3.7 1.3}
\arr{4.3 1.3} {4.7 1.7}
\arr{5.3 1.7} {5.7 1.3}
\arr{6.3 1.3} {6.7 1.7}
\arr{7.3 1.7} {7.7 1.3}

\setdots<2pt>
\plot 0.7 1  1.3 1 /
\plot 2.7 1  3.3 1 /
\plot 4.7 1  5.3 1 /
\plot 6.7 1  7.3 1 /

\setsolid
\plot 0 1.4  0 3 /
\plot 8 1.4  8 3 /
\multiput{$\vdots$} at 2 3  4 3  6 3 /
\setshadegrid span <1.5mm>
\vshade  0  1 3.5 <,,,>
         8  1 3.5 /
\endpicture}
$$
Let $B_1$ be the algebra obtained from $B_0$ by a one-point 
extension at $\rad P_6$ which is the module labelled by a * in the 
above diagram. 
Thus, $B_1$ is the pathalgebra of the quiver $Q_1$ obtained from $Q$ by deleting the point
$4'$, with the induced relations.
This algebra $B_1$ is tame domestic, and according to [5, Theorem 4.9.2], its module
category consists of the following subcategories:
(a) the preprojective component $\Cal P_0$,
(b) a tubular family $\Cal T_1$ of type $(5,3,2)$, obtained from the tubular family for 
$B_0$ by a ray insertion in the 4-tube, and (c) a preinjective component $\Cal I_1$.
The tubular family separates $\Cal P_0$ from $\Cal I_1$.
 
The component $\Cal I_1$ contains the module $\rad P[4]$ as a $\tau$-successor of $I[0]$,
so when we insert a ray at $\rad P[4]$ we obtain a component of $B$-modules;
this is the component $\Cal U$ which we picture here.
$$
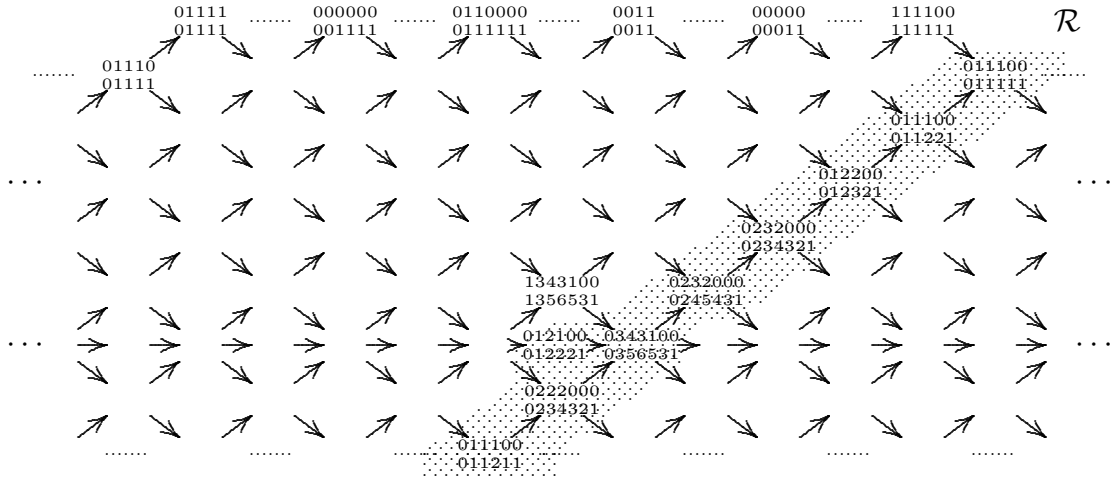
\hbox{\beginpicture
\setcoordinatesystem units <0.8cm,0.6cm>
\put{} at 0 0
\put{} at 12 7
\put{$\ssize {01110  \atop  01111}$}   at 1 7
\put{$\ssize {01111  \atop  01111}$}   at 2 8
\put{$\ssize {000000 \atop 001111}$}   at 4 8
\put{$\ssize {0110000\atop0111111}$}   at 6 8
\put{$\ssize {0011   \atop   0011}$}   at 8 8
\put{$\ssize {00000  \atop  00011}$}   at 10 8
\put{$\ssize {111100 \atop 111111}$}   at 12 8
\put{$\ssize {011100 \atop 011111}$}   at 13 7
\put{$\ssize {1343100\atop1356531}$}   at 7 3

\put{$\ssize {011100 \atop 011211}$}   at 6 0
\put{$\ssize {0222000\atop0234321}$}   at 7 1
\put{$\ssize {0343100\atop0356531}$}   at 8.1 2
\put{$\ssize {012100 \atop 012221}$}   at 6.9 2
\put{$\ssize {0232000\atop0245431}$}   at 9 3
\put{$\ssize {0232000\atop0234321}$}   at 10 4
\put{$\ssize {012200 \atop 012321}$}   at 11 5
\put{$\ssize {011100 \atop 011221}$}   at 12 6

\arr{0.3 0.3} {0.7 0.7}
\arr{2.3 0.3} {2.7 0.7}
\arr{4.3 0.3} {4.7 0.7}
\arr{6.3 0.3} {6.7 0.7}
\arr{8.3 0.3} {8.7 0.7}
\arr{10.3 0.3} {10.7 0.7}
\arr{12.3 0.3} {12.7 0.7}

\arr{1.3 0.7} {1.7 0.3}
\arr{3.3 0.7} {3.7 0.3}
\arr{5.3 0.7} {5.7 0.3}
\arr{7.3 0.7} {7.7 0.3}
\arr{9.3 0.7} {9.7 0.3}
\arr{11.3 0.7} {11.7 0.3}
\arr{13.3 0.7} {13.7 0.3}

\arr{0.3 1.7} {0.7 1.3} 
\arr{2.3 1.7} {2.7 1.3} 
\arr{4.3 1.7} {4.7 1.3} 
\arr{6.3 1.7} {6.7 1.3} 
\arr{8.3 1.7} {8.7 1.3} 
\arr{10.3 1.7} {10.7 1.3} 
\arr{12.3 1.7} {12.7 1.3} 

\arr{1.3 1.3} {1.7 1.7} 
\arr{3.3 1.3} {3.7 1.7} 
\arr{5.3 1.3} {5.7 1.7} 
\arr{7.3 1.3} {7.7 1.7} 
\arr{9.3 1.3} {9.7 1.7} 
\arr{11.3 1.3} {11.7 1.7} 
\arr{13.3 1.3} {13.7 1.7} 

\arr{0.3 2.0} {0.7 2.0}
\arr{1.3 2.0} {1.7 2.0}
\arr{2.3 2.0} {2.7 2.0}
\arr{3.3 2.0} {3.7 2.0}
\arr{4.3 2.0} {4.7 2.0}
\arr{5.3 2.0} {5.7 2.0}
\arr{6.3 2.0} {6.5 2.0}
\arr{7.4 2.0} {7.6 2.0}
\arr{8.6 2.0} {8.9 2.0}
\arr{9.3 2.0} {9.7 2.0}
\arr{10.3 2.0} {10.7 2.0}
\arr{11.3 2.0} {11.7 2.0}
\arr{12.3 2.0} {12.7 2.0}
\arr{13.3 2.0} {13.7 2.0}

\arr{0.3 2.3} {0.7 2.7} 
\arr{2.3 2.3} {2.7 2.7} 
\arr{4.3 2.3} {4.7 2.7} 
\arr{6.3 2.3} {6.7 2.7} 
\arr{8.3 2.3} {8.7 2.7} 
\arr{10.3 2.3} {10.7 2.7} 
\arr{12.3 2.3} {12.7 2.7} 

\arr{1.3 2.7} {1.7 2.3} 
\arr{3.3 2.7} {3.7 2.3} 
\arr{5.3 2.7} {5.7 2.3} 
\arr{7.3 2.7} {7.7 2.3} 
\arr{9.3 2.7} {9.7 2.3} 
\arr{11.3 2.7} {11.7 2.3} 
\arr{13.3 2.7} {13.7 2.3} 

\arr{0.3 3.7} {0.7 3.3} 
\arr{2.3 3.7} {2.7 3.3} 
\arr{4.3 3.7} {4.7 3.3} 
\arr{6.3 3.7} {6.7 3.3} 
\arr{8.3 3.7} {8.7 3.3} 
\arr{10.3 3.7} {10.7 3.3} 
\arr{12.3 3.7} {12.7 3.3} 

\arr{1.3 3.3} {1.7 3.7} 
\arr{3.3 3.3} {3.7 3.7} 
\arr{5.3 3.3} {5.7 3.7} 
\arr{7.3 3.3} {7.7 3.7} 
\arr{9.3 3.3} {9.7 3.7} 
\arr{11.3 3.3} {11.7 3.7} 
\arr{13.3 3.3} {13.7 3.7} 

\arr{0.3 4.3} {0.7 4.7} 
\arr{2.3 4.3} {2.7 4.7} 
\arr{4.3 4.3} {4.7 4.7} 
\arr{6.3 4.3} {6.7 4.7} 
\arr{8.3 4.3} {8.7 4.7} 
\arr{10.3 4.3} {10.7 4.7} 
\arr{12.3 4.3} {12.7 4.7} 

\arr{1.3 4.7} {1.7 4.3} 
\arr{3.3 4.7} {3.7 4.3} 
\arr{5.3 4.7} {5.7 4.3} 
\arr{7.3 4.7} {7.7 4.3} 
\arr{9.3 4.7} {9.7 4.3} 
\arr{11.3 4.7} {11.7 4.3} 
\arr{13.3 4.7} {13.7 4.3} 

\arr{0.3 5.7} {0.7 5.3} 
\arr{2.3 5.7} {2.7 5.3} 
\arr{4.3 5.7} {4.7 5.3} 
\arr{6.3 5.7} {6.7 5.3} 
\arr{8.3 5.7} {8.7 5.3} 
\arr{10.3 5.7} {10.7 5.3} 
\arr{12.3 5.7} {12.7 5.3} 

\arr{1.3 5.3} {1.7 5.7} 
\arr{3.3 5.3} {3.7 5.7} 
\arr{5.3 5.3} {5.7 5.7} 
\arr{7.3 5.3} {7.7 5.7} 
\arr{9.3 5.3} {9.7 5.7} 
\arr{11.3 5.3} {11.7 5.7} 
\arr{13.3 5.3} {13.7 5.7} 

\arr{0.3 6.3} {0.7 6.7} 
\arr{2.3 6.3} {2.7 6.7} 
\arr{4.3 6.3} {4.7 6.7} 
\arr{6.3 6.3} {6.7 6.7} 
\arr{8.3 6.3} {8.7 6.7} 
\arr{10.3 6.3} {10.7 6.7} 
\arr{12.3 6.3} {12.7 6.7} 

\arr{1.3 6.7} {1.7 6.3} 
\arr{3.3 6.7} {3.7 6.3} 
\arr{5.3 6.7} {5.7 6.3} 
\arr{7.3 6.7} {7.7 6.3} 
\arr{9.3 6.7} {9.7 6.3} 
\arr{11.3 6.7} {11.7 6.3} 
\arr{13.3 6.7} {13.7 6.3} 

\arr{2.3 7.7} {2.7 7.3} 
\arr{4.3 7.7} {4.7 7.3} 
\arr{6.3 7.7} {6.7 7.3} 
\arr{8.3 7.7} {8.7 7.3} 
\arr{10.3 7.7} {10.7 7.3} 
\arr{12.3 7.7} {12.7 7.3} 

\arr{1.3 7.3} {1.7 7.7} 
\arr{3.3 7.3} {3.7 7.7} 
\arr{5.3 7.3} {5.7 7.7} 
\arr{7.3 7.3} {7.7 7.7} 
\arr{9.3 7.3} {9.7 7.7} 
\arr{11.3 7.3} {11.7 7.7} 
\setdots<2pt>
\plot 0.7 0  1.3 0 /
\plot 2.7 0  3.3 0 /
\plot 4.7 0  5.3 0 /
\plot 6.7 0  7.3 0 /
\plot 8.7 0  9.3 0 /
\plot 10.7 0  11.3 0 /
\plot 12.7 0  13.3 0 /

\plot -0.3 7 0.3 7 /
\plot 2.7 8  3.3 8 /
\plot 4.7 8  5.3 8 /
\plot 6.7 8  7.3 8 /
\plot 8.7 8  9.3 8 /
\plot 10.7 8 11.2 8 /
\plot 13.7 7 14.3 7 /

\put{$\cdots$} at 14.4 2
\put{$\cdots$} at 14.4 5
\put{$\cdots$} at -0.4 2
\put{$\cdots$} at -0.4 5

\put{$\Cal R$} at 14 8

\setshadegrid span <0.6mm>
\hshade -0.5 5 7 <,,,z> 0 5 7 <,,z,z> 1.3 6.3 8.3 <,,z,z> 1.5 6.3 8.5 <,,z,z> 1.7 6.3 8.7 <,,z,z> 2 6.3 8.7 <,,z,z> 2.4 6.8 9 <,,z,z> 2.5 7.5 9.5 <,,,z> 7 12 14 <,,z,> 7.5 12.5 14 /

\endpicture} 
$$
The modules in the shaded region $\Cal R$ and on the right hand side of it
have no nonzero maps into $I[0]$ and hence are modules over the 
algebra $B_2$, which is given by the quiver $Q_2$
obtained from $Q$ by deleting the point $0$, and  by the two commutativity 
relations.  
It turns out that $\Cal R$ forms a slice in the preprojective component
for $B_2$; it follows that the right hand part of $\Cal U$ is just
the right hand part of this component.  

\smallskip
The algebra $B_2$ is tame concealed of tubular type
$(5,3,2)$; besides the modules in the preprojective component $\Cal P_2$, 
there is a family of tubes $\Cal T_2$ and a preinjective component
$\Cal I_2$. The tubular family separates $\Cal P_2$ from $\Cal I_2$,
and hence --- when considered a tubular family of $B$-modules ---
it separates $\Cal P_0$, $\Cal T_1$, and $\Cal U$ from $\Cal I_2$. 
This finishes our description of the Auslander-Reiten quiver for $B$.

\medskip
Note that none of the objects of type 
$Y[i]$ in $\Cal S_4(6)$ occurs as a $B$-module,
so we have to perform another --- final --- extension.
Consider the mouth of the big tube in $\Cal T_2$, 
which contains the module $\rad Y[2]$, labelled by a star.
$$
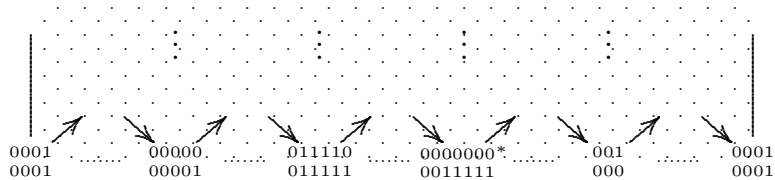
\hbox{\beginpicture
\setcoordinatesystem units <0.8cm,0.7cm>
\put{} at 0 1
\put{} at 10 2.5
\put{$\ssize {0001     \atop 0001    }$} at 0 1
\put{$\ssize {00000    \atop    00001}$} at 2 1
\put{$\ssize {011110   \atop   011111}$} at 4 1
\put{$\ssize {0000000^*  \atop  0011111^{\phantom *}}$} at 6 1
\put{$\ssize {001     \atop      000}$} at 8 1
\put{$\ssize {0001     \atop     0001}$} at 10 1

\arr{0.3 1.3} {0.7 1.7}
\arr{1.3 1.7} {1.7 1.3}
\arr{2.3 1.3} {2.7 1.7}
\arr{3.3 1.7} {3.7 1.3}
\arr{4.3 1.3} {4.7 1.7}
\arr{5.3 1.7} {5.7 1.3}
\arr{6.3 1.3} {6.7 1.7}
\arr{7.3 1.7} {7.7 1.3}
\arr{8.3 1.3} {8.7 1.7}
\arr{9.3 1.7} {9.7 1.3}

\setdots<2pt>
\plot 0.7 1  1.3 1 /
\plot 2.7 1  3.3 1 /
\plot 4.7 1  5.3 1 /
\plot 6.7 1  7.3 1 /
\plot 8.7 1  9.3 1 /

\setsolid
\plot 0 1.4  0 3 /
\plot 10 1.4  10 3 /

\multiput{$\vdots$} at 2 3  4 3  6 3  8 3 /
\setshadegrid span <1.5mm>
\vshade 0 1 3.5 <,,,> 10 1 3.5 /
\endpicture}
$$

Let $B'$ denote the algebra obtained from 
$B$ by a one-point extension at $\rad Y[2]$. 
The Auslander-Reiten quiver for $B'$ has the following structure.
$$
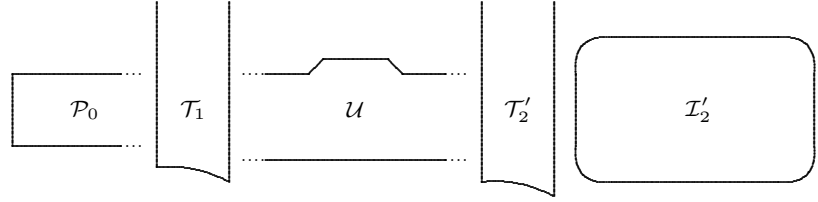
\hbox{\beginpicture
\setcoordinatesystem units <0.8cm,0.8cm>
\put{} at 0 0
\put{} at 9 2.5 

\setdots<2pt>
\plot 1.5 0.5  1.8 0.5 /
\plot 1.5 1.5  1.8 1.5 /
\plot 3.2 0.3  3.5 0.3 /
\plot 3.2 1.5  3.5 1.5 /
\plot 6 0.3  6.3 0.3 /
\plot 6 1.5  6.3 1.5 /

\setsolid
\plot 1.5 0.5  0 0.5  0 1.5  1.5 1.5 /
\plot 2 2.5  2 0.2 /
\plot 3 2.5  3 0 /
\plot 3.5 0.3  6 0.3 /
\plot 3.5 1.5  4.1 1.5  4.3 1.7  5.2 1.7  5.4 1.5  6 1.5 /
\plot 6.5 2.5  6.5 0  /
\plot 7.5 -0.2  7.5 2.5 /
\plot 10.6 2  8.3 2 /
\plot 10.6 0  8.3 0 /
\plot 11.1 0.5  11.1 1.5 /
\plot 7.8 0.5     7.8 1.5 /

\setquadratic
\plot 2 0.2  2.5 0.18  3 0 /
\plot 6.5 0  7 -0.02   7.5 -0.2 /
\plot 10.6 2 11 1.9 11.1 1.5 /
\plot 10.6 0 11 0.1 11.1 0.5 /
\plot 7.8 1.5  7.9 1.9  8.3 2 /
\plot 7.8 0.5  7.9 0.1  8.3 0 /

\put{$\ssize \Cal P_0$} at 1 1
\put{$\ssize \Cal T_1$} at 2.5 1
\put{$\ssize \Cal U$} at 4.75 1
\put{$\ssize \Cal T_2'$} at 7 1
\put{$\ssize \Cal I_2'$} at 9.5 1


\endpicture}
$$
Since $\Cal D'=\Cal P_0\vee\Cal T_1\vee\Cal U\vee\Cal T_2'$ 
contains the injective 
representation $I[0]$ 
and the relatively injective presentation $Y[2]$  but no unextended
radicals of projective representations, 
it follows that every indecomposable in $\Cal S_4(6)$
has a translate in $\Cal D'$ with respect to the shift.

\smallskip
We observe that the simple module $S'_2$ 
occurs at the mouth of the big tube in
$\Cal T_2'$ while the simple modules $S'_3$ and $S'_4$ 
lie in $\Cal I_2'$.
Thus, all the modules in $\Cal P_0\vee\Cal T_1\vee\Cal U\vee\Cal T_2'$
are representations of $\Cal S_4(6)$, 
with the exception of the modules in the ray 
starting at $S'_2$. 
Let $\Cal E$ be the tube obtained from the nonstable tube 
$\Cal E'$ in $\Cal T_2'$
by deleting the ray starting at $S'_2$.

\medskip\noindent{\bf Lemma 17.}
\item{1.} {\it The tube $\Cal E$ is a connected component
of the Auslander-Reiten quiver of $\Cal S_4(6)$.}
\item{2.} {\it Let $\Cal T$ be the family of tubes 
obtained from $\Cal T_2$ by replacing the nonstable
tube by $\Cal E$. Then $\Cal D=\Cal U\vee\Cal T$ is 
a fundamental domain for the shift in $\Cal S_4(6)$.}

\smallskip\noindent{\it Proof:\/}
1.\ After the ray deletion, the tube consists only of objects in $\Cal S_4(6)$.
We show that all source maps are source maps in $\Cal S_4(6)\cap B'$-mod. 
Consider the subgraph of $\Cal E'$ where $V_i,V_{i+1}$ are on the ray
starting at $S_2'$.
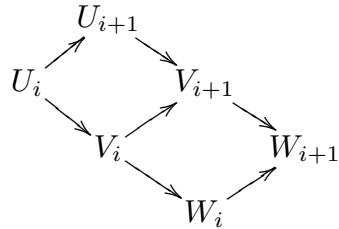
$$\xymatrix@1@=3mm{ 
        & U_{i+1}\ar[dr]\\ U_i\ar[ur]\ar[dr] & & V_{i+1} \ar[dr]\\
        & V_i\ar[ur]\ar[dr] & & W_{i+1}\\ & & W_i\ar[ur]}$$
Then the map $V_i\to W_i$ is a left approximation for $V_i$
in $\Cal S_4(6)$, and
the map $f\: U_i\to U_{i+1}\oplus W_i$, being the composition of the 
source map with the left approximation, is a source map in the 
category $B'\text{-mod}\cap S_4(6)$.
We will see later that --- as in the case of the category $\Cal S_3(7)$ ---
$f$ in fact is a source map in $S_4(6)$. 

Here is the mouth of the 5-tube in $\Cal T$, after this operation.
$$
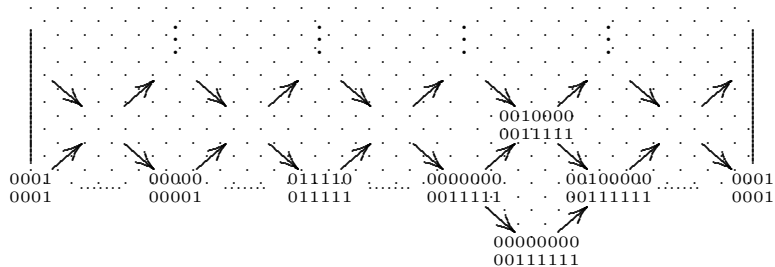
\hbox{\beginpicture
\setcoordinatesystem units <0.8cm,0.7cm>
\put{} at 0 0
\put{} at 10 3.5
\put{$\ssize {00000000 \atop 00111111}$} at 7 0

\put{$\ssize {0001     \atop 0001    }$} at 0 1
\put{$\ssize {00000    \atop    00001}$} at 2 1
\put{$\ssize {011110   \atop   011111}$} at 4 1
\put{$\ssize {0000000  \atop  0011111}$} at 6 1
\put{$\ssize {00100000 \atop 00111111}$} at 8 1
\put{$\ssize {0001     \atop     0001}$} at 10 1

\put{$\ssize {0010000  \atop  0011111}$} at 7 2

\arr{6.3 0.7} {6.7 0.3}
\arr{7.3 0.3} {7.7 0.7}

\arr{0.3 1.3} {0.7 1.7}
\arr{1.3 1.7} {1.7 1.3}
\arr{2.3 1.3} {2.7 1.7}
\arr{3.3 1.7} {3.7 1.3}
\arr{4.3 1.3} {4.7 1.7}
\arr{5.3 1.7} {5.7 1.3}
\arr{6.3 1.3} {6.7 1.7}
\arr{7.3 1.7} {7.7 1.3}
\arr{8.3 1.3} {8.7 1.7}
\arr{9.3 1.7} {9.7 1.3}

\arr{0.3 2.7} {0.7 2.3}
\arr{1.3 2.3} {1.7 2.7}
\arr{2.3 2.7} {2.7 2.3}
\arr{3.3 2.3} {3.7 2.7}
\arr{4.3 2.7} {4.7 2.3}
\arr{5.3 2.3} {5.7 2.7}
\arr{6.3 2.7} {6.7 2.3}
\arr{7.3 2.3} {7.7 2.7}
\arr{8.3 2.7} {8.7 2.3}
\arr{9.3 2.3} {9.7 2.7}
\setdots<2pt>
\plot 0.7 1  1.3 1 /
\plot 2.7 1  3.3 1 /
\plot 4.7 1  5.3 1 /

\plot 8.7 1  9.3 1 /

\setsolid
\plot 0 1.4  0 3.5 /
\plot 10 1.4  10 3.5 /
\multiput{$\vdots$} at 2 3.5  4 3.5  6 3.5  8 3.5 /
\setshadegrid span <1.5mm> 
\vshade 0 1 4 <z,,,>  6 1 4 <z,z,,>  7 0 4 <z,z,,>  8 1 4 <,z,,> 10 1 4 / 
\endpicture}
$$
2.\ We observe that the modules in $\Cal T_1$ are contained in $\Cal T[-1]$;
and the only representations of $\Cal T[-1]$ which are not in $\Cal T_1$ 
lie in the coray ending at $Y[1]$. To finish the proof one checks 
that each representation in $\Cal P_0$ has a translate in $\Cal D$. Here we 
picture $\Cal P_0$ 
and indicate for each representation to which category $\Cal T[i]$ 
or $\Cal U[i]$ it belongs. \qed
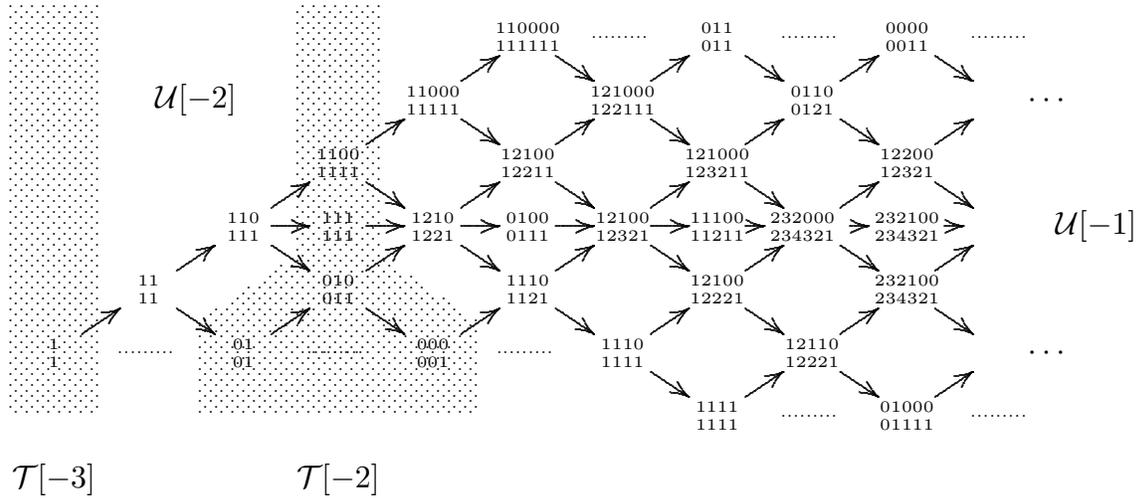
$$
\hbox{\beginpicture
\setcoordinatesystem units <1.05cm,0.7cm>
\put{} at 0 0
\put{} at 12 6

\put{$\ssize {1111    \atop    1111}$} at 7 -1
\put{$\ssize {01000   \atop   01111}$} at 9 -1

\put{$\ssize {1       \atop       1}$} at 0 0
\put{$\ssize {01 \atop           01}$} at 2 0
\put{$\ssize {000     \atop     001}$} at 4 0
\put{$\ssize {1110    \atop    1111}$} at 6 0
\put{$\ssize {12110   \atop   12221}$} at 8 0

\put{$\ssize {11      \atop      11}$} at 1 1
\put{$\ssize {010     \atop     011}$} at 3 1
\put{$\ssize {1110    \atop    1121}$} at 5 1
\put{$\ssize {12100   \atop   12221}$} at 7 1
\put{$\ssize {232100  \atop  234321}$} at 9 1

\put{$\ssize {110     \atop     111}$} at 2 2
\put{$\ssize {111     \atop     111}$} at 3 2
\put{$\ssize {1210    \atop    1221}$} at 4 2
\put{$\ssize {0100    \atop    0111}$} at 5 2
\put{$\ssize {12100   \atop   12321}$} at 6 2
\put{$\ssize {11100   \atop   11211}$} at 7 2
\put{$\ssize {232000  \atop  234321}$} at 7.9 2
\put{$\ssize {232100  \atop  234321}$} at 9 2

\put{$\ssize {1100    \atop    1111}$} at 3 3
\put{$\ssize {12100   \atop   12211}$} at 5 3
\put{$\ssize {121000  \atop  123211}$} at 7 3
\put{$\ssize {12200   \atop   12321}$} at 9 3

\put{$\ssize {11000   \atop   11111}$} at 4 4
\put{$\ssize {121000  \atop  122111}$} at 6 4
\put{$\ssize {0110    \atop    0121}$} at 8 4

\put{$\ssize {110000  \atop  111111}$} at 5 5
\put{$\ssize {011     \atop     011}$} at 7 5
\put{$\ssize {0000    \atop    0011}$} at 9 5

\arr{6.3 -0.3} {6.7 -0.7}
\arr{7.3 -0.7} {7.7 -0.3}
\arr{8.3 -0.3} {8.7 -0.7}
\arr{9.3 -0.7} {9.7 -0.3}

\arr{0.3 0.3} {0.7 0.7}
\arr{2.3 0.3} {2.7 0.7}
\arr{4.3 0.3} {4.7 0.7}
\arr{6.3 0.3} {6.7 0.7}
\arr{8.3 0.3} {8.7 0.7}

\arr{1.3 0.7} {1.7 0.3}
\arr{3.3 0.7} {3.7 0.3}
\arr{5.3 0.7} {5.7 0.3}
\arr{7.3 0.7} {7.7 0.3}
\arr{9.3 0.7} {9.7 0.3}

\arr{2.3 1.7} {2.7 1.3} 
\arr{4.3 1.7} {4.7 1.3} 
\arr{6.3 1.7} {6.7 1.3} 
\arr{8.3 1.7} {8.7 1.3} 

\arr{1.3 1.3} {1.7 1.7} 
\arr{3.3 1.3} {3.7 1.7} 
\arr{5.3 1.3} {5.7 1.7} 
\arr{7.3 1.3} {7.7 1.7} 
\arr{9.3 1.3} {9.7 1.7} 

\arr{2.3 2.0} {2.7 2.0}
\arr{3.3 2.0} {3.7 2.0}
\arr{4.3 2.0} {4.7 2.0}
\arr{5.3 2.0} {5.7 2.0}
\arr{6.3 2.0} {6.7 2.0}
\arr{7.3 2.0} {7.5 2.0}
\arr{8.4 2.0} {8.6 2.0}
\arr{9.4 2.0} {9.6 2.0}

\arr{2.3 2.3} {2.7 2.7} 
\arr{4.3 2.3} {4.7 2.7} 
\arr{6.3 2.3} {6.7 2.7} 
\arr{8.3 2.3} {8.7 2.7} 

\arr{3.3 2.7} {3.7 2.3} 
\arr{5.3 2.7} {5.7 2.3} 
\arr{7.3 2.7} {7.7 2.3} 
\arr{9.3 2.7} {9.7 2.3} 

\arr{4.3 3.7} {4.7 3.3} 
\arr{6.3 3.7} {6.7 3.3} 
\arr{8.3 3.7} {8.7 3.3} 

\arr{3.3 3.3} {3.7 3.7} 
\arr{5.3 3.3} {5.7 3.7} 
\arr{7.3 3.3} {7.7 3.7} 
\arr{9.3 3.3} {9.7 3.7} 

\arr{4.3 4.3} {4.7 4.7} 
\arr{6.3 4.3} {6.7 4.7} 
\arr{8.3 4.3} {8.7 4.7} 

\arr{5.3 4.7} {5.7 4.3} 
\arr{7.3 4.7} {7.7 4.3} 
\arr{9.3 4.7} {9.7 4.3} 


\setdots<2pt>
\plot 0.7 0  1.3 0 /
\plot 2.7 0  3.3 0 /
\plot 4.7 0  5.3 0 /
\plot 7.7 -1  8.3 -1 /
\plot 9.7 -1  10.3 -1 /

\plot 5.7 5  6.3 5 /
\plot 7.7 5  8.3 5 /
\plot 9.7 5 10.3 5 /

\put{$\cdots$} at 10.5 4
\put{$\cdots$} at 10.5 0
\setshadegrid span <0.6mm>
\hshade -1  1.5 4.5 <,,,z>  0.5 1.5 4.5 <,,z,z> 1.5 2.5 3.5 <,,z,> 5.5 2.5 3.5 /
\hshade -1  -0.5 0.5  5.5 -0.5 0.5 /

\put{$\Cal T[-3]$} at 0 -2
\put{$\Cal U[-2]$} at 1.5 4
\put{$\Cal U[-1]$} at 11 2
\put{$\Cal T[-2]$} at 3 -2
\endpicture}
$$

\smallskip
The Auslander-Reiten sequences in $\Cal U$ and in $\Cal T$ actually
are Auslander-Reiten sequences in the category $\Cal S_4(6)$:
Let $B^+$ be obtained from the algebra $B'$ by forming iterated
one point extensions at radicals of projective $\Cal S_4(6)$-modules.
Such radicals occur only in the category $\Cal I_2'$, so the components
$\Cal U$ and $\Cal T_2'$ are components of the Auslander-Reiten quiver
for $B^+$. The argument in the first statement in Lemma~17 yields that
$\Cal U$ and $\Cal T$ (which is $\Cal T_2'\cap\Cal S_4(6)$) consist
of Auslander-Reiten sequences in $\Cal S_4(6)\cap B^+$-mod. It follows,
as in the case $m=1$, that they are Auslander-Reiten sequences in $\Cal S_4(6)$.

\smallskip
The category $\Cal S_4(6)$ has separation properties similar to those for
$\Cal S_3(7)$.  The result corresponding to Proposition~13 can be shown
as in the section on $\Cal S_3(7)$, or by using the separation properties
in the big category $\Cal S(6)$ studied in [6].

\smallskip Using covering theory, we obtain the following conclusion for
$\Cal S_4(k[T]/T^6)$.

\medskip\noindent{\bf Proposition 18.} {\it
The category $\Cal S_4(k[T]/T^6)$ consists of:
\item{$\bullet$} A $\Bbb P_1(k)$ family of tubes $\Cal T$
of type $(5,3,2)$. All tubes are stable with the 
exception of the tube of circumference five 
pictured below which contains one nonstable module $Y$.
\item{$\bullet$} The connecting component $\Cal U$ 
which has stable orbit type $\tilde{\Bbb E}_8$ and a 
nonstable orbit of length 6 from $P$ to $I$.

\noindent Moreover, homomorphisms in the infinite radical of $\Cal T$
factor through any slice in $\Cal U$, and homomorphisms in the infinite
radical of $\Cal U$ factor through any of the tubes in $\Cal T$. } \qed

\bigskip
Let us now consider $\Cal S_4(k[T]/T^6)$ as a full subcategory of
$\Cal S(k[T]/T^6)$, this is a tubular category of tubular type $(6,3,2)$,
so each indecomposable occurs in a tube and the set of tubes is 
indexed by a rational angle $\gamma\in\Bbb Q/\Bbb Z$ and an irreducible
polynomial in $\Bbb P_1(k)$ (see [6]). 
All tubes are stable with the exception of the tube $\Cal T^6_0$ 
of circumference six in the family of index $\gamma=0+\Bbb Z$. 
We will need the modules $X$ and $Z$ below which occur as first term
and as endterm of the Auslander-Reiten sequence which contains $P=I$
as summand of the middle term and which forms part of the
mouth of $\Cal T^6_0$. 
$$0\;\longrightarrow\; \beginpicture\setcoordinatesystem units <10mm,10mm>
        \multiput{$\smallsq2$} at 0 -.5 *5 0 .2 /
        \multiput{$\ssize\bullet$} at .1 .3 / 
        \endpicture 
    \;\longrightarrow\; \beginpicture\setcoordinatesystem units <10mm,10mm>
        \multiput{$\smallsq2$} at 0 -.3 *4 0 .2 /
        \multiput{$\ssize\bullet$} at .1 .3 / 
        \endpicture 
    \,\oplus\beginpicture\setcoordinatesystem units <10mm,10mm>
        \multiput{$\smallsq2$} at 0 -.5 *5 0 .2 /
        \multiput{$\ssize\bullet$} at .1 .5 / 
        \endpicture 
    \;\longrightarrow\;\beginpicture\setcoordinatesystem units <10mm,10mm>
        \multiput{$\smallsq2$} at 0 -.3 *4 0 .2 /
        \multiput{$\ssize\bullet$} at .1 .5 / 
        \endpicture 
    \;\longrightarrow\; 0$$
Each indecomposable object $(A\sub B)$ for which $A$ has a summand
of type $k[T]/T^5$ admits a monomorphism from $Z$ and has a map
into $X$ which is an epimorphism on the subspaces. Note that the subspaces
in $X$ and in $Z$ are shifted against each other by one box. 
It follows that $(A\sub B)$ occurs either on the ray starting at $Z$ or on the 
coray ending at $X$ or in a tubular family of index $\gamma\neq 0+\Bbb Z$. 
Conversely, the stable tubes of index $\gamma = 0+\Bbb Z$ form the stable
tubes in the family of tubes $\Cal T$ in $\Cal S_4(k[T]/T^6)$ above.
And the nonstable tube in $\Cal T$ is obtained from the nonstable tube
in $\Cal T^6_0$ by deleting $P=I$ and the modules on the ray starting at $Z$
and the modules on the coray ending at $X$, as follows.

\def\bigtubefoursix{%
\beginpicture\setcoordinatesystem units <1cm,.9cm>
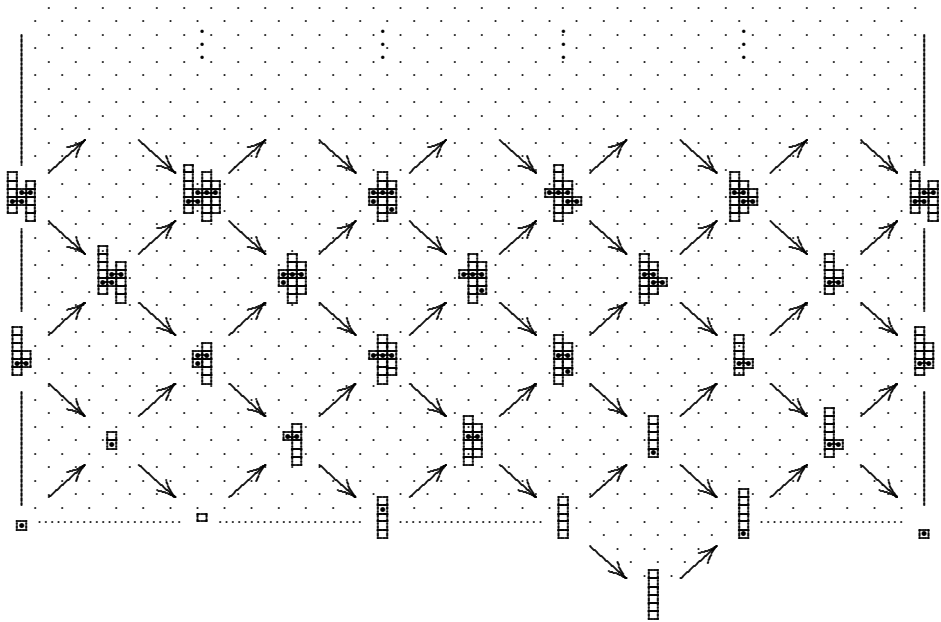
\put{The Tube of Circumference Five in $\Cal S_4(k[T]/T^6)$} at 5 -2
\multiput{$\smallsq1$} at 6.95 -1.15 *5 0 .1 /

\multiput{$\smallsq1$} at -.05 -.05 /
\multiput{$\boldkey.$} at 0 -.05 /

\multiput{$\smallsq1$} at 1.95 .05 /

\multiput{$\smallsq1$} at 3.95 -.15 *4 0 .1 /
\multiput{$\boldkey.$} at 4 .15 /

\multiput{$\smallsq1$} at 5.95 -.15 *4 0 .1 /

\multiput{$\smallsq1$} at 7.95 -.15 *5 0 .1  /
\multiput{$\boldkey.$} at 8 -.15 /

\multiput{$\smallsq1$} at 9.95 -.15 /
\multiput{$\boldkey.$} at 10 -.15 /

\multiput{$\smallsq1$} at  .95 .95 *1 0 .1 /
\multiput{$\boldkey.$} at  1 .95 /

\multiput{$\smallsq1$} at  2.9 1.05 /
\multiput{$\smallsq1$} at  3 .75 *4 0 .1 /
\multiput{$\boldkey.$} at  2.95 1.05  3.05 1.05 /
\plot  2.95 1.05  3.05 1.05 /

\multiput{$\smallsq1$} at  4.9 .75 *5 0 .1 /
\multiput{$\smallsq1$} at  5 .85 *3 0 .1 /
\multiput{$\boldkey.$} at 4.95 1.05 5.05 1.05 /
\plot  4.95 1.05 5.05 1.05 /

\multiput{$\smallsq1$} at  6.95 .85 *4 0 .1 /
\multiput{$\boldkey.$} at 7 .85 /

\multiput{$\smallsq1$} at 8.9 .85 *5  0 .1 /
\multiput{$\smallsq1$} at 9 .95 /
\multiput{$\boldkey.$} at 8.95 .95  9.05 .95 /
\plot  8.95 .95  9.05 .95 /

\multiput{$\smallsq1$} at  -.1 1.85 *5 0 .1 /
\multiput{$\smallsq1$} at  0 1.95 *1 0 .1 /
\multiput{$\boldkey.$} at -.05 1.95  .05 1.95 /
\plot -.05 1.95  .05 1.95 /

\multiput{$\smallsq1$} at  1.9 1.95 *1 0 .1 /
\multiput{$\smallsq1$} at 2 1.75 *4  0 .1 /
\multiput{$\boldkey.$} at 1.95 1.95  1.95 2.05  2.05 2.05 /
\plot   1.95 2.05  2.05 2.05 /

\multiput{$\smallsq1$} at 3.85 2.05 /
\multiput{$\smallsq1$} at 3.95 1.75 *5  0 .1 /
\multiput{$\smallsq1$} at 4.05 1.85 *3  0 .1 /
\multiput{$\boldkey.$} at 3.9 2.05  4 2.05  4.1 2.05 /
\plot 3.9 2.05 4.1 2.05 /

\multiput{$\smallsq1$} at 5.9 1.75 *5  0 .1 /
\multiput{$\smallsq1$} at 6 1.85 *3  0 .1 /
\multiput{$\boldkey.$} at 5.95 2.05  6.05 2.05  6.05 1.85 /
\plot  5.95 2.05  6.05 2.05 /

\multiput{$\smallsq1$} at  7.9 1.85 *4 0 .1 /
\multiput{$\smallsq1$} at  8 1.95 /
\multiput{$\boldkey.$} at  7.95 1.95  8.05 1.95 /
\plot  7.95 1.95  8.05 1.95 /

\multiput{$\smallsq1$} at  9.9 1.85 *5 0 .1 /
\multiput{$\smallsq1$} at  10 1.95 *2 0 .1 /
\multiput{$\boldkey.$} at  9.95 1.95  10.05 1.95 /
\plot   9.95 1.95  10.05 1.95 /

\multiput{$\smallsq1$} at .85 2.85 *5  0 .1 /
\multiput{$\smallsq1$} at .95 2.95 *1  0 .1 /
\multiput{$\smallsq1$} at 1.05 2.75 *4  0 .1 /
\multiput{$\boldkey.$} at .9 2.95  1 2.95  1 3.05  1.1 3.05 /
\plot  1 3.05  1.1 3.05 /
\plot  .9 2.95  1 2.95  /

\multiput{$\smallsq1$} at 2.85 2.95 *1  0 .1 /
\multiput{$\smallsq1$} at 2.95 2.75 *5  0 .1 /
\multiput{$\smallsq1$} at 3.05 2.85 *3  0 .1 /
\multiput{$\boldkey.$} at 2.9 2.95  2.9 3.05  3 3.05  3.1 3.05 /
\plot  2.9 3.05   3.1 3.05 /

\multiput{$\smallsq1$} at 4.85 3.05 /
\multiput{$\smallsq1$} at 4.95 2.75 *5  0 .1 /
\multiput{$\smallsq1$} at 5.05 2.85 *3  0 .1 /
\multiput{$\boldkey.$} at 4.9 3.05  5 3.05  5.1 3.05  5.1 2.85 /
\plot  4.9 3.05   5.1 3.05 /

\multiput{$\smallsq1$} at 6.85 2.75 *5  0 .1 /
\multiput{$\smallsq1$} at 6.95 2.85 *3  0 .1 /
\multiput{$\smallsq1$} at 7.05 2.95 /
\multiput{$\boldkey.$} at 6.9 3.05  7 3.05  7 2.95  7.1 2.95 /
\plot  6.9 3.05  7 3.05 /
\plot   7 2.95  7.1 2.95 /

\multiput{$\smallsq1$} at 8.9 2.85 *4  0 .1 /
\multiput{$\smallsq1$} at 9 2.95 *1  0 .1 /
\multiput{$\boldkey.$} at 8.95 2.95  9.05 2.95 /
\plot  8.95 2.95  9.05 2.95 /

\multiput{$\smallsq1$} at -.15 3.85 *4  0 .1 /
\multiput{$\smallsq1$} at -.05 3.95 *1  0 .1 /
\multiput{$\smallsq1$} at  .05 3.75 *4 0 .1 /
\multiput{$\boldkey.$} at -.1 3.95  0 3.95  0 4.05  .1 4.05 /
\plot  -.1 3.95  0 3.95 /
\plot   0 4.05  .1 4.05 /

\multiput{$\smallsq1$} at 1.8 3.85 *5  0 .1 /
\multiput{$\smallsq1$} at 1.9 3.95 *1  0 .1 /
\multiput{$\smallsq1$} at 2 3.75 *5  0 .1 /
\multiput{$\smallsq1$} at 2.1 3.85 *3  0 .1 /
\multiput{$\boldkey.$} at 1.85 3.95  1.95 3.95  1.95 4.05  2.05 4.05  2.15 4.05 /
\plot  1.85 3.95  1.95 3.95 /
\plot   1.95 4.05   2.15 4.05 /

\multiput{$\smallsq1$} at 3.85 3.95 *1  0 .1 /
\multiput{$\smallsq1$} at 3.95 3.75 *5  0 .1 /
\multiput{$\smallsq1$} at 4.05 3.85 *3  0 .1 /
\multiput{$\boldkey.$} at 3.9 3.95  3.9 4.05  4 4.05  4.1 4.05  4.1 3.85 /
\plot   3.9 4.05   4.1 4.05 /

\multiput{$\smallsq1$} at  5.8 4.05 /
\multiput{$\smallsq1$} at  5.9 3.75 *6 0 .1 /
\multiput{$\smallsq1$} at  6 3.85 *3  0 .1 /
\multiput{$\smallsq1$} at  6.1 3.95 /
\multiput{$\boldkey.$} at  5.85 4.05  5.95 4.05  6.05 4.05  6.05 3.95  6.15 3.95 /
\plot  6.05 3.95  6.15 3.95 /
\plot  5.85 4.05  6.05 4.05 /

\multiput{$\smallsq1$} at  7.85 3.75 *5 0 .1 /
\multiput{$\smallsq1$} at  7.95 3.85 *3 0 .1 /
\multiput{$\smallsq1$} at  8.05 3.95 *1 0 .1 /
\multiput{$\boldkey.$} at  7.9 4.05  8 4.05  8 3.95  8.1 3.95 /
\plot  7.9 4.05  8 4.05 /
\plot    8 3.95  8.1 3.95 /

\multiput{$\smallsq1$} at 9.85 3.85 *4  0 .1 /
\multiput{$\smallsq1$} at 9.95 3.95 *1  0 .1 /
\multiput{$\smallsq1$} at 10.05 3.75 *4  0 .1 /
\multiput{$\boldkey.$} at 9.9 3.95  10 3.95  10 4.05  10.1 4.05 /
\plot  10 4.05  10.1 4.05 /
\plot  9.9 3.95  10 3.95 /

\arr{7.3 -.7}{7.7 -.3}
\arr{6.3 -.3}{6.7 -.7}
\arr{0.3 0.3}{0.7 0.7}
\arr{2.3 0.3}{2.7 0.7}
\arr{4.3 0.3}{4.7 0.7}
\arr{6.3 0.3}{6.7 0.7}
\arr{8.3 0.3}{8.7 0.7}
\arr{1.3 0.7}{1.7 0.3}
\arr{3.3 0.7}{3.7 0.3}
\arr{5.3 0.7}{5.7 0.3}
\arr{7.3 0.7}{7.7 0.3}
\arr{9.3 0.7}{9.7 0.3}
\arr{0.3 1.7}{0.7 1.3}
\arr{2.3 1.7}{2.7 1.3}
\arr{4.3 1.7}{4.7 1.3}
\arr{6.3 1.7}{6.7 1.3}
\arr{8.3 1.7}{8.7 1.3}
\arr{1.3 1.3}{1.7 1.7}
\arr{3.3 1.3}{3.7 1.7}
\arr{5.3 1.3}{5.7 1.7}
\arr{7.3 1.3}{7.7 1.7}
\arr{9.3 1.3}{9.7 1.7}

\arr{0.3 2.3}{0.7 2.7}
\arr{2.3 2.3}{2.7 2.7}
\arr{4.3 2.3}{4.7 2.7}
\arr{6.3 2.3}{6.7 2.7}
\arr{8.3 2.3}{8.7 2.7}
\arr{1.3 2.7}{1.7 2.3}
\arr{3.3 2.7}{3.7 2.3}
\arr{5.3 2.7}{5.7 2.3}
\arr{7.3 2.7}{7.7 2.3}
\arr{9.3 2.7}{9.7 2.3}
\arr{0.3 3.7}{0.7 3.3}
\arr{2.3 3.7}{2.7 3.3}
\arr{4.3 3.7}{4.7 3.3}
\arr{6.3 3.7}{6.7 3.3}
\arr{8.3 3.7}{8.7 3.3}
\arr{1.3 3.3}{1.7 3.7}
\arr{3.3 3.3}{3.7 3.7}
\arr{5.3 3.3}{5.7 3.7}
\arr{7.3 3.3}{7.7 3.7}
\arr{9.3 3.3}{9.7 3.7}

\arr{0.3 4.3}{0.7 4.7}
\arr{2.3 4.3}{2.7 4.7}
\arr{4.3 4.3}{4.7 4.7}
\arr{6.3 4.3}{6.7 4.7}
\arr{8.3 4.3}{8.7 4.7}
\arr{1.3 4.7}{1.7 4.3}
\arr{3.3 4.7}{3.7 4.3}
\arr{5.3 4.7}{5.7 4.3}
\arr{7.3 4.7}{7.7 4.3}
\arr{9.3 4.7}{9.7 4.3}

\multiput{$\vdots$} at 2 6 *3 2 0 /

\setdots <2pt>
\plot .2 0  1.8 0 /
\plot 2.2 0  3.8 0 /
\plot 4.2 0  5.8 0 /
\plot 8.2 0  9.8 0 /

\setsolid
\plot 0 .2  0 1.6 /
\plot 0 2.6  0 3.6 /
\plot 0 4.4  0 6 /
\plot 10 .2  10 1.6 /
\plot 10 2.6 10 3.6 /
\plot 10 4.4 10 6 /

\setshadegrid span <1.5mm>
\vshade    0  0 6.5  <,z,,> 
           6  0 6.5  <z,z,,>
           7 -1 6.5  <z,z,,>
           8  0 6.5  <z,,,>
          10  0 6.5 /
\endpicture}

$$\bigtubefoursix$$

The inclusion map $Z\to X$ factors over some module in each of 
the tubes with index $\gamma\neq 0+\Bbb Z$. Thus, each tube with
$\gamma\neq 0+\Bbb Z$ has a module on its mouth which has $k[T]/T^5$
as summand of the subspace --- with one possible exception:
The Auslander-Reiten sequence in $\Cal S(k[T]/T^6)$ 
(involving the modules $I$ and $P$ from $\Cal S_4(k[T]/T^6)$)
$$0\;\longrightarrow\; \beginpicture\setcoordinatesystem units <10mm,10mm>
        \multiput{$\smallsq2$} at 0 -.5 *5 0 .2 /
        \multiput{$\ssize\bullet$} at .1 .1 / 
        \endpicture 
    \;\longrightarrow\; \beginpicture\setcoordinatesystem units <10mm,10mm>
        \multiput{$\smallsq2$} at 0 -.5 *5 0 .2 /
        \multiput{$\smallsq2$} at 0.2 -.3 *3 0 .2 /
        \multiput{$\ssize\bullet$} at .1 .3  .3 .3  .3 .1 /
        \plot .1 .3  .3 .3 / 
        \endpicture 
    \;\longrightarrow\;\beginpicture\setcoordinatesystem units <10mm,10mm>
        \multiput{$\smallsq2$} at 0 -.3 *3 0 .2 /
        \multiput{$\ssize\bullet$} at .1 .3 / 
        \endpicture 
    \;\longrightarrow\; 0$$
has the property that the sequence of subspaces is not split exact, 
and it is the only Auslander-Reiten 
sequence with this property in which a summand 
of type $k[T]/T^4$ in the subspace is involved.  
This sequence occurs on the mouth of the tube $\Cal T_{1/2}^6$
of circumference six and index $\gamma=\frac12+\Bbb Z$, and hence 
$\Cal T_{1/2}^6\cap S_4(k[T]/T^6)$ consists of all
modules on the wing bounded by the mouth of the tube, the ray 
starting at the endterm and the coray ending at the first term of the 
above sequence. There are three indecomposables in 
$\Cal T_{1/2}^3\cap S_4(k[T]/T^6)$ which form an Auslander-Reiten sequence along the
mouth and a single module in $\Cal T_{1/2}^2\cap S_4(k[T]/T^6)$.
This wing triple of size $(6,2,1)$ forms the diamond shaped center piece of the
connecting component $\Cal U$, as pictured below.  

\def\diamondfoursix{%
\beginpicture\setcoordinatesystem units <.43cm,.4cm>
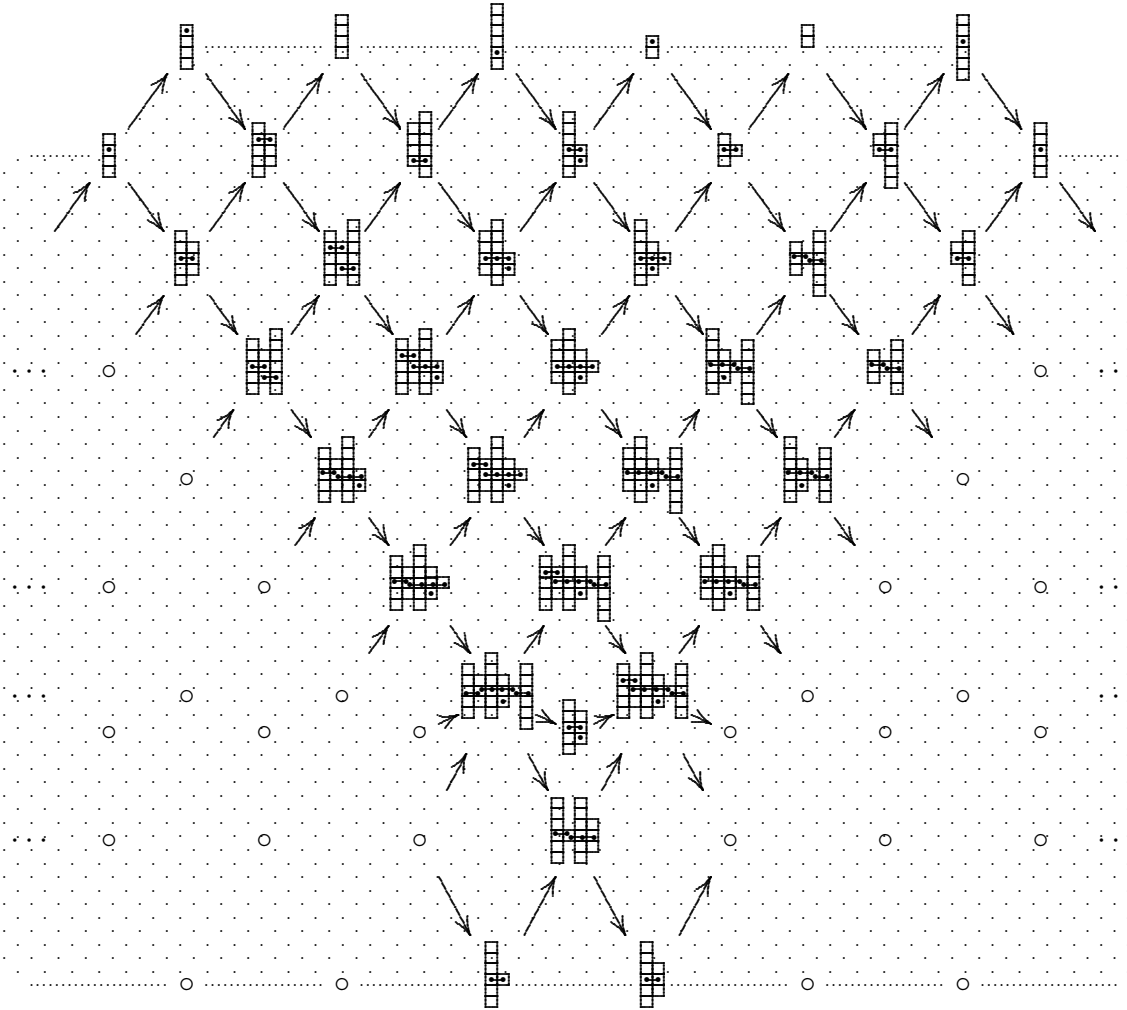
\put{The Diamond in the Connecting Component in  $\Cal S_4(k[T]/T^6)$} at 14 -3
\multiput{$\smallsq3$} at 3.85 25.55  *3 0 .3 /
\multiput{$\boldkey.$} at 4 26.45 /

\multiput{$\smallsq3$} at 7.85 25.85  *3 0 .3 /

\multiput{$\smallsq3$} at 11.85 25.55 *5 0 .3 /
\multiput{$\boldkey.$} at 12 25.85 /

\multiput{$\smallsq3$} at 15.85 25.85 *1 0 .3 /
\multiput{$\boldkey.$} at 16 26.15 /

\multiput{$\smallsq3$} at 19.85 26.15 *1 0 .3 /

\multiput{$\smallsq3$} at 23.85 25.25 *5 0 .3 /
\multiput{$\boldkey.$} at 24 26.15 /

\multiput{$\smallsq3$} at 1.85 22.55  *3 0 .3 /
\multiput{$\boldkey.$} at 2 23.15 /

\multiput{$\smallsq3$} at 5.7 22.55 *4 0 .3 /
\multiput{$\smallsq3$} at 6 22.85 *2 0 .3 /
\multiput{$\boldkey.$} at 5.85 23.45  6.15 23.45 /
\plot 5.85 23.45  6.15 23.45 /

\multiput{$\smallsq3$} at 9.7 22.85 *3 0 .3 /
\multiput{$\smallsq3$} at 10 22.55 *5  0 .3 /
\multiput{$\boldkey.$} at 9.85 22.85  10.15 22.85 /
\plot 9.85 22.85  10.15 22.85 /

\multiput{$\smallsq3$} at 13.7 22.55  *5 0 .3 /
\multiput{$\smallsq3$} at 14 22.85  *1 0 .3 /
\multiput{$\boldkey.$} at 13.85 23.15  14.15 23.15  14.15 22.85 /
\plot 13.85 23.15  14.15 23.15 /

\multiput{$\smallsq3$} at 17.7 22.85 *2 0 .3 /
\multiput{$\smallsq3$} at 18 23.15 /
\multiput{$\boldkey.$} at 17.85 23.15  18.15 23.15 /
\plot 17.85 23.15  18.15 23.15 /

\multiput{$\smallsq3$} at 21.7 23.15 *1 0 .3 /
\multiput{$\smallsq3$} at 22 22.25 *5 0 .3 /
\multiput{$\boldkey.$} at 21.85 23.15  22.15 23.15 /
\plot  21.85 23.15  22.15 23.15 /

\multiput{$\smallsq3$} at 25.85 22.55 *4 0 .3 /
\multiput{$\boldkey.$} at 26 23.15 /

\multiput{$\smallsq3$} at 3.7 19.55 *4 0 .3 /
\multiput{$\smallsq3$} at 4 19.85 *2 0 .3 /
\multiput{$\boldkey.$} at 3.85 20.15  4.15 20.15 /
\plot 3.85 20.15  4.15 20.15 /

\multiput{$\smallsq3$} at 7.55 19.55 *4 0 .3 /
\multiput{$\smallsq3$} at 7.85 19.85 *2 0 .3 /
\multiput{$\smallsq3$} at 8.15 19.55 *5 0 .3 /
\multiput{$\boldkey.$} at 7.7 20.45  8 20.45  8 19.85  8.3 19.85  /
\plot  7.7 20.45  8 20.45 /
\plot  8 19.85  8.3 19.85  /

\multiput{$\smallsq3$} at 11.55 19.85 *3 0 .3 /
\multiput{$\smallsq3$} at 11.85 19.55 *5 0 .3 /
\multiput{$\smallsq3$} at 12.15 19.85 *1 0 .3 /
\multiput{$\boldkey.$} at 11.7 20.15  12.3 20.15  12 20.15  12.3 19.85 /
\plot  11.7 20.15  12.3 20.15 /

\multiput{$\smallsq3$} at 15.55 19.55 *5 0 .3 /
\multiput{$\smallsq3$} at 15.85 19.85 *2 0 .3 /
\multiput{$\smallsq3$} at 16.15 20.15 /
\multiput{$\boldkey.$} at 15.7 20.15  16.3 20.15  16 20.15  16 19.85 /
\plot 15.7 20.15  16.3 20.15 /

\multiput{$\smallsq3$} at 19.55 19.85  *2 0 .3 /
\multiput{$\smallsq3$} at 19.85 20.15 /
\multiput{$\smallsq3$} at 20.15 19.25 *5 0 .3 /
\multiput{$\boldkey.$} at 19.65 20.2  19.95 20.2  20.05 20.1  20.35 20.1 /
\plot 20.05 20.1  20.35 20.1 /
\plot   19.65 20.2  19.95 20.2 /

\multiput{$\smallsq3$} at 23.7 20.15  *1 0 .3 /
\multiput{$\smallsq3$} at 24 19.55 *4 0 .3 /
\multiput{$\boldkey.$} at 23.85 20.15  24.15 20.15 /
\plot 23.85 20.15  24.15 20.15 /

\multiput{$\smallsq3$} at 5.55 16.55 *4 0 .3 /
\multiput{$\smallsq3$} at 5.85 16.85 *2 0 .3 /
\multiput{$\smallsq3$} at 6.15 16.55 *5 0 .3 /
\multiput{$\boldkey.$} at 5.7 17.15  6 17.15  6 16.85  6.3 16.85 /
\plot 6 16.85  6.3 16.85 /
\plot 5.7 17.15  6 17.15 /

\multiput{$\smallsq3$} at 9.4 16.55 *4 0 .3 /
\multiput{$\smallsq3$} at 9.7 16.85 *2 0 .3 /
\multiput{$\smallsq3$} at 10 16.55  *5 0 .3 /
\multiput{$\smallsq3$} at 10.3 16.85 *1 0 .3 /
\multiput{$\boldkey.$} at 9.55 17.45  9.85 17.45  9.85 17.15  10.45 17.15 
                                10.15 17.15  10.45 16.85 /
\plot 9.55 17.45  9.85 17.45 /
\plot 9.85 17.15  10.45 17.15 /

\multiput{$\smallsq3$} at 13.4 16.85 *3 0 .3 /
\multiput{$\smallsq3$} at 13.7 16.55 *5 0 .3 /
\multiput{$\smallsq3$} at 14   16.85 *2 0 .3 /
\multiput{$\smallsq3$} at 14.3 17.15 /
\multiput{$\boldkey.$} at 13.55 17.15  14.45 17.15  13.85 17.15  14.15 17.15
                                14.15 16.85 /
\plot 13.55 17.15  14.45 17.15 /

\multiput{$\smallsq3$} at 17.4 16.55 *5 0 .3 /
\multiput{$\smallsq3$} at 17.7 16.85 *2 0 .3 /
\multiput{$\smallsq3$} at 18 17.15 /
\multiput{$\smallsq3$} at 18.3 16.25 *5 0 .3 /
\multiput{$\boldkey.$} at 17.5 17.2  17.8 17.2  18.1 17.2  18.2 17.1  18.5 17.1
                                17.85 16.85 /
\plot 17.5 17.2  18.1 17.2 /
\plot 18.2 17.1  18.5 17.1 /

\multiput{$\smallsq3$} at 21.55 16.85 *2 0 .3 /
\multiput{$\smallsq3$} at 21.85 17.15 /
\multiput{$\smallsq3$} at 22.15 16.55 *4 0 .3 /
\multiput{$\boldkey.$} at 21.65 17.2  21.95 17.2  22.05 17.1  22.35 17.1 /
\plot 21.65 17.2  21.95 17.2 /
\plot 22.05 17.1  22.35 17.1 /

\multiput{$\smallsq3$} at 7.4 13.55 *4 0 .3 /
\multiput{$\smallsq3$} at 7.7 13.85 *2 0 .3 /
\multiput{$\smallsq3$} at 8   13.55 *5 0 .3 /
\multiput{$\smallsq3$} at 8.3 13.85 *1 0 .3 /
\multiput{$\boldkey.$} at 7.5 14.2  7.8 14.2  7.9 14.1  8.2 14.1  8.5 14.1
                                        8.45 13.85  /
\plot 7.5 14.2  7.8 14.2 /
\plot 7.9 14.1  8.5 14.1 /

\multiput{$\smallsq3$} at 11.25 13.55 *4 0 .3 /
\multiput{$\smallsq3$} at 11.55 13.85 *2 0 .3 /
\multiput{$\smallsq3$} at 11.85 13.55 *5 0 .3 /
\multiput{$\smallsq3$} at 12.15 13.85 *2 0 .3 /
\multiput{$\smallsq3$} at 12.45 14.15 /
\multiput{$\boldkey.$} at 11.4 14.45  11.7 14.45  11.7 14.15  12.6 14.15
                                12 14.15  12.3 14.15  12.3 13.85   /
\plot 11.4 14.45  11.7 14.45 /
\plot  11.7 14.15  12.6 14.15 /

\multiput{$\smallsq3$} at 15.25 13.85 *3 0 .3 /
\multiput{$\smallsq3$} at 15.55 13.55 *5 0 .3 /
\multiput{$\smallsq3$} at 15.85 13.85 *2 0 .3 /
\multiput{$\smallsq3$} at 16.15 14.15 /
\multiput{$\smallsq3$} at 16.45 13.25 *5 0 .3 /
\multiput{$\boldkey.$} at 15.35 14.2  15.65 14.2  15.95 14.2  16.25 14.2
                                16.35 14.1  16.65 14.1  16 13.85   /
\plot 15.35 14.2  16.25 14.2 /
\plot 16.35 14.1  16.65 14.1 /

\multiput{$\smallsq3$} at 19.4 13.55 *5 0 .3 /
\multiput{$\smallsq3$} at 19.7 13.85 *2 0 .3 /
\multiput{$\smallsq3$} at 20   14.15 /
\multiput{$\smallsq3$} at 20.3 13.55 *4 0 .3 /
\multiput{$\boldkey.$} at 19.5 14.2  20.1 14.2  19.8 14.2  20.2 14.1  20.5 14.1
                                19.85 13.85  /
\plot 19.5 14.2  20.1 14.2 /
\plot  20.2 14.1  20.5 14.1 /


\multiput{$\smallsq3$} at 9.25 10.55 *4 0 .3 /
\multiput{$\smallsq3$} at 9.55 10.85 *2 0 .3 /
\multiput{$\smallsq3$} at 9.85 10.55 *5 0 .3 /
\multiput{$\smallsq3$} at 10.15 10.85 *2 0 .3 /
\multiput{$\smallsq3$} at 10.45 11.15 /
\multiput{$\boldkey.$} at 9.35 11.2  9.65 11.2  9.75 11.1  10.65 11.1  
                                10.05 11.1  10.35 11.1  10.3 10.85  /
\plot 9.35 11.2  9.65 11.2 /
\plot 9.75 11.1  10.65 11.1 /

\multiput{$\smallsq3$} at 13.1 10.55 *4 0 .3 /
\multiput{$\smallsq3$} at 13.4 10.85 *2 0 .3 /
\multiput{$\smallsq3$} at 13.7 10.55 *5 0 .3 /
\multiput{$\smallsq3$} at 14   10.85 *2 0 .3 /
\multiput{$\smallsq3$} at 14.3 11.15 /
\multiput{$\smallsq3$} at 14.6 10.25 *5 0 .3 /
\multiput{$\boldkey.$} at 13.25 11.45  13.55 11.45  13.5 11.2  
                           14.4 11.2  13.8 11.2
                           14.1 11.2  14.5 11.1  14.8 11.1  14.15 10.85  /
\plot 13.25 11.45  13.55 11.45 /
\plot 13.5 11.2  14.4 11.2  /
\plot  14.5 11.1  14.8 11.1 /

\multiput{$\smallsq3$} at 17.25  10.85 *3  0 .3 /
\multiput{$\smallsq3$} at 17.55  10.55 *5 0 .3 /
\multiput{$\smallsq3$} at 17.85  10.85 *2 0 .3 /
\multiput{$\smallsq3$} at 18.15  11.15 /
\multiput{$\smallsq3$} at 18.45  10.55 *4 0 .3 /
\multiput{$\boldkey.$} at 17.35 11.2  18.25 11.2  17.65 11.2  17.95 11.2
                          18.35 11.1  18.65 11.1  18 10.85  /
\plot 17.35 11.2  18.25 11.2 /
\plot 18.35 11.1  18.65 11.1 /

\multiput{$\smallsq3$} at 11.1 7.55 *4 0 .3 /
\multiput{$\smallsq3$} at 11.4 7.85 *2 0 .3 /
\multiput{$\smallsq3$} at 11.7 7.55 *5 0 .3 /
\multiput{$\smallsq3$} at 12   7.85 *2 0 .3 /
\multiput{$\smallsq3$} at 12.3 8.15 /
\multiput{$\smallsq3$} at 12.6 7.25 *5 0 .3 /
\multiput{$\boldkey.$} at 11.2 8.1  11.5 8.1  11.6 8.2  12.4 8.2  11.87 8.2
                                    12.13 8.2  12.5 8.1  12.8 8.1  12.15 7.85 /
\plot 11.2 8.1  11.5 8.1 /
\plot 11.6 8.2  12.4 8.2 /
\plot  12.5 8.1  12.8 8.1 /

\multiput{$\smallsq3$} at 15.1 7.55 *4 0 .3 /
\multiput{$\smallsq3$} at 15.4 7.85 *2 0 .3 /
\multiput{$\smallsq3$} at 15.7 7.55 *5 0 .3 /
\multiput{$\smallsq3$} at 16   7.85 *2 0 .3 /
\multiput{$\smallsq3$} at 16.3 8.15 /
\multiput{$\smallsq3$} at 16.6 7.55 *4 0 .3 /
\multiput{$\boldkey.$} at 15.25 8.45  15.55 8.45  15.5 8.2  16.4 8.2  15.8 8.2
                                    16.1 8.2  16.5 8.1  16.8 8.1  16.15 7.85 /
\plot 15.25 8.45  15.55 8.45 /
\plot 15.5 8.2  16.4 8.2 /
\plot  16.5 8.1  16.8 8.1 /

\multiput{$\smallsq3$} at 13.7 6.55 *4 0 .3 /
\multiput{$\smallsq3$} at 14 6.85 *2 0 .3 /
\multiput{$\boldkey.$} at 13.85 7.15  14.15 7.15  14.15 6.85  /
\plot 13.85 7.15  14.15 7.15 /

\multiput{$\smallsq3$} at 13.4 3.55 *5 0 .3 /
\multiput{$\smallsq3$} at 13.7 4.15 /
\multiput{$\smallsq3$} at 14   3.55 *5 0 .3 /
\multiput{$\smallsq3$} at 14.3 3.85 *2 0 .3 /
\multiput{$\boldkey.$} at 13.5 4.2  13.8 4.2  13.9 4.1  14.5 4.1  14.2 4.1  /
\plot 13.5 4.2  13.8 4.2 / 
\plot 13.9 4.1  14.5 4.1 /

\multiput{$\smallsq3$} at 11.7 -.45 *5 0 .3 /
\multiput{$\smallsq3$} at 12 .15 /
\multiput{$\boldkey.$} at 11.85 .15  12.15 .15  /
\plot 11.85 .15  12.15 .15  /

\multiput{$\smallsq3$} at 15.7 -.45 *5 0 .3 /
\multiput{$\smallsq3$} at 16 -.15 *2 0 .3 /
\multiput{$\boldkey.$} at 15.85 .15  16.15 .15  /
\plot 15.85 .15  16.15 .15  /

\arr{2.5 23.75} {3.5 25.25}
\arr{6.5 23.75} {7.5 25.25}
\arr{10.5 23.75} {11.5 25.25}
\arr{14.5 23.75} {15.5 25.25}
\arr{18.5 23.75} {19.5 25.25}
\arr{22.5 23.75} {23.5 25.25}
\arr{4.5 25.25} {5.5 23.75}
\arr{8.5 25.25} {9.5 23.75}
\arr{12.5 25.25} {13.5 23.75}
\arr{16.5 25.25} {17.5 23.75}
\arr{20.5 25.25} {21.5 23.75}
\arr{24.5 25.25} {25.5 23.75}

\arr{2.5 22.25} {3.4 20.9}
\arr{6.5 22.25} {7.4 20.9}
\arr{10.5 22.25} {11.4 20.9}
\arr{14.5 22.25} {15.4 20.9}
\arr{18.5 22.25} {19.4 20.9}
\arr{22.5 22.25} {23.4 20.9}
\arr{26.5 22.25} {27.4 20.9}
\arr{0.6 20.9} {1.5 22.25}
\arr{4.6 20.9} {5.5 22.25}
\arr{8.6 20.9} {9.5 22.25}
\arr{12.6 20.9} {13.5 22.25}
\arr{16.6 20.9} {17.5 22.25}
\arr{20.6 20.9} {21.5 22.25}
\arr{24.6 20.9} {25.5 22.25}

\arr{2.7 18.05} {3.4 19.1}
\arr{6.7 18.05} {7.4 19.1}
\arr{10.7 18.05} {11.4 19.1}
\arr{14.7 18.05} {15.4 19.1}
\arr{18.7 18.05} {19.4 19.1}
\arr{22.7 18.05} {23.4 19.1}
\arr{4.6 19.1} {5.3 18.05}
\arr{8.6 19.1} {9.3 18.05}
\arr{12.6 19.1} {13.3 18.05}
\arr{16.6 19.1} {17.3 18.05}
\arr{20.6 19.1} {21.3 18.05}
\arr{24.6 19.1} {25.3 18.05}

\arr{6.7 15.95} {7.2 15.2}
\arr{10.7 15.95} {11.2 15.2}
\arr{14.7 15.95} {15.2 15.2}
\arr{18.7 15.95} {19.2 15.2}
\arr{22.7 15.95} {23.2 15.2}
\arr{4.7 15.2} {5.2 15.95}
\arr{8.7 15.2} {9.2 15.95}
\arr{12.7 15.2} {13.2 15.95}
\arr{16.7 15.2} {17.2 15.95}
\arr{20.7 15.2} {21.2 15.95}

\arr{8.7 12.95}{9.2 12.2}
\arr{12.7 12.95}{13.2 12.2}
\arr{16.7 12.95}{17.2 12.2}
\arr{20.7 12.95}{21.2 12.2}
\arr{6.8 12.2}{7.3 12.95}
\arr{10.8 12.2}{11.3 12.95}
\arr{14.8 12.2}{15.3 12.95}
\arr{18.8 12.2}{19.3 12.95}

\arr{10.8 9.95}{11.3 9.2}
\arr{14.8 9.95}{15.3 9.2}
\arr{18.8 9.95}{19.3 9.2}
\arr{8.7 9.2} {9.2 9.95}
\arr{12.7 9.2} {13.2 9.95}
\arr{16.7 9.2} {17.2 9.95}

\arr{10.5 7.25}{11 7.5}
\arr{14.5 7.25}{15 7.5}
\arr{13 7.5}{13.5 7.25}
\arr{17 7.5}{17.5 7.25}

\arr{12.8 6.4}{13.3 5.4}
\arr{16.8 6.4}{17.3 5.4}
\arr{14.7 5.4} {15.2 6.4}
\arr{10.7 5.4} {11.2 6.4}

\arr{10.5 3} {11.3 1.4}
\arr{14.5 3} {15.3 1.4}
\arr{12.7 1.4} {13.5 3}
\arr{16.7 1.4} {17.5 3}

\multiput{$\circ$} at   2 4  2 7  2 11  2 17
                        4 14  4 8  4 0  6 4  6 7  6 11  8 0  8 8  10 4  10 7
                        18 4  18 7  20 0  20 8  22 4  22 7  22 11  24 0  24 8 
                        24 14  26 4  26 7  26 11  26 17 /
\multiput{$\cdots$} at 0 4  0 8  0 11  0 17 
                        28 4  28 8  28 11  28 17  /

\setdots <2pt>
\plot 0 23  1.5 23 /
\plot 4.5 26  7.5 26 /
\plot 8.5 26  11.5 26 /
\plot 12.5 26  15.5 26 /
\plot 16.5 26  19.5 26 /
\plot 20.5 26  23.5 26 /
\plot 26.5 23  28 23 /
\plot 0 0  3.5 0 /
\plot 4.5 0  7.5 0 /
\plot 8.5 0  11.5 0 /
\plot 12.5 0  15.5 0 /
\plot 16.5 0  19.5 0 /
\plot 20.5 0  23.5 0 /
\plot 24.5 0  28 0 /

\setsolid 

\setshadegrid span <1.5mm>
\vshade   -1  0 23  <,z,,> 
           2  0 23  <z,z,,>
           4  0 26  <z,z,,>
          24  0 26  <z,z,,>
          26  0 23  <z,,,>
          28.5 0 23 /
\endpicture}

$$\diamondfoursix$$

Of particular interest is the module on the intersection of the ray
starting at $P$ and the coray ending at $I$.  It is the only indecomposable
object in $\Cal S_4(k[T]/T^6)$ in which the subgroup has two summands of
type $k[T]/T^4$ which are shifted against each other. 

We conclude with the remark that the connecting component $\Cal U$ is in fact
made up from many more and smaller diamonds.  For example, for each fraction
of type $\gamma=\frac 1n$ or of type $\gamma=\frac{n-1}n$, the tubular 
family in $\Cal S(k[T]/T^6)$ 
of index $\gamma+\Bbb Z$ intersected with $\Cal S_4(k[T]/T^6)$
consists exactly of three wings of size$(5,2,1)$.  These wing triples are
aligned in $\Cal U$ according to their index.  For example, the 
irreducible predecessor of $P$ (pictured on the left) occurs in the 
wing of sidelength five from $\Cal T_{1/3}^6$ and the irreducible
successor of $I$ (pictured on the right) occurs in the wing of length five
in $\Cal T_{2/3}^6$.

        \bigskip\bigskip
\frenchspacing
\noindent
{\bf References}        {\baselineskip=9pt \rmk
\parindent=.8truecm
                                         \medskip\smallskip\noindent
\item{[1]} D.\ M.\ Arnold: {\itk Abelian Groups and Representations of Finite 
  Partially Ordered Sets,} Springer CMS Books in Mathematics (2000).
\item{[2]} M.~Auslander and S.~O.~Smal\o: 
  {\itk Almost Split Sequences in Subcategories,}
  Journal of Algebra 69 (1981), 426-454.
\item{[3]} D.\ Happel: 
 {\itk Triangulated Categories in the 
        Representation Theory of Finite Dimensional Algebras,}
        {\rmk London Mathematical Society Lecture Notes series {\bfk 119}},
        Cambridge University Press, Cambridge, 1988.
\item{[4]} F.~Richman and E.~A.~Walker: 
        {\itk Subgroups of p$\ssize{}^5$-Bounded Groups,}
  in: Trends in Mathematics, Birk\-h\"auser 1999, 55--73.
\item{[5]} C.~M.~Ringel: 
        {\itk Tame Algebras and Integral Quadratic Forms,}
  Springer LNM 1099 (1984).
\item{[6]} C.~M.~Ringel and M.~Schmidmeier: 
        {\itk Submodules of Modules
  over a Uniserial Ring,} manuscript (2004), 1--24.
\item{[7]} C.~M.~Ringel and M.~Schmidmeier: 
        {\itk The Auslander-Reiten Translation
   in Submodule Categories,} manuscript (2004), 1--25.
\item{[8]} C.~M.~Ringel and M.~Schmidmeier: {\itk Subgroup Categories
        of Wild Representation Type,} manuscript (2004), 1--10.
\item{[9]} W.~Rump: 
        {\itk Irreduzible und Unzerlegbare Darstellungen Klassischer
  Ordnungen,} Bayreuther Mathematische Schriften {\bfk 32} (1990), 
  viii+405pp.  
\par}

\bye